\documentclass{amsart}
\usepackage{amssymb}
\usepackage{amsfonts}

\newtheorem{theorem}{Theorem}
\theoremstyle{plain}

\newtheorem{example}{Example}

\newtheorem{lemma}{Lemma}

\newtheorem{proposition}{Proposition}
\newtheorem{remark}{Remark}

\numberwithin{equation}{section}
\numberwithin{example}{section}
\numberwithin{lemma}{section}
\numberwithin{proposition}{section}
\numberwithin{remark}{section}

\begin{document}
\title[Hypergeometric functions]{Variations on hypergeometric functions}
\author{Micha\l\ Zakrzewski}
\address{Laboratory of Mathematics, Jan Kochanowski University\\
ul. Uniwersytecka 7, 25-406 Kielce, Poland}
\email{zakrzewski@mimuw.edu.pl}
\author{Henryk \.{Z}o\l \c{a}dek}
\address{Institute of Mathematics, University of Warsaw\\
ul. Banacha 2, 02-097 Warsaw, Poland}
\email{zoladek@mimuw.edu.pl}

\begin{abstract}
We prove new integral formulas for generalized hypergeometric functions and
their confluent variants. We apply them, via stationary phase formula, to
study WKB expansions of solutions: for large argument in the confluent case
and for large parameter in the general case. We also study variations of
hypergeometric functions for small perturbations of hypergeometric
equations, i.e., in expansions of solutions in powers of a small parameter.
Next, we present a new proof of a theorem due to Wasow about equivalence of
the Airy equation with its perturbation; in particular, we explain that this
result does not deal with the WKB solutions and the Stokes phenomenon.
Finally, we study hypergeometric equations, one of second order and another
of third order, which are related with two generating functions for MZVs,
one $\Delta _{2}(\lambda )$ for $\zeta (2,\ldots ,2)$'s and another $\Delta
_{3}(\lambda )$ for $\zeta (3,\ldots ,3)$'s; in particular, we correct a
statement from \cite{ZZ3} that the function $\Delta _{3}(\lambda )$ admits a
regular WKB expansion.
\end{abstract}

\maketitle

\section{Introduction}

Hypergeometric functions play dominating role in the theory of linear
differential equations. Many special functions in Mathematical Physics are
related with hypergeometric functions (see \cite{BE1, BE2, GM}). They appear
as entries in many representations of Lie groups (see \cite{Vil}) and in the
discrete groups of symmetries (see \cite{Der}). They appear also in the
Mirror Symmetry (see \cite{COGP, CK, LYau} and Example 2.2 below).\bigskip

Recently it has turned out that some generating functions for Multiple Zeta
Values (MZVs) \footnote{%
The Multiple Zeta Values (MZVs) are defined by%
\begin{equation*}
\zeta (d_{1},\ldots ,d_{k})=\sum_{0<n_{1}<\ldots <n_{k}}\frac{1}{%
n_{1}^{d_{1}}\ldots n_{k}^{d_{k}}},
\end{equation*}%
$d_{j}\geq 1,\ d_{k}\geq 2.$ They equal $\mathrm{Li}_{d_{1},\ldots
,d_{k}}\left( 1\right) ,$ where%
\begin{equation*}
\mathrm{Li}_{d_{1},\ldots ,d_{k}}(t)=\sum_{0<n_{1}<\ldots <n_{k}}\frac{%
t^{n_{k}}}{n_{1}^{d_{1}}\ldots n_{k}^{d_{k}}}
\end{equation*}%
are the polylogarithms.} are expressed via hypergeometric series (see \cite%
{AOW, CFR, Li, KoZa, OhZa, Zag, Zud}). Two such series have attracted our
attention in \cite{ZZ1, ZZ3}:%
\begin{eqnarray}
\Delta _{2}\left( \lambda \right)  &:&=1-\zeta \left( 2\right) \lambda
^{2}+\zeta \left( 2,2\right) \lambda ^{4}-\ldots ,  \label{1.1} \\
\Delta _{3}\left( \lambda \right)  &:&=1-\zeta \left( 3\right) \lambda
^{3}+\zeta \left( 3,3\right) \lambda ^{6}-\ldots .  \label{1.2}
\end{eqnarray}%
They are values at $t=1$ of corresponding hypergeometric functions, which
are solutions to hypergeometric equations: $(1-t)\mathcal{D}%
_{t}^{2}u=\lambda ^{2}tu$, of order 2,and $(1-t)\mathcal{D}_{t}^{3}u=\lambda
^{3}tu$, of order 3, respectively, where $\mathcal{D}_{t}=t\frac{\partial }{%
\partial t}$ (see Example 2.1 and Chapter 5 below).

The series (1.1) is simple,%
\begin{equation*}
\Delta _{2}(\lambda )=\frac{\sin (\pi \lambda )}{\pi \lambda }=1-\frac{\pi
^{2}}{6}\lambda ^{2}+\ldots .
\end{equation*}%
Anyway, we succeeded in \cite{ZZ1} to apply the WKB analysis to obtain two
new proofs of the above formula. Here crucial is the fact that the basic
solutions of the equation near $t=1$ are also expressed via hypergeometric
series in powers of $s=1-t.$

But the hypergeometric equation associated with $\Delta _{3}(\lambda )$ has
turned out more delicate.\footnote{%
Recall that the irrationality of $\zeta (3)$ was proved firstly by R. Ap\'{e}%
ry \cite{Ap} (see also \cite{vPo}). F. Beukers and C. Peters in \cite{BePe}
found a third order Picard--Fuchs equation for periods of some K3 surface
and associated it with the Ap\'{e}ry's recurrence for approximations of $%
\zeta (3).$} One cannot get a simple recurrence for coefficients in series
in powers of $s=1-t$ defining basic solutions. In \cite{ZZ3} we tried to
apply the stationary phase formula as $\lambda \rightarrow \infty $ to
integrals defining solutions$.$ Since near $s=0$ there are no such integral
formulas, we approximated basic solutions by Bessel type series, with
explicit integral representations. For this we needed a sort of equivalence
between a hypergeometric equation and a Bessel type equation. We have found
corresponding statement in the book of W. Wasow \cite[Chapter VIII]{Was} in
the case of a perturbation of the Airy equation (see Example 2.5 below).
Wasow directly refers to the WKB analysis.

Using such equivalence we claimed in \cite[Theorem 5.7]{ZZ3} that the
function $\Delta _{3}(\lambda )$ admits a WKB\ type representation near $%
\lambda =\infty $ (with suitable Stokes operators) and satisfies a sixth
order linear differential equation near $\lambda =\infty $ with meromorphic
coefficients.

Unfortunately, this is not true. There is no WKB type representation, at
least as claimed; see Proposition 5.1 below. The function $\Delta
_{3}(\lambda )$ satisfies a third order linear differential equation, but
with more complex singularity at $\lambda =\infty .$

Our mistake relied upon an improper interpretation of the Wasow theorem. In
Section 4 below we demonstrate that the equivalence between the Airy
equation and its perturbation is quite standard and does not involve the WKB
analysis and the Stokes phenomena. One only should properly normalize the
time.\bigskip

When working on the subject we have elaborated some properties of
hypergeometric functions and equations, which seem to be new and which
motivated this work. \medskip

We begin with a new integral formula for the generalized hypergeometric
function $_{p}F_{q}(t),$ with $p=q+1.$ In Theorem 1 (in Section 2.1) we
present a general integral formula for the function 2.1. That formula
combined a multidimensional residuum with an integral along a hypercube. It
is different from some recurrent formula, which can be found in the
Wikipedia and which is the formula (2.10) in Remark 2.4 (in Section 2.1). In
the same section, in Example 2.1, we present integral formulas, via residua,
of the hypergeometric functions associated with the generating series (1.1)
and (1.2). In Example 2.2 we present a relation of the famous mirror
symmetry from the String Theory with some generalized hypergeometric
functions.

In Theorem 2 (in Section 2.2) we present an analogous integral formula (with
a multidimensional residuum and a hypercube) for a confluent hypergeometric
function $_{p}F_{q},$ with $p>q+1.$ In Examples 2.3 -- 2.5 we discuss the
case of rather standard confluent hypergeometric equations, like the Bessel
and Airy equations. In Example 2.6 we analyze some confluent hypergeometric
equation associated with the series (1.2); some interesting integral
formulas are presented.\medskip

The formulas from Theorem 1 and Theorem 2 are suitable to get WKB expansions
of the hypergeometric functions. In Section 2.3 we compute the WKB expansion
for completely confluent hypergeometric functions $_{0}F_{q}(t)$ as $%
t\rightarrow \infty .$

It is known that the WKB expansions are related with the so-called Stokes
phenomenon. The series defining the (formal) WKB solutions are divergent in
general, but they are asymptotic for some analytic functions defined in some
sectors in $\mathbb{C}$ near $t=\infty .$ The relations of the corresponding
solutions in adjacent sectors are expressed via some constant matrices, the
so-called Stokes matrices. These Stokes matrices were computed in an
important work of C. Duval and A. Mitschi \cite{DuMi}.

In Theorem 3 (in Section 2.3) we expand the confluent hypergeometric
function $_{0}F_{q}$ in the basis of WKB solutions, with the agreement that
only the coefficients before dominating solutions are correct. The
coefficients before subdominant solutions should be determined using the
Stokes matrices, which we skip (because of complexity). In the same section,
in Example 2.7, we analyze the Stokes phenomenon in for the standard
hypergeometric function $_{1}F_{1}$

In Section 2.4 we expand the hypergeometric function $_{q+1}F_{q}$ with a
large parameter $A,A\rightarrow +\infty ,$ in the basis of other WKB
solutions $\mathrm{e}^{AS(t)}\psi (t;A)$ (like in the Schr\"{o}dinger
equation). One arrives to a Hamilton--Jacobi equation for so-called action $%
S(t)$ and to series of transport equations for so-called amplitude $\psi
(t;A)=\psi _{0}(t)+A^{-1}\psi _{1}(t)+\ldots .$

Here arises a question of the integration constants for the solutions of the
transport equations. In the case of the Schr\"{o}dinger equations this
problem is solved using the known Born--Sommerfeld quantization conditions.
In our case one introduces so-called testing WKB\ solutions and principal
WKB solutions (determined by integral formulas).

In Theorem 4 we present corresponding expansion; again, with the agreement
that only the coefficients before the dominating WKB solutions are correct.

In Example 2.9 we discuss the WKB solutions of the hypergeometric equation
associated with the series $\Delta _{3}(\lambda )$ for large $\lambda .$%
\bigskip

When considering solutions to hypergeometric equation near another singular
point $t=1$ one finds that they are not given by simple series; the
recurrence for coefficients is not standard. We approach the problem by
considering the hypergeometric equation (near $s=1-t=0)$ as a perturbation
of a confluent hypergeometric equation. So, the solutions are expanded in
powers of a small parameter $\varepsilon $ with coefficients being
variations of confluent hypergeometric functions. Those variations admit
explicit power series expansions (in powers of $s)$ and integral
representations. This is done in Section 3. In particular, in Theorem 5 we
present a general formula for these variations with a claim that those
variations admit integral representations.

In Example 3.2 we find variations of some confluent hypergeometric functions
associated with the generating series $\Delta _{3}(\lambda ).$\bigskip

Section 4 is devoted to a new (and plausibly simpler) proof of the above
mentioned Wasow theorem. In Example 4.1 we analyze variations of solutions
of some perturbation of the Airy equation.\bigskip

In Section 5 we prove the relations between some hypergeometric functions
and the generating series $\Delta _{2}(\lambda )$ and $\Delta _{3}(\lambda ).
$

In particular, in Proposition 5.1 we show that some statement from \cite[%
Theorem 5.7]{ZZ3} is definitely wrong.

We also we discuss some linear differential equations satisfied by the above
series. Plausibly new phenomena are revealed.

\section{Some new integral formulas}

\subsection{Generalized hypergeometric series}

Recall that this series equals

\begin{equation}
F\left( t\right) =_{p}F_{q}\left( \alpha _{1},\ldots ,\alpha _{p};\beta
_{1},\ldots ,\beta _{q};t\right) =\sum \frac{\left( \alpha _{1}\right)
_{n}\cdots \left( \alpha _{p}\right) _{n}}{\left( \beta _{1}\right)
_{n}\cdots \left( \beta _{q}\right) _{n}}\frac{t^{n}}{n!},  \label{2.1}
\end{equation}%
where $\left( \alpha \right) _{0}=1,$%
\begin{equation*}
\left( \alpha \right) _{n}=\alpha \left( \alpha +1\right) \cdots \left(
\alpha +n-1\right) =\frac{\Gamma \left( \alpha +n\right) }{\Gamma \left(
\alpha \right) }
\end{equation*}%
is the Pochhammer symbol and we assume that $t\in \mathbb{R}$ (for
simplicity) and $\beta _{j}\not\in \mathbb{Z}_{-}\cup 0.$ Above $p\leq q+1$;
otherwise, the series is divergent. In the case $p<q+1$ the series will be
called the confluent hypergeometric series; it is convergent in the whole
complex plane.

The series (2.1) with $p=q+1=2$ were studied by K. F. Gauss \cite{Gaus}.
Their generalization was introduced by J. Thomae \cite{Tho}.

Usually we will write $F\left( \alpha _{1},\ldots ,\alpha _{p};\beta
_{1},\ldots ,\beta _{q};t\right) $ or $F\left( t\right) $; in the case $p=0$
we will write $F\left( \emptyset ;\beta _{1},\ldots ,\beta _{q};t\right) $
and in the case $q=0$ we will write $F\left( \alpha _{1},\ldots ,\alpha
_{p};\emptyset ;t\right) .$

\begin{theorem}
\label{t1}Assume that $p=q+1$ and $\mathrm{Re}\beta _{j}>1.$ Then function
(2.1) has the following representation:%
\begin{equation}
F\left( t\right) =C\int \mathrm{d}^{q}\tau \prod_{i=1}^{q}\left( 1-\tau
_{i}\right) ^{\beta _{j}-2}\mathrm{Res}\prod\limits_{j=1}^{q+1}\left(
1-a_{j}\eta \right) ^{-\alpha _{j}}\frac{\mathrm{d}^{q+1}\mathrm{\ln }a}{%
\mathrm{d\ln }\left( a_{1}\cdots a_{q+1}\right) },  \label{2.2}
\end{equation}%
where%
\begin{equation*}
\eta =\left( \tau _{1}\cdots \tau _{q}t\right) ^{\frac{1}{q+1}},
\end{equation*}%
\begin{equation}
C=\prod \left( \beta _{i}-1\right) ,  \label{2.3}
\end{equation}%
the integral $\int $\ runs over the hypercube $\left\{ 0<\tau
_{i}<1:i=1,\ldots ,q\right\} $, the residuum is treated as integration along
the $q-$dimension fundamental cycle in the hypersurface%
\begin{equation*}
\left\{ a_{1}\cdots a_{q+1}=1\right\} \simeq \left( \mathbb{C}^{\ast
}\right) ^{q}
\end{equation*}%
of the Gelfand--Leray $q-$form $\mathrm{d}^{q+1}\mathrm{\ln }a/\mathrm{d\ln }%
\left( a_{1}\cdots a_{q+1}\right) $.\footnote{%
If $\omega $ is a holomorphic differential $k-$form in $\mathbb{C}^{k}$ and $%
\left\{ f=0\right\} \subset \mathbb{C}^{k}$ is an analytic hypersurface then
the Gelfand--Leray form $\omega /\mathrm{d}f$ in $\left\{ f=0\right\} $ is
represented by a $\left( k-1\right) -$form $\eta $ such that $\eta \wedge
\mathrm{d}f=\omega .$%
\par
We also use the notation $\mathrm{d}\ln ^{q+1}\ln a$ in place of $\left(
\mathrm{d}a_{1}/a_{1}\right) \wedge \cdots \wedge \left( \mathrm{d}%
a_{q+1}/a_{q+1}\right) .$}
\end{theorem}

\begin{remark}
\label{r21}When we parametrize the hypersurface $\left\{ a_{1}\cdots
a_{q+1}=1\right\} $ by
\begin{equation*}
a_{1}=b_{1}^{q},a_{2}=\frac{b_{2}^{q-1}}{b_{1}},a_{3}=\frac{b_{3}^{q-2}}{%
b_{1}b_{2}},\ldots ,a_{q+1}=\frac{1}{b_{1}\cdots b_{q}},
\end{equation*}%
then the corresponding form becomes%
\begin{equation*}
\frac{\mathrm{d}^{q+1}\mathrm{\ln }a}{\mathrm{d\ln }\left( a_{1}\cdots
a_{q+1}\right) }=q!\cdot \mathrm{d}^{q}\ln b.
\end{equation*}%
On the other hand, the map $b\longmapsto a$ is a covering of degree $q!$. In
application we will deal either with the variables $a_{j}$ or with the
variables $b_{k}$.
\end{remark}

\begin{remark}
\label{r22}In the case some $\mathrm{Re}\beta _{j}<1$ (but $\beta
_{j}\not\in \mathbb{Z})$ the corresponding integral $\int \mathrm{d}\tau _{j}
$ is replaced with $\left( 1-\exp \left( -2\pi \mathrm{i}\beta _{j}\right)
\right) ^{-1}\times $ integral along a loop with vertex at $\tau _{j}=0$ and
surrounding $\tau _{j}=1$ in negative direction.
\end{remark}

\begin{remark}
\label{r23}One can replace the terms $a_{j}\left( \tau _{1}\cdots \tau
_{q}t\right) ^{1/\left( q+1\right) }$ in Eq. (2.2) with%
\begin{equation}
a_{j}\left( \tau _{1}\cdots \tau _{q}\right) ^{1/\left( q+1\right) }\times
\nu _{j}t^{\mu _{j}},  \label{2.4}
\end{equation}%
where $\mu _{j}\geq 0$, $\nu _{j}\in \mathbb{C}$\ and $\sum_{j}\mu _{j}=%
\frac{1}{q+1},$ $\prod \nu _{j}=1;$ we will use such replacements when
writing integral formulas for confluent hypergeometric functions in Theorem
2 below. In the non-confluent case the choice from Eq. (2.2) seems to be the
most natural.\medskip
\end{remark}

When viewing various integral formulas for special functions (see \cite{BE1,
BE2}) one can observe that there are only few ideas behind them. \medskip

One is the residuum formula. For example, the simplest \textbf{Bessel
function}

\begin{equation}
J_{0}(t)=\sum_{n\geq 0}\frac{\left( -t^{2}/4\right) ^{n}}{\left( n!\right)
^{2}}=\sum_{n}\frac{\left( t/2\right) ^{n}}{n!}\frac{\left( -t/2\right) ^{n}%
}{n!}  \label{2.5}
\end{equation}%
can be expressed as\footnote{%
In the Complex Analysis the residuum is treated as suitable coefficient in
the Laurent expansion of a function. Here we prefer treating it as an
invariant of a suitable holomorphic form, with a natural geometric meaning.}%
\begin{equation}
J_{0}\left( t\right) =\mathrm{Res}_{b=0}\left( \sum_{m}\frac{\left(
bt/2\right) ^{m}}{m!}\right) \left( \sum_{n}\frac{\left( -t/2b\right) ^{n}}{%
n!}\right) \frac{\mathrm{d}b}{b}=\mathrm{Res}_{b=0}\mathrm{e}^{\frac{t}{2}%
\left( b-\frac{1}{b}\right) }\mathrm{d}\ln b.  \label{2.6}
\end{equation}

We have also the \textbf{Euler Gamma} and \textbf{Beta functions}%
\begin{eqnarray*}
\Gamma \left( \mu \right) &=&\int_{0}^{\infty }\tau ^{\mu -1}\mathrm{e}%
^{-\tau }\mathrm{d}\tau , \\
B\left( \mu ,\nu \right) &=&\frac{\Gamma \left( \mu \right) \Gamma \left(
\nu \right) }{\Gamma \left( \mu +\nu \right) }=\int_{0}^{1}\tau ^{\mu
-1}\left( 1-\tau \right) ^{\nu -1}\mathrm{d}\tau ;
\end{eqnarray*}%
here $\mathrm{Re}\mu ,$ $\mathrm{Re}\nu >0,$ otherwise we act like in Remark
2.2. For example, the Euler formula (2.9) below for $F\left( \alpha
_{1},\alpha _{2};\beta ;t\right) $ is essentially based on the Beta function.

Finally, recall also the \textbf{binomial formula}%
\begin{equation}
_{1}F_{0}\left( \gamma ;z\right) =F\left( \gamma ;\emptyset ;z\right)
=\sum_{n}\frac{\left( \gamma \right) _{n}}{n!}z^{n}=\left( 1-z\right)
^{-\gamma }.  \label{2.7}
\end{equation}
\medskip

\textit{Proof of Theorem 1}. We give the proof of Theorem 1 in the case $%
q=2; $ it will be clear how to generalize it for greater $q$'s. Series (2.1)
for $q=2,$ i.e.,%
\begin{equation*}
F\left( t\right) =\sum_{n}\frac{n!}{\left( \beta _{1}\right) _{n}}\frac{n!}{%
\left( \beta _{2}\right) _{n}}\frac{\left( \alpha _{1}\right) _{n}t^{n/3}}{n!%
}\frac{\left( \alpha _{2}\right) _{n}t^{n/3}}{n!}\frac{\left( \alpha
_{3}\right) _{n}t^{n/3}}{n!},
\end{equation*}%
can be expressed via the following iterated residuum formula:%
\begin{eqnarray*}
F\left( t\right) &=&\mathrm{Res}\left( \sum_{j}\frac{j!}{\left( \beta
_{1}\right) _{j}}c_{1}^{4j}\right) \left( \sum_{k}\frac{k!}{\left( \beta
_{2}\right) _{k}}\frac{c_{2}^{2k}}{c_{1}^{k}}\right) \left( \sum_{l}\frac{%
\left( \alpha _{1}\right) _{l}}{l!}\frac{b_{1}^{2l}t^{l/3}}{%
c_{1}^{l}c_{2}^{l}}\right) \\
&&\left( \sum_{m}\frac{\left( \alpha _{2}\right) _{m}}{m!}\frac{%
b_{2}^{m}t^{m/3}}{c_{1}^{m}c_{2}^{m}b_{1}^{m}}\right) \left( \sum_{n}\frac{%
\left( \alpha _{3}\right) _{n}}{n!}\frac{t^{n/3}}{%
c_{1}^{n}c_{2}^{n}b_{1}^{n}b_{2}^{n}}\right) \mathrm{d}^{2}\ln c\mathrm{d}%
^{2}\ln b.
\end{eqnarray*}%
Indeed, taking the residuum at $b_{2}=0$ means $n=m,$ taking the residuum at
$b_{1}=0$ means $m=l,$ etc.

Since $j!/\left( \beta _{1}\right) _{j}=\Gamma \left( j+1\right) \Gamma
\left( \beta _{1}\right) /\Gamma \left( \beta _{1}+j\right) =\left( \beta
_{1}-1\right) \times \Gamma \left( j+1\right) $\linebreak $\Gamma \left(
\beta _{1}-1\right) /\Gamma \left( \beta _{1}+j\right) =\left( \beta
_{1}-1\right) \cdot B\left( j+1,\beta _{1}-1\right) ,$ the first factor in
the above formula equals%
\begin{equation*}
\left( \beta _{1}-1\right) \sum_{j}\int_{0}^{1}\left( c_{1}^{4}\tau
_{1}\right) ^{j}\left( 1-\tau _{1}\right) ^{\beta _{1}-2}\mathrm{d}\tau
_{1}=\left( \beta _{1}-1\right) \int_{0}^{1}\left( 1-s\right) ^{\beta
_{1}-2}\left( 1-c_{1}^{4}\tau _{1}\right) ^{-1};
\end{equation*}%
analogously we rewrite the second factor. This together with Eq. (2.7) gives%
\begin{eqnarray*}
F &=&C\int \mathrm{d}^{2}\tau \left( 1-\tau _{1}\right) ^{\beta
_{1}-2}\left( 1-\tau _{2}\right) ^{\beta _{2}-2}\times \mathrm{Res}_{c=0}%
\frac{\mathrm{d}c_{1}\mathrm{d}c_{2}}{\left( 1-c_{1}^{4}\tau _{1}\right)
\left( c_{1}-c_{2}^{3}\tau _{2}\right) c_{2}}\times \\
&&\mathrm{Res}_{b=0}\left( 1-\frac{b_{1}^{2}t^{1/3}}{c_{1}c_{2}}\right)
^{-\alpha _{1}}\left( 1-\frac{b_{2}t^{1/3}}{c_{1}c_{2}b_{1}}\right)
^{-\alpha _{2}}\left( 1-\frac{t^{1/3}}{c_{1}c_{2}b_{1}b_{2}}\right)
^{-\alpha _{3}}\mathrm{d}^{2}\ln b,
\end{eqnarray*}%
where $C=\left( \beta _{1}-1\right) \left( \beta _{2}-1\right) $. Above the
residua at $c_{1,2}=0$\ are calculated via the integration along the circles
$\left\vert c_{1,2}\right\vert =r_{1,2}$ with $0<r_{1}<1\leq \tau
_{1}^{-1/4} $ and $0<r_{2}<r_{1}^{1/3}\tau _{2}^{-1/3}.$

The residuum at $c_{2}=0$ is replaced with the minus residuum at $%
c_{2}=\infty $ and the minus residua at $c_{2}^{\left( k\right) }=\mathrm{e}%
^{2\pi \mathrm{i}k/3}c_{1}^{1/3}\tau _{2}^{-1/3},$ $k=0,1,2.$ The residuum
at $c_{2}=\infty $ vanishes and there remain minus the other three residua.

We are left with calculation of the residuum at $c_{1}=0$ of%
\begin{equation*}
\frac{\mathrm{d}c_{1}}{3\left( 1-c_{1}^{4}\tau _{1}\right) c_{1}}\cdot
\mathrm{Res}_{b=0}\left\{ \cdot \right\} .
\end{equation*}%
Again this residuum is replaced with the minus residuum at $c_{1}=\infty $
and the minus residua at $c_{1}^{\left( l\right) }=\mathrm{e}^{2\pi \mathrm{i%
}l/4}\tau _{1}^{-1/12}\tau _{2}^{-1/3},$ $l=0,1,2,3.$ The residuum at $%
c_{2}=\infty $ vanishes and there remain minus other three residua. Note that%
\begin{equation*}
\frac{1}{c_{1}^{\left( l\right) }c_{2}^{\left( k\right) }}=\mathrm{e}^{2\pi
\mathrm{i}(k+l)/3}\left( \tau _{1}\tau _{2}\right) ^{1/3}.
\end{equation*}

This leads to the formula%
\begin{eqnarray*}
F &=&\frac{1}{12}C\int \mathrm{d}^{2}\tau \prod \left( 1-\tau _{j}\right)
^{\beta _{j}-2}\times \sum_{k,l}\mathrm{Res}_{b=0}\mathrm{d}^{2}\ln b \\
&&\left( 1-b_{1}^{2}\zeta ^{k+l}\left( \tau _{1}\tau _{2}t\right)
^{1/3}\right) ^{-\alpha _{1}}\left( 1-\frac{b_{2}}{b_{1}}\zeta ^{k+l}\left(
\tau _{1}\tau _{2}t\right) ^{1/3}\right) ^{-\alpha _{2}} \\
&&\left( 1-\frac{1}{b_{1}b_{2}}\zeta ^{k+l}\left( \tau _{1}\tau _{2}t\right)
^{1/3}\right) ^{-\alpha _{3}},
\end{eqnarray*}%
where $\zeta ^{k+l}=\mathrm{e}^{2\pi \mathrm{i}(k+l)/3}.$ But now we have a
freedom to replace $b_{m},$ $m=1,2,$ with $b_{m}e^{i\theta _{m}}.$ It is
easy that such change allows to reduce the factors $\zeta ^{k+l}.$

Finally, one could wonder what happens when the power $t^{n}$ in Eq. (2.1)
become divided into five factors instead of three factors $t^{n/3}.$ In
fact, the result remains the same. \qed

\begin{remark}
\label{r24}In the case when some of the parameters $\beta _{j}$ is a
positive integer formula (2.2) can be simplified, i.e., the number of
integrations can be reduced.

In the case of the standard Gauss hypergeometric series $F\left( \alpha
_{1},\alpha _{2};\beta ;t\right) $ formula (2.1), i.e.,%
\begin{equation}
F\left( t\right) =\left( \beta -1\right) \int_{0}^{1}\mathrm{d}\tau \left(
1-\tau \right) ^{\beta -2}\mathrm{Res}_{b=0}\left( 1-b\sqrt{\tau t}\right)
^{-\alpha _{1}}\left( 1-b^{-1}\sqrt{\tau t}\right) ^{-\alpha _{2}}\mathrm{d}%
\ln b,  \label{2.8}
\end{equation}%
for $\beta =k+1\in \mathbb{Z}_{+},$ can be simplified as follows. We have%
\begin{eqnarray*}
F\left( t\right) &=&\sum_{n\geq 0}\frac{\left( \alpha _{1}\right) _{n}t^{n/2}%
}{\left( n+k\right) !}\frac{\left( \alpha _{2}\right) _{n}t^{n/2}}{n!} \\
&=&Dt^{-k/2}\sum_{n\geq 0}\frac{\left( \alpha _{1}-k\right) _{n+k}t^{(n+k)/2}%
}{\left( n+k\right) !}\frac{\left( \alpha _{2}\right) _{n}t^{n/2}}{n!} \\
&=&Dt^{-k/2}\mathrm{Res}_{b=0}\left( \sum_{m\geq 0}\frac{\left( \alpha
_{1}-k\right) _{m}b^{m}t^{m/2}}{m!}\right) \left( \sum_{n\geq 0}\frac{\left(
\alpha _{2}\right) _{n}t^{n/2}}{n!b^{n}}\right) \frac{\mathrm{d}b}{b^{k+1}}
\\
&=&Dt^{-k/2}\mathrm{Res}_{b=0}\left( 1-b\sqrt{t}\right) ^{k-\alpha
_{1}}\left( 1-b^{-1}\sqrt{t}\right) ^{-\alpha _{2}}\frac{\mathrm{d}b}{b^{k+1}%
},
\end{eqnarray*}%
where $D=1/\left( \alpha _{1}-1\right) \cdots \left( \alpha _{1}-k\right) .$

It is easy to generate corresponding formulas in the cases of general
hypergeometric series with some $\beta _{k}$ positive integers.

Recall also the standard Euler representation%
\begin{equation}
F\left( \alpha _{1},\alpha _{2};\beta ;t\right) =\frac{\Gamma \left( \beta
\right) }{\Gamma \left( \alpha _{2}\right) \Gamma \left( \beta -\alpha
_{2}\right) }\int_{0}^{1}\tau ^{\alpha _{2}-1}\left( 1-\tau \right) ^{\beta
-\alpha _{2}-1}\left( 1-\tau t\right) ^{-\alpha _{1}}\mathrm{d}\tau .
\label{2.9}
\end{equation}

Note that many other special functions are (more or less directly) related
with the Gauss hypergeometric function: Legendre functions, complete
elliptic integrals, Chebyshev polynomials, Jacobi functions, Jacobi
polynomials, Gegenbauer polynomials.

Finally, we note that Eq. (2.9) is generalized as follows:%
\begin{eqnarray}
&&F(\alpha _{1},\ldots ,\alpha _{q+1};\beta _{1},\ldots ,\beta _{q};t)
\label{2.10} \\
&=&\frac{\Gamma \left( \beta _{q}\right) }{\Gamma \left( \alpha
_{q+1}\right) \Gamma \left( \beta _{q}-\alpha _{q+1}\right) }%
\int_{0}^{1}\tau ^{\alpha _{q+1}-1}\left( 1-\tau \right) ^{\beta _{q}-\alpha
_{q+1}-1}\widetilde{F}\left( \tau t\right) \mathrm{d}\tau   \notag
\end{eqnarray}%
where $\widetilde{F}(z)=F(\alpha _{1},\ldots ,\alpha _{q};\beta _{1},\ldots
,\beta _{q-1};z).$ This leads to another integral formula for the
generalized hypergeometric function.
\end{remark}

\begin{example}
\label{e21}The hypergeometric series $F\left( t\right) =F\left( \lambda
,-\lambda ;1;t\right) $ takes the following form:%
\begin{equation}
F\left( t\right) =\mathrm{Res}_{b=0}\left( \frac{1-\sqrt{t}b}{1-\sqrt{t}/b}%
\right) ^{\lambda }\mathrm{d}\ln b.  \label{2.11}
\end{equation}%
According to \cite{ZZ1} it associated with the generating series%
\begin{equation*}
F\left( 1\right) =\Delta _{2}\left( \lambda \right) :=1-\zeta \left(
2\right) \lambda ^{2}+\zeta \left( 2,2\right) \lambda ^{4}-\ldots
\end{equation*}%
from Eq. (1.1) for Multiple Zeta Values (see also Section 5.1 below).

The hypergeometric series $F\left( t\right) =F\left( -\lambda ,\epsilon
\lambda ,\bar{\epsilon}\lambda ;1,1;t\right) ,$%
\begin{equation*}
\epsilon =\mathrm{e}^{\pi \mathrm{i}/3},
\end{equation*}%
takes the form%
\begin{equation}
F\left( t\right) =\mathrm{Res}_{b=0}\left( \frac{1-a_{1}\sqrt[3]{t}}{\left(
1-a_{2}\sqrt[3]{t}\right) ^{\epsilon }\left( 1-a_{3}\sqrt[3]{t}\right) ^{%
\bar{\epsilon}}}\right) ^{\lambda }\frac{\mathrm{d}^{3}\ln a}{\mathrm{d}\ln
\left( a_{1}a_{2}a_{3}\right) }.  \label{2.12}
\end{equation}%
It is associated with the generating series%
\begin{equation*}
F\left( 1\right) =\Delta _{3}\left( \lambda \right) :=1-\zeta \left(
3\right) \lambda ^{3}+\zeta \left( 3,3\right) \lambda ^{6}-\ldots
\end{equation*}%
from Eq (1.2) (compare Section 5.2 below).
\end{example}

\begin{example}
\label{e22}The function $u_{0}\left( t\right) =F\left( \frac{1}{5},\frac{2}{5%
},\frac{3}{5},\frac{4}{5};1,1,1;t\right) =\sum_{n\geq 0}\frac{\left(
5n\right) !}{\left( n!\right) ^{5}}\left( 5^{-5}t\right) ^{n}$ plays
important role in the mirror symmetry. It is one of the basic solutions of
the hypergeometric equation (2.31) below (with the parameters $\alpha _{1}=%
\frac{1}{5},,$ $\alpha _{2}=\frac{2}{5},$ $\alpha _{3}=\frac{3}{5},\alpha
_{4}=\frac{4}{5},$ $\beta _{1}=\beta _{2}=\beta _{3}=1).$ Another solution
is of the form $u_{1}\left( t\right) =u_{1}\left( t\right) \ln \left(
5^{-5}t\right) +\tilde{u}_{1}\left( t\right) ,$ where $\tilde{u}_{1}\left(
t\right) $ is analytic near $t=0.$

If $z=-5^{-5}t=a^{5}$ and $q=q\left( z\right) =\exp \left(
u_{1}/u_{0}\right) ,$\ then we have the following mirror symmetry relation:%
\begin{equation*}
5+\sum_{d\geq 1}N\left( d\right) d^{3}\frac{q^{d}}{1-q^{d}}=\frac{5}{\left(
1-t\right) u_{0}^{2}}\left( \frac{q}{z}\frac{\mathrm{d}z}{\mathrm{d}q}%
\right) ^{3},
\end{equation*}%
where $N\left( d\right) $ is the number of degree $d$ rational curves in
generic quintic hypersurface $M$ in $\mathbb{P}_{\mathbb{C}}^{4}.$

On the other hand, the functions $u_{1}$ and $u_{2}$ are some periods in the
mirror dual to $M$ variety $M^{\prime }=\left\{
x_{1}^{5}+x_{2}^{5}+x_{3}^{5}+x_{4}^{5}+x_{5}^{5}+ax_{1}x_{2}x_{3}x_{4}x_{5}=0\right\} /\left(
\mathbb{Z}_{5}\right) ^{3}.$ See \cite{COGP, CK} for more informations.

Note also that in \cite{LYau} a third order equation of the form $\left(
\mathcal{D}_{t}^{3}-tP\left( \mathcal{D}_{t}\right) \right) u=0,$ treated as
Picard--Fuchs equations for a one-parameter deformations of K3 surfaces, is
used in the mirror symmetry property for K3 surfaces.
\end{example}

\subsection{Confluent hypergeometric series}

Recall that hypergeometric series (2.1) is called \textbf{confluent
hypergeometric series} if%
\begin{equation*}
p<q+1.
\end{equation*}%
\medskip

\begin{theorem}
\label{t2}In this case we have%
\begin{equation}
F\left( t\right) =C\int \mathrm{d}^{q}\tau \prod_{i=1}^{q}\left( 1-\tau
_{i}\right) ^{\beta _{i}-2}\mathrm{Res}\prod\limits_{j=1}^{p}\left(
1-a_{j}\eta \right) ^{-\alpha _{j}}\exp [\eta t^{\kappa }\left(
a_{p+1}+\ldots +a_{q+1}\right) ],  \label{2.13}
\end{equation}%
where $\eta =\eta \left( \tau \right) =\left( \tau _{1}\cdots \tau
_{q}\right) ^{1/\left( q+1\right) }$, $\kappa =\frac{1}{q+1-p}$, $C=\prod
\Gamma \left( \beta _{i}-1\right) $ and the integral $\int \mathrm{d}%
^{q}\tau $ and the residuum are the same as in Theorem 1.\medskip
\end{theorem}

\textit{Proof}. Recall that the name `confluent' comes from the fact that
those functions are obtained from the general hypergeometric function via
some limit procedure. This means that some singular points, e.g., $t=1,$
tend to infinity; thus $t=\infty $ becomes irregular singular point of the
corresponding differential equation.

We rewrite the formula from the previous section for the non-confluent
hypergeometric series as follows. We put%
\begin{eqnarray*}
&&F\left( \alpha _{1},\ldots ,\alpha _{p},A,\ldots ,A;\beta _{1}\ldots
,\beta _{q};t/A^{q-p+1}\right)  \\
&=&\sum \frac{n!}{\left( \beta _{1}\right) _{n}}\ldots \frac{n!}{\left(
\beta _{q}\right) _{n}}\frac{\left( \alpha _{1}\right) _{n}}{n!}\cdots \frac{%
\left( \alpha _{p}\right) _{n}}{n!}\frac{\left( A\right) _{n}\left(
t^{\kappa }/A\right) ^{n}}{n!}\cdots \frac{\left( A\right) _{n}\left(
t^{\kappa }/A\right) }{n!},
\end{eqnarray*}%
$\kappa =\frac{1}{q+1-p}.$ As $A\rightarrow \infty $ it tends to the
confluent hypergeometric series. By Theorem 1 its integral form is following:%
\begin{equation*}
C\int \mathrm{d}^{q}\tau \prod_{i=1}^{q}\left( 1-\tau _{i}\right) ^{\beta
_{i}-2}\mathrm{Res}\left\{ \prod\limits_{j=1}^{p}\left( 1-a_{j}\eta (\tau
\right) ^{-\alpha _{j}}\prod\limits_{k=p+1}^{q+1}\left( 1-a_{k}\eta \left(
\tau \right) \frac{t^{\kappa }}{A}\right) ^{-A}\right\} ,
\end{equation*}%
and, as $A\rightarrow \infty $, it gives the formula from the thesis of
Theorem 2. \qed

\begin{example}
\label{e23}The classical \textbf{confluent hypergeometric function}
(sometimes denoted by the Humbert symbol $\Phi \left( \alpha ,\beta
;t\right) )$ equals%
\begin{equation*}
F\left( \alpha ;\beta ;t\right) =\left( \beta -1\right) \int_{0}^{1}\mathrm{d%
}\tau \left( 1-\tau \right) ^{\beta -2}\mathrm{Res}_{b=0}\left( 1-b\sqrt{%
\tau }\right) ^{-\alpha }\mathrm{e}^{\sqrt{\tau }t/b}\mathrm{d}\ln b.
\end{equation*}%
This representation differs from the formula
\begin{equation}
F\left( \alpha ;\beta ;t\right) =\frac{\Gamma \left( \beta \right) }{\Gamma
\left( \alpha \right) \Gamma \left( \beta -\alpha \right) }\int_{0}^{1}%
\mathrm{e}^{t\tau }\tau ^{\alpha -1}\left( 1-\tau \right) ^{\beta -\alpha -1}%
\mathrm{d}\tau ,  \label{2.14}
\end{equation}%
which is obtained as a suitable limit from formula (2.9) in Remark 2.4 above
(see also \cite[Section 6.5]{BE1}).

We note also that some special functions, like Whittaker functions,
Weber--Hermite functions, Laguerre polynomials and Hermite polynomials, the
probabilistic $\mathrm{Erf}$ function, are associated with the confluent
hypergeometric function.
\end{example}

\begin{example}
\label{e24}The \textbf{Bessel function}%
\begin{equation}
J_{\nu }\left( z\right) =\sum_{n\geq 0}\frac{\left( -1\right) ^{n}\left(
z/2\right) ^{2n+\nu }}{\Gamma \left( \nu +n+1\right) n!}=\frac{\left(
z/2\right) ^{\nu }}{\Gamma \left( \nu +1\right) }\times F\left( \emptyset
;\nu +1;-z^{2}/4\right)  \label{2.15}
\end{equation}%
is expressed via the doubly confluent function%
\begin{eqnarray*}
F\left( \emptyset ;\nu +1;t\right) &=&\sum \frac{n!}{\left( \nu +1\right)
_{n}}\frac{t^{n/2}}{n!}\frac{t^{n/2}}{n!} \\
&=&\nu \int_{0}^{1}\mathrm{d}\tau \left( 1-\tau \right) ^{\nu }\mathrm{Res}%
_{b=0}\mathrm{e}^{\sqrt{\tau t}\left( b+1/b\right) }\mathrm{d}\ln b.
\end{eqnarray*}%
So we have%
\begin{equation}
J_{\nu }\left( z\right) =\frac{\left( z/2\right) ^{\nu }}{\Gamma \left( \nu
\right) }\int_{0}^{1}\mathrm{d}\tau \left( 1-\tau \right) ^{\nu }\mathrm{Res}%
_{b=0}\mathrm{e}^{\frac{z}{2}\sqrt{\tau }\left( b-1/b\right) }\mathrm{d}\ln
b.  \label{2.16}
\end{equation}%
It is different from the Schl\"{a}fli formula from \cite{BE2, ZZ1}.

Related with the Bessel functions are Hankel functions $H_{\nu }^{\left(
j\right) }$, modified Bessel functions $K_{\nu }$ and $I_{\nu },$ Kelvin
functions and Airy functions (see the next example).
\end{example}

\begin{example}
\label{e25}The \textbf{Airy equation}%
\begin{equation}
\ddot{u}=tu  \label{2.17}
\end{equation}%
is the simplest equation not solvable in generalized quadratures (see \cite%
{Zol}); this equation and its perturbations are studied in greater detail in
Section 4 below. We rewrite it as follows:%
\begin{equation*}
\left( \mathcal{D}_{t}^{2}-\mathcal{D}_{t}-t^{3}\right) u=0,
\end{equation*}%
where $\mathcal{D}_{t}=t\frac{\partial }{\partial t}=t\partial _{t}$ is the
Euler derivative. From this it follows that its basic solutions are%
\begin{eqnarray*}
u_{1}\left( t\right) &=&\sum_{n\geq 0}\frac{\left( t^{3}/9\right) ^{n}}{%
\left( 2/3\right) _{n}n!}=F\left( \emptyset ;2/3;t^{3}/9\right) , \\
u_{2}\left( t\right) &=&t\sum_{n\geq 0}\frac{\left( t^{3}/9\right) }{\left(
4/3\right) _{n}n!}=t\cdot F\left( \emptyset ;4/3;t^{3}/9\right) ,
\end{eqnarray*}%
i.e., are expressed via doubly confluent hypergeometric series; in the
literature (see \cite{Was}) they are represented via Bessel functions.

But there are other solutions defined via integral Airy functions. For
example, we have%
\begin{eqnarray}
\mathrm{Ai}\left( t\right)  &=&\frac{1}{\pi }\int_{0}^{\infty }\cos \left(
\frac{1}{3}\tau ^{3}+\tau t\right) \mathrm{d}\tau   \label{2.18} \\
&=&\frac{1}{2\pi }3^{-1/6}\Gamma \left( \frac{1}{3}\right) u_{1}\left(
t\right) -\frac{1}{2\pi }3^{1/6}\Gamma \left( \frac{2}{3}\right) u_{2}\left(
t\right) .  \notag
\end{eqnarray}%
In \cite{Was} it is underlined that also the functions $\mathrm{Ai}\left(
\mathrm{e}^{\pm 2\pi \mathrm{i}/3}t\right) $ are also solutions; $\mathrm{Ai}%
\left( t\right) $ and $\mathrm{Ai}\left( \mathrm{e}^{2\pi \mathrm{i}%
/3}t\right) $ are chosen as basic solutions.
\end{example}

\begin{example}
\label{e26}The hypergeometric function $u_{0}\left( t\right) =F\left(
-\lambda ,\epsilon \lambda ,\bar{\epsilon}\lambda ;1,1;t\right) $ from
Example 2.1 satisfies the hypergeometric equation%
\begin{equation}
\left\{ \left( 1-t\right) \mathcal{D} _{t}^{3}+\lambda ^{3}t\right\} u=0,
\label{2.19}
\end{equation}%
where $\mathcal{D} _{t}=t\frac{\partial }{\partial t}=t\partial _{t}.$

$u_{0}\left( t\right) $ is one of three basic solutions near $t=0;$ the
other basic solutions near $t=0$ are of the form $u_{1}\left( t\right)
=u_{0}\left( t\right) \ln t+\tilde{u}_{1}\left( t\right) $ and $u_{2}\left(
t\right) =\frac{1}{2}u_{0}\left( t\right) \ln ^{2}t+\tilde{u}_{1}\left(
t\right) \ln t+\tilde{u}_{2}\left( t\right) ,$ where the functions $\tilde{u}%
_{1}\left( t\right) $ and $\tilde{u}_{2}\left( t\right) $ are analytic near $%
t=0$ (see also Section 5.2 below).

Eq. (2.19) has $t=1$ as another regular singular point. With the variable%
\begin{equation*}
s=1-t
\end{equation*}%
it takes the following form:%
\begin{equation}
\left\{ \mathcal{D} _{s}\left[ (1-s)\partial _{s}\right] ^{2}-\lambda
^{3}\right\} v=0,  \label{2.20}
\end{equation}%
$\mathcal{D} _{s}=s\frac{\partial }{\partial s}=s\partial _{s}.$ We have
three independent solutions of the form%
\begin{equation}
\begin{array}{l}
v_{1}(s;\lambda )=\lambda ^{3/2}s+O\left( s^{2}\right) , \\
v_{2}(s;\lambda )=\lambda ^{3}s^{2}+O(s^{3}), \\
v_{3}(s;\lambda )=\frac{1}{4}v_{2}\ln \lambda ^{3}s^{2}+1+w\left( s;\lambda
\right) ,%
\end{array}
\label{2.21}
\end{equation}%
where $v_{1}$ $v_{2}$ and $w=O(s^{3})$ are analytic near $s=0;$ but none of
them is a hypergeometric function and does not have integral representation.

But we can try another expansion of the above solutions. To this aim we
assume $\lambda $ large, but%
\begin{equation*}
z=\lambda ^{3}s^{2}
\end{equation*}%
small. Then we get%
\begin{equation}
v_{j}\left( s;\lambda \right) =V_{j}\left( z\right) +O(\lambda ^{-3/2}),
\label{2.22}
\end{equation}%
where $V_{j}\left( z\right) $ are basic solutions to the following triply
confluent hypergeometric equation (a generalization of the Bessel equation)%
\begin{equation}
\left( \mathcal{Q}_{0}-z\right) V=0,\text{ }\mathcal{Q}_{0}=2\mathcal{D}
_{z}\left( 2\mathcal{D} _{z}-1\right) \left( 2\mathcal{D} _{z}-2\right) .
\label{2.23}
\end{equation}%
Namely%
\begin{eqnarray}
V_{1}(z) &=&\sqrt{z}\left( 1+\sum_{n=1}^{\infty }\frac{z^{n}}{\left(
2n+1\right) !(2n-1)!!}\right) =\sqrt{z}\cdot F\left( \emptyset ;\frac{3}{2},%
\frac{1}{2};\frac{z}{8}\right) ,  \label{2.24} \\
V_{2}(z) &=&2\sum_{n=1}^{\infty }\frac{z^{n}}{\left( 2n\right) !(2n-2)!!}%
=z\cdot F\left( \emptyset ;2,\frac{3}{2};\frac{z}{8}\right) ,  \label{2.25}
\\
V_{3}(z) &=&\frac{1}{4}V_{2}(z)\ln z+1+\widetilde{V}_{3}(z),  \label{2.26}
\end{eqnarray}%
\textit{where}%
\begin{equation}
\widetilde{V}_{3}=\frac{z^{2}}{16}\lim {}_{\mu \rightarrow 0}\frac{1}{\mu }%
\left\{ \frac{\tilde{F}_{1}\left( z\right) }{\left( 3-2\mu \right) \left(
2-\mu \right) }-\frac{\tilde{F}_{2}\left( z\right) }{3\left( 2+2\mu \right) }%
\right\} =-\frac{13}{64\cdot 9}z^{2}+\ldots ,  \label{2.27}
\end{equation}%
\begin{equation*}
\tilde{F}_{1}=F\left( 1;2,\frac{5}{2}-\mu ,3-\mu ;\frac{z}{8}\right) ,\text{
\ }\tilde{F}_{2}=F\left( 1;2+\mu ,\frac{5}{2},3;\frac{z}{8}\right) .
\end{equation*}

In the proof of Eq. (2.27) we take the perturbation $\mathcal{Q}_{\mu
}=8z^{-1}\left( \mathcal{D} _{z}-\mu \right) $\linebreak $\left( \mathcal{D}
_{z}-\frac{1}{2}-\mu \right) \left( \mathcal{D} _{z}-1\right) $ of the
operator $Q_{0}$ from Eq. (2.23). The corresponding equation $\left(
\mathcal{Q}_{\mu }-z\right) V$ $=0$ has solutions: $V_{1}\left( z;\mu
\right) =z^{1/2+\mu }\left( 1+\ldots \right) ,$ $V_{2}\left( z;\mu \right)
=zF\left( \emptyset ;2-\mu ,\frac{3}{2}-\mu ;\frac{z}{8}\right) =z+\frac{%
1+\left( 7/6\right) \mu }{24}z^{2}+\frac{1+\left( 47/30\right) \mu }{32\cdot
9\cdot 5}z^{3}+\ldots $ and $V_{3}\left( z;\mu \right) =$\linebreak $z^{\mu
}F\left( \emptyset ;\mu ,\frac{1}{2};\frac{z}{8}\right) =z^{\mu }\left( 1+%
\frac{1}{4\mu }z+\frac{1-\mu }{4\cdot 24\mu }z^{2}+\ldots \right) =\frac{1}{%
4\mu }V_{2}\left( z\right) +\{\frac{1}{4}V_{2}\left( z\right) \ln z+1-\frac{1%
}{96}z^{2}$\linebreak $+\ldots \}+O\left( \mu \right) $. It follows that the
third solution is $V_{3}\left( z\right) =\lim {}_{\mu \rightarrow
0}\{V_{3}\left( z;\mu \right) $\linebreak $-\frac{1}{4\mu }V_{2}\left( z,\mu
\right) \}=\frac{1}{4}V_{2}\left( z\right) \ln z+1+\widetilde{V}_{3}\left(
z\right) ,$ where $1+\widetilde{V}_{3}=\lim {}_{\mu \rightarrow 0}\{F\left(
\emptyset ;\mu ,\frac{1}{2};\frac{z}{8}\right) -\frac{z}{4\mu }F\left(
\emptyset ;2-\mu ,\frac{3}{2}-\mu ;\frac{z}{8}\right) \}=1-\frac{13}{64\cdot
9}z^{2}+\ldots .$ The latter formula is not good for its integral
representation; we take $F\left( \emptyset ;\mu ,\frac{1}{2};\frac{z}{8}%
\right) =1+\frac{z}{4\mu }+\frac{\left( z/8\right) ^{2}}{\mu \left( \mu
+1\right) \cdot 3/2}F\left( 1;\mu +2,\frac{5}{2};3;\frac{z}{8}\right) $ and $%
\frac{z}{4\mu }F\left( \emptyset ;2-\mu ,\frac{3}{2}-\mu ;\frac{z}{8}\right)
=\frac{z}{4\mu }+\frac{z^{2}/8}{4\mu \left( 2-\mu \right) \left( 3/2-\mu
\right) }F\left( 1;3-\mu ,\frac{5}{2}-\mu ,2;\frac{z}{8}\right) .$ From this
formula (2.27) follows.

We have the following integral representations:%
\begin{eqnarray}
V_{1} &=&z^{1/2}+\frac{z^{1/2}}{2}\int_{0}^{1}\frac{\mathrm{d}\tau }{\sqrt{%
1-\tau }}\mathrm{Res}_{b=0}\sinh \left( z^{1/3}b\right) \exp \left( z^{1/3}%
\frac{\tau }{2b^{2}}\right) \frac{\mathrm{d}b}{b^{4}};  \label{2.28} \\
V_{2} &=&2z^{1/3}\mathrm{Res}_{b=0}\cosh \left( z^{1/3}b\right) \exp \left(
z^{1/3}/2b^{2}\right) \frac{\mathrm{d}b}{b^{3}}\mathrm{;}  \label{2.29} \\
\widetilde{V}_{3} &=&-\frac{z^{2}}{32}\int \sqrt{\tilde{\tau}_{2}}\tilde{\tau%
}_{3}\mathrm{\ln }\left( \tilde{\tau}_{1}\tilde{\tau}_{2}\tilde{\tau}%
_{3}\right) \mathrm{Res}_{b=0}\exp \left( z^{1/3}\phi \right) \mathrm{d}%
^{2}\ln b,  \label{2.30}
\end{eqnarray}%
\textit{where}%
\begin{equation*}
\phi =\frac{1}{2}\left( \tau _{1}b_{1}^{2}+\tau _{2}b_{2}/b_{1}+\tau
_{3}/b_{1}b_{2}\right) ,
\end{equation*}%
$\tilde{\tau}_{j}=1-\tau _{j}$ \textit{and the integral is} $%
\int_{0}^{1}\int_{0}^{1}\int_{0}^{1}\mathrm{d}^{3}\tau $. The formulas for $%
V_{1}$ (in slightly different form) and $V_{2}$ come from \cite{ZZ2} and the
formula for $\widetilde{V}_{3}$ is new.

In the proof of formula (2.28) we write $z^{-1/2}V_{1}\left( z\right)
=\sum_{n\leq 0}\frac{\left( z/2\right) ^{n}}{\left( 2n+1\right) !\left(
1/2\right) _{n}}$ $=\mathrm{Res}_{b=0}\left( \sinh \left( z^{1/3}b\right)
/z^{1/3}b\right) \times F\left( 1;\frac{1}{2};z^{1/3}/2b^{2}\right) $ $%
\mathrm{d\ln }a$ \ Here we represent the confluent hypergeometric function
as $F\left( 1;\frac{1}{2};x\right) =1+2xF\left( 1;\frac{3}{2};x\right) ,$ $%
x=z^{1/3}/2b^{2},$ where $F\left( 1;\frac{3}{2};x\right) =\frac{\Gamma
\left( 3/2\right) }{\Gamma \left( 1\right) \Gamma \left( 1/2\right) }%
\int_{0}^{1}\sqrt{1-\tau }\mathrm{e}^{x\tau }\mathrm{d}\tau $ $=\frac{1}{2}%
\int_{0}^{1}\sqrt{1-\tau }\mathrm{e}^{x\tau }\mathrm{d}\tau .$

Next, $V_{2}=2z\cdot \frac{\left( z/2\right) ^{n}}{(2n+2)!n!}=2z\cdot
\mathrm{Res}_{b=0}\left\{ \left( \cosh \left( z^{1/3}b\right) -1\right)
/z^{2/3}b^{2}\right\} \exp \left( z^{1/3}/2b^{2}\right) $ $\mathrm{d\ln }b,$
where $\mathrm{Res}_{b=0}\exp \left( cb^{-2}\right) b^{-3}\mathrm{d}b=0$.
From this formula (2.29) follows.

In the proof of formula (2.30) one uses the representations $F\left( 1;\beta
_{1},\beta _{2},\beta _{3};x\right) =\prod \left( \beta _{j}-1\right) \times
\int \int \int \mathrm{d}^{3}\tau \prod \left( 1-\tau _{j}\right) ^{\beta
_{j}-2}\mathrm{Res}_{b=0}\exp \left( 2x^{1/3}\phi \right) \mathrm{d}^{2}\ln b
$ with $\beta _{1}=2,\beta _{2}=\frac{5}{2}-\mu ,$ $\beta _{3}=3-\mu $ and $%
\beta _{1}=2+\mu ,$ $\beta _{2}=\frac{5}{2}$, $\beta _{3}=3$ (compare
Theorem 2). Therefore, the expression in the braces in Eq. (2.27) equals $-%
\frac{1}{2}\int \int \int \mathrm{d}^{3}\tau \sqrt{\tilde{\tau}_{2}}\tilde{%
\tau}_{3}\times \mathrm{Res}_{b=0}\exp \left( z^{1/3}\phi \right) \lim
{}_{\mu \rightarrow 0}\frac{1}{\mu }\left\{ \tilde{\tau}_{2}^{-\mu }\tilde{%
\tau}_{3}^{-\mu }-\tilde{\tau}_{1}^{\mu }\right\} \mathrm{d}^{2}\ln b,$
where the latter limit is\linebreak $-\ln \left( \tilde{\tau}_{1}\tilde{\tau}%
_{2}\tilde{\tau}_{3}\right) $.
\end{example}

\subsection{WKB approximations for confluent hypergeometric functions}

The name WKB approximation for solutions (or WKB solutions or WKB
expansions) comes from the works of G. Wentzel \cite{Wen}, H. A. Kramers
\cite{Kra} and L. Brillouin \cite{Bri} who found approximate solutions to
the Schr\"{o}dinger equation assuming that the Planck constant $h$ is small.
One should add the name of H. Jeffereys (see \cite{Jef}) and sometimes it is
said about the WKBJ method. See also \cite{Fed} for more informations. The
hypergeometric equations with a large parameter are considered in the next
subsection.

But, besides the cases with a large parameter, like $\lambda =h^{-1},$ one
often deals with the cases with large argument (e.g., the time variable $t$)
in differential equations. Here some solutions have exponential behavior
(like the WKB solutions) and the point $t=\infty $ is an irregular
singularity of the equation. This occurs for the confluent hypergeometric
equations.

\subsubsection{Analytic/formal solutions near singular points}

The hypergeometric function $u_{0}\left( t\right) =F\left( \alpha
_{1},\ldots ,\alpha _{p}:\beta _{1},\ldots ,\beta _{q};t\right) $ satisfies
the linear differential equation (the hypergeometric equation)%
\begin{equation}
\left( \mathcal{Q}-t\mathcal{P}\right) u=0,  \label{2.31}
\end{equation}%
where%
\begin{eqnarray}
\mathcal{Q} &=&Q\left( \mathcal{D}_{t}\right) =\mathcal{D}_{t}\left(
\mathcal{D}_{t}+\beta _{1}-1\right) \cdots \left( \mathcal{D}_{t}+\beta
_{q}-1\right) ,  \label{2.32} \\
\mathcal{P} &=&P\left( \mathcal{D}_{t}\right) =\left( \mathcal{D}_{t}+\alpha
_{1}\right) \cdots \left( \mathcal{D}_{t}+\alpha _{p}\right) ,  \label{2.33}
\end{eqnarray}%
$\mathcal{D}_{t}=t\frac{\partial }{\partial t},$ are differential operators.
Indeed, assuming a solution of the form%
\begin{equation}
u\left( t\right) =t^{\gamma }\left( 1+a_{1}t+\ldots \right)  \label{2.34}
\end{equation}%
we find that%
\begin{equation*}
Q\left( \gamma \right) =0
\end{equation*}%
and the recurrent relations%
\begin{equation}
Q\left( \gamma +n+1\right) a_{n+1}=P\left( \gamma +n\right) a_{n}.
\label{2.35}
\end{equation}

From this follows the following

\begin{proposition}
\label{p21}(\textbf{Analytic solutions near} $t=0)$ The hypergeometric
equation (2.31), with $p\leq q+1,$ has solutions%
\begin{equation}
u_{j}\left( t\right) =t^{1-\beta _{j}}\times F_{j}(t),\text{ }j=0,\ldots ,q,
\label{2.36}
\end{equation}%
where $\beta _{0}=1$ and $F_{j}$ are some hypergeometric functions (with
corresponding parameters).
\end{proposition}

If $\beta _{j}$'s are pairwise different then these functions are
functionally independent and form a fundamental system of solutions to Eq.
(2.32). Otherwise, some solutions contain summands of the form $v_{k}\left(
t\right) \ln ^{m}t$ with analytic functions $v_{k}.$ These new solutions can
be obtained as limits of suitable `difference rations' of the functions $%
u_{i}\left( t\right) $ with perturbed some parameters $\beta _{l}$; so, they
also admit integral representations (see Example 2.6 above).

Also in the confluent situation, $p<q+1,$ we have the same solutions as in
Proposition 1.\bigskip

\textbf{Near }$t=1.$ Recall that a linear differential equation has a
singular point, say $t=0,$ of \textbf{regular type}, i.e., its local
solutions have power type growth as $t\rightarrow 0,$ iff it is of the form%
\begin{equation*}
\left\{ P_{0}\left( \mathcal{D}_{t}\right) +tP_{1}\left( \mathcal{D}%
_{t}\right) +\ldots \right\} u=0,
\end{equation*}%
where $P_{0}$ is a polynomial of degree equal the order of the equation, say
$q+1,$ and the series in the braces (with polynomials $P_{j}$ of degree $%
\leq q+1$) is convergent (see \cite{Koh, Was, Zol}). (Of course Eq. (2.31)
has regular singularity at $t=0.)\smallskip $

\textit{Also }$t=1$\textit{\ is a regular singular point}$.\smallskip $

Indeed, we can write $\mathcal{Q}-t\mathcal{P}=\left( 1-t\right) \mathcal{D}%
_{t}^{q+1}+\ldots =\left( 1-t\right) \left( \frac{\partial }{\partial \left(
t-1\right) }\right) +\ldots $, where \ldots\ mean lower order derivatives in
$t-1$ with analytic coefficients; this implies the regularity at $t=1.$
However, for $p=q+1>2$,\ there are no simple formulas for basic solutions
near $t=1;$ the recurrent relations are of greater length (greater than 2
like in Eq. (2.34)).\footnote{%
F. Beukers and G. Heckman in \cite{BeHe} proved that, for $p=q+1,$ there are
$q$ analytic solutions of the form $\left( t-1\right) ^{k}\left( 1+\ldots
\right) ,$ $k=0,\ldots ,q-1,$ and one analytic solution of the form $\left(
t-1\right) ^{\gamma }\left( 1+\ldots \right) .$ Using this they calculated
the monodromy group of the equation (2.31) under some additional assumptions.%
}\bigskip

\textbf{Near} $t=\infty .$ Next, for $\xi =\frac{1}{t},$ we have $\mathcal{D}%
_{t}=-\mathcal{D}_{\xi }=-\xi \frac{\partial }{\partial \xi }$ and hence $%
\frac{1}{t}\left( \mathcal{Q}-t\mathcal{P}\right) =-\{P\left( -\mathcal{D}%
_{\xi }\right) -\xi Q\left( -\mathcal{D}_{\xi }\right) \}$, it has similar
form as $\mathcal{Q}-t\mathcal{P}$. The basic solutions near $t=\infty $ are
of the form
\begin{equation}
t^{-\alpha _{k}}\times \widetilde{F}_{k}(1/t),\text{ }k=1,\ldots ,p=q+1,
\label{2.37}
\end{equation}%
where $\widetilde{F}_{k}$ are some hypergeometric \ functions; if $p=q+1,$
then they are analytic.\bigskip

Let us consider the \textbf{case}%
\begin{equation*}
p<q+1.
\end{equation*}%
Here the solutions $u_{i}\left( t\right) ,$ with generic parameters $\beta
_{k}$, are defined via series defining entire functions in the complex
plane. It follows that the only singular points are $t=0$ (regular) and $%
t=\infty $ (irregular); this holds also in the cases of non-generic $\beta
_{k}$.

Let us look for basic solutions near $t=\infty .$ Firstly consider regular
series%
\begin{equation*}
v\left( t\right) =t^{\gamma }\left( 1+a_{1}t^{-1}+\ldots \right) ,\text{ }%
t\rightarrow \infty .
\end{equation*}%
We find that $P\left( \gamma \right) =0,$ i.e.,%
\begin{equation*}
\gamma =-\alpha _{j},\text{ }j=1,\ldots ,p.
\end{equation*}%
Moreover, we have the recurrent relations $P\left( \gamma -n\right)
a_{n}=Q\left( \gamma -n+1\right) a_{n-1},$ i.e.,%
\begin{equation*}
P\left( -\gamma +n\right) a_{n}=\left( -1\right) ^{q+1-p}Q\left( -\gamma
-1+n\right) a_{n-1},
\end{equation*}%
which implies the solutions%
\begin{equation}
v_{j}\left( t\right) =t^{-\alpha _{j}}\times G_{j}\left( z\right) ,\text{ }%
z=\left( -1\right) ^{q+1-p}/t,  \label{2.38}
\end{equation}%
$j=1,\ldots ,p,$ where $G_{j}\left( z\right) =_{q+1}F_{p}\left( \tilde{\alpha%
}_{1},\ldots \tilde{\alpha}_{q+1};\tilde{\beta}_{1},\ldots ,\tilde{\beta}%
_{p};z\right) $ are suitable hypergeometric series; they are divergent,
because the number of $\tilde{\alpha}$'s is greater than $1$ plus the number
of $\tilde{\beta}$'s.

Other solutions are WKB type series. Indeed, let us postulate a solution of
the form%
\begin{equation*}
v\left( t\right) =\mathrm{e}^{ct^{\kappa }}\times \psi (t)=\mathrm{e}%
^{ct^{\kappa }}\times t^{\mu }(1+\ldots ).
\end{equation*}%
We have $\mathcal{D}_{t}v=\mathcal{D}_{t}\left( \mathrm{e}^{ct^{\kappa
}}\times \psi (t)\right) =\mathrm{e}^{ct^{\kappa }}\times \left( c\kappa
t^{\kappa }\psi +\mathcal{D}_{t}\psi \right) $ and generally%
\begin{equation*}
\mathcal{D}_{t}^{m}v=\mathrm{e}^{ct^{\kappa }}\times \left\{ \left( c\kappa
t^{\kappa }\right) ^{m}\psi +\left( c\kappa t^{\kappa }\right) ^{m-1}\left[
\frac{m\left( m-1\right) }{2}\kappa \psi +m\mathcal{D}_{t}\psi \right]
+\ldots \right\} ,
\end{equation*}%
where the dots mean lower powers of $t$. This implies (by induction)%
\begin{eqnarray*}
\mathcal{Q}v &=&\mathrm{e}^{ct^{\kappa }}\left\{ \left( c\kappa t^{\kappa
}\right) ^{q+1}+\left[ \frac{q(q+1)}{2}\kappa +\left( q+1\right) \mu +\beta %
\right] \left( c\kappa t^{\kappa }\right) ^{q}\right\} t^{\mu }(1+\ldots ),
\\
t\mathcal{P}v &=&\mathrm{e}^{ct^{\kappa }}\left\{ \left( c\kappa t^{\kappa
}\right) ^{p}t+\left[ \frac{p\left( p-1\right) }{2}\kappa +p\mu +\alpha %
\right] \left( c\kappa t^{\kappa }\right) ^{p-1}t\right\} t^{\mu }(1+\ldots
),
\end{eqnarray*}%
where%
\begin{equation*}
\alpha =\sum \alpha _{j}\text{ and }\beta =\sum \left( \beta _{i}-1\right)
\end{equation*}%
are defined by $Q\left( \mathcal{D}_{t}\right) =\mathcal{D}_{t}^{q+1}+\beta
\mathcal{D}_{t}^{q}+\ldots $ and $P\left( \mathcal{D}_{t}\right) =\mathcal{D}%
_{t}^{p}+\alpha \mathcal{D}_{t}^{p-1}+\ldots .$ By comparing terms with the
highest powers of $t$ we get: $\kappa \left( q+1\right) =\kappa p+1$ and $%
\left( c\kappa \right) ^{q+1}=\left( c\kappa \right) ^{p},$ i.e.,%
\begin{equation}
\kappa =\frac{1}{q+1-p},\text{ }c=c_{k}=\frac{1}{\kappa }\mathrm{e}^{2\pi
\mathrm{i}k\kappa }=\left( q+1-p\right) \zeta ^{k},\text{ }k=1,\ldots ,q+1-p.
\label{2.39}
\end{equation}%
Comparing the next terms we get $\frac{1}{2}q\left( q+1\right) \kappa
+\left( q+1\right) \mu +\beta =\frac{1}{2}p\left( p-1\right) \kappa +p\mu
+\alpha ,$ i.e.,%
\begin{equation}
\mu =\kappa \left( \alpha -\beta \right) -\frac{1}{2}\kappa ^{2}\left\{
q\left( q+1\right) -p\left( p-1\right) \right\} =\kappa \left( \alpha -\beta
-\frac{q}{2}-\frac{p}{2}\right) .  \label{2.40}
\end{equation}%
It means that we have $q+1-p$ formal solutions of the WKB form%
\begin{equation}
v_{k}\left( t\right) =\mathrm{e}^{c_{k}t^{\kappa }}\times t^{\mu }\cdot
H_{k}(t^{-\kappa }),\text{ }k=p+1,\ldots ,q+1.  \label{2.41}
\end{equation}%
where $H_{k}\left( t^{-\kappa }\right) =1+\ldots $ are formal series in
powers of $t^{-\kappa }.$ We can summarize the above in the following
result, which we refer to the work \cite{DuMi} of A. Duval and C. Mitschi.

\begin{proposition}
\label{p22}(\textbf{Formal solutions near} $t=\infty )$ The hypergeometric
equation (2.31) with $p<q+1$ has $p$ formal solutions of the form (2.38) and
$q+1-p$ formal solutions of the form (2.41) as $\left\vert t\right\vert
\rightarrow \infty $.
\end{proposition}

One can expect that solutions (2.41) could be obtained from the regular
solutions (2.36), i.e., for $p=q+1,$ in some limit procedure (confluence).
We did not succeed to realize this.

\subsubsection{Stokes phenomenon}

We know that for $p<q+1$ the series $G_{j}$ in Eq. (2.31) and $H_{k}$ in Eq.
(2.41) are divergent. Unfortunately, they are not easily expressed, at least
for $p=q+1>2,$ because the corresponding recurrent relations for
coefficients are of length greater than 2. There exist formulas which use
Meijer's G--functions (Mellin transformations of products of Euler Gamma
functions), see also \cite{DuMi, Me}.

Recall also that the solutions $\nu _{l}\left( t\right) $ are subject to the
so-called Stokes phenomenon. The neighborhood of $t=\infty $ is covered by a
finite collection of sectors $S$, where the series $H_{k}$ and/or $G_{j}$
are asymptotic for well defined holomorphic functions. The corresponding
analytic solutions are denoted $v_{j}^{S}$ and $v_{k}^{S}$. These analytic
solutions are defined either via reduction of the corresponding first order
2--dimensional linear systems to diagonal normal form (as it is done in \cite%
{Was, Zol}), or via Meijer's G--functions (as it is done in \cite{DuMi}) or
via Borel/Laplace transformations (see \cite{Bal, Ram}).

But when one passes to an adjacent sector $S^{\prime }$ the collection $%
\left\{ v_{1}^{S}(t),\ldots ,v_{q+1}^{S}(t)\right\} $ of solutions (which
are analytic in $S$) undergoes a linear change; to a solution $v_{l}^{S}$,
which is large (dominant), a combination of solutions, which are small
(subdominant), is added. Thus we have $v_{l}^{S^{\prime
}}(t)=v_{l}^{S}(t)+\sum c_{i}v_{i}^{S}(t)$ for $t\in S\cap S^{\prime },$
i.e., the change is defined by a constant triangularized Stokes matrix.

The Stokes matrices for the confluent hypergeometric equation were
calculated by A. Duval and C. Mitschi \cite{DuMi, Mit}.\footnote{%
The authors of \cite{DuMi, Mit} used those Stokes matrices to compute the
differential Galois group of the confluent hypergeometric equations.} The
corresponding formulas are quite complex and we cannot present them here.

\subsubsection{WKB solutions at infinity for completely confluent equation}

One of the aims of our research is to express the confluent hypergeometric
function in the bases $\left\{ v_{1}^{S}(t),\ldots ,v_{q+1}^{S}(t)\right\} $%
;  there are no such formulas in \cite{DuMi} (only the Stokes matrices). But
to do it completely, one needs to know the Stokes matrices.

In the below Theorem 3 we present a corresponding formula with concrete
coefficients $K_{l}.$ But these coefficients $K_{l}$ are correct only for
WKB solutions which are dominant in a given sector. See also the previous
subsection and Example 2.7 below. We also skip the upper index $S$ in the
basic solutions $v_{l}^{S}.$

\begin{theorem}
\label{t3}In the case $p=0,$ i.e., complete confluence, as $t\rightarrow
\infty $\ we have%
\begin{equation}
F\left( \emptyset ,\beta _{1},\ldots ,\beta _{q};t\right) \sim
\sum_{l=1}^{q+1}K_{l}v_{l}\left( t\right) ,  \label{2.42}
\end{equation}%
where $v_{l}\left( t\right) =\mathrm{e}^{c_{l}t^{\kappa }}\times t^{\mu
}\cdot H_{l}(t^{\kappa })$, $\kappa =\frac{1}{q+1},$ $c_{l}=\left(
q+1\right) \zeta ^{l},$ $\zeta =\mathrm{e}^{2\pi \mathrm{i}/\left(
q+1\right) },$ $\mu =-\frac{q/2+\beta }{q+1}$, $\beta =\sum \left( \beta
_{j}-1\right) ,$ are functions (2.41) and the constants, before dominating
functions $v_{l}\left( t\right) $, equal%
\begin{equation}
K_{l}=\frac{1}{\sqrt{q+1}}\cdot \left( \prod \Gamma \left( \beta _{j}\right)
\right) \cdot \left( 2\pi \right) ^{-q/2}\cdot \left( \zeta ^{l}\right)
^{-q/2-\beta }.  \label{2.43}
\end{equation}%
\medskip
\end{theorem}

\textit{Proof}. The integral in Eq. (2.13) from Theorem 2 is of double type.
We have the exterior integral $\int \mathrm{d}^{q}\tau $ and the interior
integral $\mathrm{Res}$.

We firstly apply the stationary phase formula to the mountain pass integral $%
\mathrm{Res}$. The phase equals%
\begin{equation*}
\phi \left( a\right) =\sum_{j}^{q+1}a_{j},
\end{equation*}%
$a=\left( a_{1},\ldots ,a_{q+1}\right) ,$ where $a_{i}$'s are subject to the
restriction%
\begin{equation*}
\varphi \left( a\right) =a_{1}\cdots a_{q+1}=1.
\end{equation*}%
The critical points of the phase are given by the equations%
\begin{equation*}
\left( \phi -\rho \varphi \right) _{a_{j}}^{\prime }=1-\rho /a_{j}=0,
\end{equation*}%
where $\rho $ is the Lagrange multiplier. It follows that the critical
points are%
\begin{equation*}
a^{\left( l\right) }=\left( \zeta ^{l},\ldots ,\zeta ^{l}\right) ,\text{ }%
l=0,\ldots ,q,
\end{equation*}%
where $\zeta =\zeta _{q+1}=\mathrm{e}^{2\pi \mathrm{i}\kappa },$ $\kappa =%
\frac{1}{q+1}$, is the corresponding root of unity. Putting%
\begin{equation*}
a_{j}=a_{j}^{\left( l\right) }\mathrm{e}^{\mathrm{i}\theta _{j}}=\zeta ^{l}%
\mathrm{e}^{\mathrm{i}\theta _{j}}
\end{equation*}%
with small $\theta _{j}$ (such that $\theta _{1}+\ldots +\theta _{q+1}=0$)
we get%
\begin{equation*}
\phi =\left( q+1\right) \zeta ^{l}-\frac{1}{2}\zeta
^{l}\sum_{j=1}^{q+1}\theta _{j}^{2}+\ldots .
\end{equation*}

The following lemma is proved in the Appendix.

\begin{lemma}
\label{l21}Consider a quadratic form $\sum_{j=1}^{q+1}\lambda _{j}\theta
_{j}^{2}$ restricted to the hypersurface $\sum_{j=1}^{q+1}\theta _{j}=0.$
Then the determinant of the corresponding symmetric matrix equals%
\begin{equation*}
e_{q}\left( \lambda \right) =\sum_{I:\left\vert I\right\vert =q}\lambda ^{I},
\end{equation*}%
where $\lambda ^{I}=\lambda _{i_{1}}\cdots \lambda _{i_{q}}$ for $I=\left\{
i_{1},\ldots ,i_{q}\right\} \subset \left\{ 1,\ldots ,q+1\right\} $; i.e.,
it is the $q^{\mathrm{th}}$ elementary symmetric polynomial in $\lambda _{j}$%
's.
\end{lemma}

From this lemma we find that the quadratic form $\theta _{1}^{2}+\ldots
+\theta _{q+1}^{2}$, i.e., with the eigenvalues $\lambda _{j}=1,$ on the
hyperplane $\theta _{1}+\ldots +\theta _{q+1}=0$ is reduced to $\sum_{1}^{q}%
\tilde{\lambda}_{j}\tilde{\theta}^{2}$ with $\prod \tilde{\lambda}_{j}=1.$

It follows that the contribution to the integral--residuum from a
neighborhood of $a^{\left( l\right) }$, via the stationary phase formula,\
equals%
\begin{eqnarray*}
&&\left( \frac{1}{2\pi \mathrm{i}}\right) ^{q}\times \exp \left[ \left(
q+1\right) \zeta ^{l}\eta t^{\kappa }\right] \times
\prod\limits_{l=1}^{q}\int \exp \left( -\frac{1}{2}\zeta ^{l}\eta t^{\kappa }%
\tilde{\lambda}_{l}\tilde{\theta}_{l}^{2}\right) \mathrm{id}\tilde{\theta}%
_{l} \\
&=&\left( \frac{1}{2\pi }\right) ^{q}\times \exp \left[ \left( q+1\right)
\zeta ^{l}\eta t^{\kappa }\right] \times \prod\limits_{l=1}^{q}\left( \frac{%
2\pi }{\zeta ^{l}\eta t^{\kappa }\tilde{\lambda}_{l}}\right) ^{1/2} \\
&=&\left( q+1\right) ^{-1/2}\left( 2\pi \right) ^{-q/2}\zeta ^{-ql/2}\left(
\tau _{1}\cdots \tau _{q}t\right) ^{-\kappa q/2}\exp \left[ \left(
q+1\right) \zeta ^{l}\left( \tau _{1}\cdots \tau _{q}t\right) ^{\kappa }%
\right] ,
\end{eqnarray*}%
where $\eta =\left( \tau _{1}\cdots \tau _{q}\right) ^{\kappa },$ $\kappa =%
\frac{1}{q+1}$.

Now we consider the exterior integration over $\tau .$ It has the form%
\begin{equation*}
\int \prod \tau _{i}^{-\kappa q/2}\tilde{\tau}_{i}^{\beta _{j}-2}\times \exp %
\left[ \Lambda \left( \tau _{1}\cdots \tau _{q}\right) ^{\kappa }\right]
\mathrm{d}^{q}\tau ,
\end{equation*}%
where $\tilde{\tau}_{i}=1-\tau _{i}$ and $\Lambda =\left( q+1\right) \zeta
^{l}t^{\kappa }$ is a large parameter (as $t$ is large). It is not an
oscillatory integral (nor a mountain pass integral), because the `phase'
takes maximal value at the point $\tau _{1}=\ldots =\tau _{q}=1$ in the
boundary of the integration domain. Nevertheless, an asymptotic analysis is
applicable.

The leading contribution comes from the neighborhood of the corner point,
where $\tilde{\tau}_{i}$ are small. We have%
\begin{equation*}
\left( \tau _{1}\cdots \tau _{q}\right) ^{\kappa }\approx 1-\kappa \left(
\tilde{\tau}_{1}+\ldots +\tilde{\tau}_{q}\right)
\end{equation*}%
and we arrive at the product of integrals%
\begin{eqnarray*}
\prod \int_{0}^{\infty }\tilde{\tau}_{i}^{\beta _{j}-2}\mathrm{e}^{-\Lambda
\kappa \tilde{\tau}_{i}}\mathrm{d}\tilde{\tau}_{i} &=&\prod \left( \Lambda
\kappa \right) ^{1-\beta _{j}}\Gamma \left( \beta _{i}-1\right)  \\
&=&\zeta ^{-l\beta }t^{-\kappa \beta }\prod \Gamma \left( \beta
_{i}-1\right) ,
\end{eqnarray*}%
where $\Lambda \kappa =\zeta ^{l}t^{\kappa }.$ This and Eq. (2.13) give
formula (2.43). \qed

\begin{remark}
\label{r25}One could ask why we do not have an analogue of Theorem 3 in the
cases of not completely confluent hypergeometric functions, i.e., with $%
0<p<q+1.$ In fact, we have tried to apply the stationary phase formula, but
without definite conclusion. We explain what happens here in Appendix 6.2.
\end{remark}

\begin{example}
\label{e27}(Standard \textbf{confluent hypergeometric function} revisited)
We have $p=q=1$ and $\kappa =\frac{1}{q+1-p}=1$. Here we can use the
representation%
\begin{equation*}
u_{0}\left( t\right) =F\left( \alpha ;\beta ;t\right) =\frac{\Gamma \left(
\beta \right) }{\Gamma \left( \alpha \right) \Gamma \left( \beta -\alpha
\right) }\int_{0}^{1}\mathrm{e}^{t\tau }\tau ^{\alpha -1}\left( 1-\tau
\right) ^{\beta -\alpha -1}\mathrm{d}\tau
\end{equation*}%
from Example 2.3. Let $\mathrm{Re}t>0$ and $\left\vert t\right\vert
\rightarrow \infty .$ Then the leading part of the integral arises from
neighborhood of $\tau =1.$ With $s=t\left( 1-\tau \right) $ we get $%
u_{0}\left( t\right) \approx \frac{\Gamma \left( \beta \right) }{\Gamma
\left( \alpha \right) \Gamma \left( \beta -\alpha \right) }\mathrm{e}%
^{t}t^{\alpha -\beta }\int_{0}^{\infty }\mathrm{e}^{-s}s^{\beta -\alpha -1}%
\mathrm{d}s\approx \frac{\Gamma \left( \beta \right) }{\Gamma \left( \alpha
\right) }\mathrm{e}^{t}t^{\alpha -\beta }\approx \frac{\Gamma \left( \beta
\right) }{\Gamma \left( \alpha \right) }v_{2}\left( t\right) $ (compare Eq.
(2.41). Let $\mathrm{Re}t<0$ and $\left\vert t\right\vert \rightarrow \infty
.$ Then the lading part of the integral arises from neighborhood of $\tau
=0. $ With $\sigma =-\tau t$ we get $u_{0}\left( t\right) \approx \frac{%
\Gamma \left( \beta \right) }{\Gamma \left( \alpha \right) \Gamma \left(
\beta -\alpha \right) }\left( -t\right) ^{-\alpha }\int_{0}^{\infty }\mathrm{%
e}^{-\sigma }\sigma ^{\alpha -1}\mathrm{d}\sigma \approx \frac{\Gamma \left(
\beta \right) }{\Gamma \left( \beta -\alpha \right) }\left( -t\right)
^{-\alpha }\approx \frac{\Gamma \left( \beta \right) }{\Gamma \left( \beta
-\alpha \right) }\mathrm{e}^{-\pi \mathrm{i}\alpha }v_{1}\left( t\right) $
(compare Eq. (2.38)); this agrees with \cite[Section 6.13.1]{BE1}.

The other basic solution near $t=0$ is $u_{1}\left( t\right) =t^{1-\beta
}F\left( \alpha -\beta +1;2-\beta ;t\right) $ (compare Eq.(2.36)). Thus%
\begin{equation*}
u_{1}=\frac{\Gamma \left( 2-\beta \right) }{\Gamma \left( \alpha -\beta
\right) \Gamma \left( 2-\alpha \right) }t^{1-\beta }\int_{0}^{1}\mathrm{e}%
^{t\tau }\tau ^{\alpha -\beta }\left( 1-\tau \right) ^{1-\alpha }\mathrm{d}%
\tau ,
\end{equation*}%
with the asymptotic $u_{1}\left( t\right) \approx \frac{\Gamma \left(
2-\beta \right) }{\Gamma \left( \alpha -\beta +1\right) }v_{2}\left(
t\right) $ as $\mathrm{Re}t>0$ and $\left\vert t\right\vert \rightarrow
\infty $ and $u_{1}\left( t\right) \approx \frac{\Gamma \left( 2-\beta
\right) }{\Gamma \left( 1-\alpha \right) }\mathrm{e}^{-\pi \mathrm{i}\alpha
}v_{1}\left( t\right) $ as $\mathrm{Re}t<0$ and $\left\vert t\right\vert
\rightarrow \infty .$

But we can get more informations about the asymptotic behavior of these
solutions using the Stokes operators; here we will follow the method from
the book \cite{He} of J. Heading. We can cover a neighborhood of $t=\infty $
by two sectors: $S^{+}$ with the upper imaginary half-line $\left\{ t_{0}=%
\mathrm{i}s:s>0\right\} $ (the Stokes line or the line of division) as its
bisectrix and with an angle $2\pi -\varepsilon $ and $S^{-}$ with the lower
imaginary half-line. These sectors intersect at two smaller sectors: $S_{\pi
}$ with the left real half-line (conjugate Stokes line) as its bisectrix and
$S_{0}$ with the right half-line (conjugate Stokes line).

In $S^{+}$ (respectively, in $S^{-})$ one has two analytic solutions $%
v_{j}^{+}\left( t\right) ,$ $j=1,2$, (respectively, $v_{j}^{-}$) with
asymptotic behavior as the formal solutions $v_{1,2}\left( t\right) $; see
Subsection   The relations between the two fundamental systems in the
sectors $S_{\pi }$ and $S_{0}$ are defined by constant triangular Stokes
matrices: we have%
\begin{eqnarray*}
v_{1}^{+} &=&v_{1}^{-}+cv_{2}^{-},\text{ }v_{2}^{+}=v_{1}^{+}\text{ in }%
S_{\pi }, \\
v_{1}^{-} &=&v_{1}^{+},\text{ }v_{2}^{-}=v_{2}^{+}+dv_{1}^{+}\text{ in }%
S_{0}.
\end{eqnarray*}%
Sometimes it is said that the above changes take place at the conjugate
Stokes lines (see \cite{He}).

Denote%
\begin{equation*}
\zeta =\mathrm{e}^{2\pi \mathrm{i}\alpha },\text{ }\nu =\mathrm{e}^{-2\pi
\mathrm{i}\beta }
\end{equation*}%
Assume that at the upper imaginary half-line we have the expansions for $%
u_{0}\left( t_{0}\right) =F_{1}\left( t_{0}\right) $ and $u_{1}\left(
t_{0}\right) =t_{0}^{1-\beta }F_{2}\left( t_{0}\right) =\mathrm{i}\nu
^{1/2}s^{1-\beta }F_{2}\left( is\right) $, $t_{0}=\mathrm{i}s$ with $s>0:$ $%
u_{0}\left( t_{0}\right) =Av_{1}^{+}\left( t_{0}\right) +Av_{2}^{+}\left(
t_{0}\right) $ and $u_{1}\left( t_{0}\right) =Cv_{1}^{+}\left( t_{0}\right)
+Dv_{2}^{+}\left( t_{0}\right) .$ Thus%
\begin{eqnarray*}
F_{1}\left( t_{0}\right) &=&At_{0}^{-\alpha }G^{+}\left( t_{0}\right)
+Bt_{0}^{\alpha -\beta }\mathrm{e}^{t_{0}}H^{+}\left( t_{0}\right) , \\
t_{0}^{1-\beta }F_{2}\left( t_{0}\right) &=&Ct_{0}^{-\alpha }G^{+}\left(
t_{0}\right) +Dt_{0}^{\alpha -\beta }\mathrm{e}^{t_{0}}H^{+}\left(
t_{0}\right) ,
\end{eqnarray*}%
where $G^{+}$ and $H^{+}$ are analytic functions with formal power series in
$t_{0}^{-1}$ like $G_{1}$ and $H_{1}$ in Eqs. (2.38) and (2.41).

Let us rotate $t_{0}$ to $\mathrm{e}^{\pi \mathrm{i}}t_{0}=-t_{0}.$ In the
sector $S_{\pi }$ we have $t=\mathrm{e}^{\mathrm{i}\theta }t_{0}$ and $%
u_{1}\left( t\right) =\mathrm{e}^{\mathrm{i}\left( 1-\beta \right) \theta
}t_{0}^{1-\beta }F_{2}\left( t\right) ,$ $v_{1}^{\pm }\left( t\right) \sim
\mathrm{e}^{-\mathrm{i}\alpha \theta }t_{0}^{-\alpha }G_{1}\left( t\right) ,$
$v_{2}^{\pm }\sim \mathrm{e}^{\mathrm{i}\left( \alpha -\beta \right) \theta
}t_{0}^{\alpha -\beta }\mathrm{e}^{t}H_{1}\left( t\right) $ and we get%
\begin{eqnarray*}
F_{1}\left( t\right) &=&A\left( v_{1}^{-}+cv_{2}^{-}\right) +Bv_{2}^{-} \\
&\sim &A\mathrm{e}^{-\mathrm{i}\alpha \theta }t_{0}^{-\alpha }G_{1}\left(
t\right) +\left( B+cA\right) \mathrm{e}^{\mathrm{i}\left( \alpha -\beta
\right) \theta }t_{0}^{\alpha -\beta }\mathrm{e}^{t}H_{1}\left( t\right) , \\
\mathrm{e}^{\mathrm{i}\left( 1-\beta \right) \theta }t_{0}^{1-\beta
}F_{2}\left( t\right) &=&C\left( v_{1}^{-}+cv_{2}^{-}\right) +Dv_{2}^{-} \\
&\sim &C\mathrm{e}^{-\mathrm{i}\alpha \theta }t_{0}^{-\alpha }G_{1}\left(
t\right) +\left( D+cC\right) \mathrm{e}^{\mathrm{i}\left( \alpha -\beta
\right) \theta }t_{0}^{\alpha -\beta }\mathrm{e}^{t}H_{1}\left( t\right) .
\end{eqnarray*}%
In particular, for $t=-t_{0}$ in the lower Stokes line we have%
\begin{eqnarray*}
F_{1}\left( t\right) &\sim &A\zeta ^{-1/2}t_{0}^{-\alpha }G_{1}\left(
t\right) +\left( B+cA\right) \left( \zeta \nu \right) ^{1/2}t_{0}^{\alpha
-\beta }\mathrm{e}^{t}H_{1}\left( t\right) , \\
-v^{1/2}t_{0}^{1-\beta }F_{2}\left( t\right) &\sim &C\zeta
^{-1/2}t_{0}^{-\alpha }G_{1}\left( t\right) +\left( D+cC\right) \left( \zeta
\nu \right) ^{1/2}t_{0}^{\alpha -\beta }\mathrm{e}^{t}H_{1}\left( t\right) .
\end{eqnarray*}%
Finally, we rotate $-t_{0}$ to $\mathrm{e}^{\pi \mathrm{i}}(-t_{0})=t_{0}.$
The analogues of the latter formulas are following:%
\begin{eqnarray*}
F_{1}\left( t_{0}\right) &=&\left( A+d\left( B+cA\right) \right) \zeta
^{-1}t_{0}^{-\alpha }G^{+}\left( t_{0}\right) +\left( B+cA\right) \zeta \nu
t_{0}^{\alpha -\beta }\mathrm{e}^{t_{0}}H^{+}\left( t_{0}\right) , \\
vt_{0}^{1-\beta }F_{2}\left( t_{0}\right) &=&\left( C+d\left( D+cC\right)
\right) \zeta ^{-1}t_{0}^{-\alpha }G^{+}\left( t_{0}\right) +\left(
D+cC\right) \zeta \nu t_{0}^{\alpha -\beta }\mathrm{e}^{t_{0}}H^{+}\left(
t_{0}\right) .
\end{eqnarray*}%
Therefore, we get the relations%
\begin{eqnarray*}
A &=&\left( A+d\left( B+cA\right) \right) \zeta ^{-1},\text{ }B=\left(
B+cA\right) \nu \zeta , \\
\nu C &=&\left( C+d\left( D+cC\right) \right) \zeta ^{-1},\text{ }\nu
D=\left( D+cC\right) \nu \zeta .
\end{eqnarray*}%
The vanishing of the determinants of each of the above two linear systems
leads to the condition%
\begin{equation*}
cd=\left( \zeta -1\right) \left( 1-\nu \zeta \right) ;
\end{equation*}%
so, we still do not know the constants $c$ and $d$ defining the Stokes
matrices. But we know the values of the constants: $A=\frac{\Gamma \left(
\beta \right) }{\Gamma \left( \beta -\alpha \right) }\zeta ^{-1/2},$ $B=%
\frac{\Gamma \left( \beta \right) }{\Gamma \left( \alpha \right) },$ $C=%
\frac{\Gamma \left( \beta \right) }{\Gamma \left( \beta -\alpha \right) }%
\zeta ^{-1/2}$ and $D=\frac{\Gamma \left( 2-\beta \right) }{\Gamma \left(
\alpha -\beta +1\right) }.$ Using this we find%
\begin{eqnarray*}
c &=&\frac{1-\nu \zeta }{\nu \zeta }\frac{B}{A}=\frac{\Gamma \left( \beta
-\alpha \right) }{\Gamma \left( \alpha \right) }\left( 1-\nu \zeta \right)
\zeta ^{-1/2}\nu ^{-1}, \\
d &=&\frac{\Gamma \left( \alpha \right) }{\Gamma \left( \beta -\alpha
\right) }\left( \zeta -1\right) \zeta ^{1/2}\nu .
\end{eqnarray*}%
Moreover, at the upper Stokes line we have%
\begin{equation*}
F\left( \alpha ;\beta ;\mathrm{i}s\right) =\frac{\Gamma \left( \beta \right)
}{\Gamma \left( \beta -\alpha \right) }\zeta ^{-1/4}s^{-\alpha }G^{+}\left(
\mathrm{i}s\right) +\frac{\Gamma \left( \beta \right) }{\Gamma \left( \alpha
\right) }\left( \nu \zeta \right) ^{1/4}s^{\alpha -\beta }e^{\mathrm{i}%
s}H^{+}\left( \mathrm{i}s\right) \text{ }(s>0)
\end{equation*}%
and on the lower Stokes line we have%
\begin{equation*}
\begin{array}{ll}
F\left( \alpha ;\beta ;-\mathrm{i}s\right) & =\frac{\Gamma \left( \beta
\right) }{\Gamma \left( \beta -\alpha \right) }\zeta ^{-3/4}s^{-\alpha
}G^{-}\left( -\mathrm{i}s\right) \\
& +\frac{\Gamma \left( \beta \right) }{\Gamma \left( \alpha \right) }\left[
1+\left( 1-\nu \zeta \right) \zeta ^{-1/2}\nu ^{-1}\right] \left( \nu \zeta
\right) ^{3/4}s^{\alpha -\beta }e^{-\mathrm{i}s}H^{-}\left( -\mathrm{i}%
s\right) \text{ }(\text{ }s>0).%
\end{array}%
\end{equation*}

Finally, we would like to note that we did not found such analysis in the
literature.
\end{example}

\begin{example}
\label{e28}(Example 2.6 revisited) Calculations like in the proof of Theorem
3 give the following WKB expansions for the basic solutions $V_{j}:$%
\begin{eqnarray*}
V_{1} &\sim &\frac{2}{\sqrt{3}}z^{1/6}\left( -\mathrm{e}^{\frac{3}{2}%
z^{1/3}}+\bar{\epsilon}\mathrm{e}^{-\frac{3}{2}\epsilon z^{1/3}}+\epsilon
\mathrm{e}^{-\frac{3}{2}\bar{\epsilon}z^{1/3}}\right) , \\
V_{2} &\sim &\frac{2}{\sqrt{6\pi }}z^{1/6}\left( \mathrm{e}^{\frac{3}{2}%
z^{1/3}}+\bar{\epsilon}\mathrm{e}^{-\frac{3}{2}\epsilon z^{1/3}}+\epsilon
\mathrm{e}^{-\frac{3}{2}\bar{\epsilon}z^{1/3}}\right) , \\
V_{3} &\sim &\frac{\mathrm{i}}{2}\sqrt{\frac{2\pi }{3}}z^{1/6}\left( \bar{%
\epsilon}\mathrm{e}^{-\frac{3}{2}\epsilon z^{1/3}}-\epsilon \mathrm{e}^{-%
\frac{3}{2}\bar{\epsilon}z^{1/3}}\right) \text{ }\mathrm{mod}\text{ }V_{2},
\end{eqnarray*}%
where $\epsilon =\mathrm{e}^{\pi \mathrm{i}/3}.$

In \cite{ZZ3} also the Stokes matrices were found, but in a different way
than in \cite{DuMi}. We have combined a Heading method with a perturbation
of the confluent equation to a non-confluent one (where some Stokes matrix
becomes a monodromy matrix).
\end{example}

\subsection{WKB expansions of hypergeometric functions with a large parameter%
}

Here we consider the hypergeometric function (2.1) and the hypergeometric
equation (2.31) in the case when $p=q+1$ (i.e., no confluence) and
\begin{equation}
\alpha _{1}=A\nu _{1},\ldots ,\alpha _{q+1}=A\nu _{q+1},  \label{2.44}
\end{equation}%
where $A>0$ is a large parameter,%
\begin{equation}
A\rightarrow +\infty ,  \label{2.45}
\end{equation}%
and $\nu _{i},\beta _{j}$ are fixed. (One could consider the case when only $%
r<p=q+1$ of the parameters $\alpha _{i}$ are large, but this leads to a
non-conclusive analysis, like in Remark 2.5 above).

\subsubsection{Hamilton--Jacobi and transport equations}

We look for solutions in the form of \textbf{WKB series}%
\begin{equation}
U\left( t\right) =\mathrm{e}^{AS\left( t\right) }\psi \left( t;A\right) ,
\label{2.46}
\end{equation}%
where $S\left( t\right) $ is an `\textbf{action}' and%
\begin{equation}
\psi \left( t;A\right) =\psi _{0}\left( t\right) +\psi _{1}\left( t\right)
A^{-1}+\ldots   \label{2.47}
\end{equation}%
is an `\textbf{amplitude}'. As a rule the series defining $\psi $ is
divergent.

In this subsection we determine the action $S$ and the leading amplitude
term $\psi _{0}$ directly from Eq. (2.31). (In Subsection 2.4.5 we do the
same by applying the stationary phase formula to the integral in Eq. (2.2),
which is a mountain pass integral.)

We have%
\begin{eqnarray*}
\mathcal{D}U &=&\mathrm{e}^{AS\left( t\right) }\left\{ A\cdot \mathcal{D}%
S\cdot \psi +\mathcal{D}\psi \right\} , \\
\left( \mathcal{D}+A\nu \right) U &=&\mathrm{e}^{AS\left( t\right) }\left\{
A\cdot \left( \mathcal{D}S+\nu \right) \cdot \psi +\mathcal{D}\psi \right\} ,
\end{eqnarray*}%
where $\mathcal{D}=\mathcal{D}_{t}=t\frac{\partial }{\partial t}.$ Recall
that $Q\left( \mathcal{D}_{t}\right) =\mathcal{D}_{t}^{q+1}+\beta \mathcal{D}%
_{t}^{q}+\ldots $, where%
\begin{equation*}
\beta =\sum \left( \beta _{j}-1\right) ,
\end{equation*}%
By induction one proves that%
\begin{eqnarray*}
\mathcal{Q}U &=&\mathrm{e}^{AS\left( t\right) }\times \{A^{q+1}\cdot \left(
\mathcal{D}S\right) ^{q+1}\cdot \psi +A^{q}\cdot  \\
&&[\left( q+1\right) \left( \mathcal{D}S\right) ^{q}\mathcal{D}\psi +\frac{1%
}{2}q\left( q+1\right) \left( \mathcal{D}S\right) ^{q-1}\left( \mathcal{D}%
^{2}S\right) \psi +\beta \left( \mathcal{D}S\right) ^{q}\psi ]+\ldots \}
\end{eqnarray*}%
and%
\begin{eqnarray*}
\mathcal{P}U &=&\mathrm{e}^{AS\left( t\right) }\times \{A^{q+1}\cdot P\left(
\mathcal{D}S\right) \cdot \psi  \\
&&+A^{q}[\frac{1}{2}P^{\prime \prime }\left( \mathcal{D}S\right) \cdot
\mathcal{D}^{2}S\cdot \psi +P^{\prime }\left( \mathcal{D}S\right) \cdot
\mathcal{D}\psi ]+\ldots \}.
\end{eqnarray*}%
By comparing the terms before $A^{q+1}$ and $A^{q}$ in the equation $\left(
\mathcal{Q}-t\mathcal{P}\right) U=0,$ we get the \textbf{Hamilton--Jacobi
equation}%
\begin{equation}
\left( \mathcal{D}S\right) ^{q+1}-tp\left( \mathcal{D}S\right) =0,
\label{2.48}
\end{equation}%
where%
\begin{equation*}
p\left( x\right) =\prod \left( x+\nu _{j}\right)
\end{equation*}%
is a renormalized polynomial, $P\left( Ax\right) =A^{q+1}p\left( x\right) $,
and the \textbf{transport equation}%
\begin{eqnarray}
&&\left\{ \left( q+1\right) \left( \mathcal{D}S\right) ^{q}-tp^{\prime
}\left( \mathcal{D}S\right) \right\} \cdot \mathcal{D}\psi _{0}  \label{2.49}
\\
&&+\left\{ \left[ \frac{1}{2}q\left( q+1\right) \left( \mathcal{D}S\right)
^{q-1}-\frac{t}{2}p^{\prime \prime }\left( \mathcal{D}S\right) \right] \cdot
\left( \mathcal{D}^{2}S\right) +\beta \left( \mathcal{D}S\right)
^{q}\right\} \cdot \psi _{0}.  \notag
\end{eqnarray}

The Hamilton--Jacobi equation is solved in two steps. Firstly, it is an
algebraic equation for $\mathcal{D}S$ with $q+1$ solutions
\begin{equation*}
\mathcal{D}S=R^{\left( l\right) }\left( t\right) ,\text{ }l=1,\ldots ,q+1,
\end{equation*}%
where $R^{\left( l\right) }\left( t\right) =D\zeta ^{l}t^{\kappa }+\ldots $
as $t\rightarrow 0$ (here and below the upper index $(l)$ does not mean the $%
l^{\mathrm{th}}$ derivative), $\kappa =\frac{1}{q+1}$, $D=\left( \prod \nu
_{j}\right) ^{\kappa }$, $\zeta =\mathrm{e}^{2\pi \mathrm{i}\kappa })$,
which are algebraic functions. Next, we solve the initial value problems%
\begin{equation*}
t\frac{\mathrm{d}S}{\mathrm{d}t}=R^{\left( l\right) }\left( t\right) ,\text{
}S\left( 0\right) =0,
\end{equation*}%
with the solutions%
\begin{equation}
S^{\left( l\right) }\left( t\right) =\int_{0}^{t}s^{-1}R^{\left( l\right)
}\left( s\right) \mathrm{d}s,  \label{2.50}
\end{equation}%
where $S^{\left( l\right) }\left( t\right) =\kappa ^{-1}D\zeta ^{l}t^{\kappa
}+\ldots $ as $t\rightarrow 0$. The functions $S^{\left( l\right) }$ belong
to the class of functions expressed via generalized quadratures (see \cite%
{Zol}).

Knowing these solutions we solve the transport equations (but not just an
initial value problem)%
\begin{equation}
\psi _{0}^{\left( l\right) }\left( t\right) =\exp \int^{t}T^{\left( l\right)
}\left( s\right) \mathrm{d}s,  \label{2.51}
\end{equation}%
where $T^{\left( j\right) }\left( t\right) $ are suitable functions
expressed in generalized quadratures following from Eq. (2.49). The
functions $\psi _{0}^{\left( l\right) }\left( t\right) $ are also expressed
via generalized quadratures. As $t\rightarrow 0$ we have $T^{\left( l\right)
}\left( t\right) \approx -\frac{1}{2}q\left( q+1\right) \left( \mathcal{D}%
S^{\left( l\right) }\right) ^{q-1}\left( \mathcal{D}^{2}S^{\left( l\right)
}\right) /\left( q+1\right) \left( \mathcal{D}S^{\left( l\right) }\right)
^{q}$\linebreak $=-\frac{q}{2}\cdot \mathcal{D}^{2}S^{\left( l\right) }/%
\mathcal{D}S^{\left( l\right) };$ thus $\psi _{0}^{\left( l\right) }\left(
t\right) \approx \mathrm{const}\cdot \left( \mathcal{D}S^{\left( l\right)
}\right) ^{-q/2}\approx \mathrm{const}\cdot t^{-q/2\left( q+1\right) }$.

\subsubsection{Testing and principal solutions}

The terms $\psi _{j}\left( t\right) $ in the expansion $\psi \left(
t;A\right) =\psi _{0}\left( t\right) +\psi _{1}\left( t\right) A^{-1}+\ldots
$ also satisfy suitable transport type equations. They are linear
non-homogeneous equations $\dot{\psi}_{j}=a\left( t\right) \psi _{j}+b\left(
t\right) $, where the corresponding homogeneous equations are the same as in
Eq. (2.49), i.e., $a\left( t\right) =a^{\left( l\right) }\left( t\right) $
do not depend on $j,$ and $b\left( t\right) =b_{j}^{\left( l\right) }\left(
t\right) $ are expressed via the $a_{i}^{\left( l\right) }$'s for $i<j.$ So,
the solutions $\psi _{j}$ to those non-homogeneous equations are defined
modulo $\mathrm{const}\cdot \psi _{0}\left( t\right) ;$ in fact, this
non-uniqueness holds only for $j$'s from some arithmetic sequence (see
Example 2.9 below).

So, it is natural to define the \textbf{testing WKB solutions} $U_{\mathrm{%
test}}^{\left( l\right) }\left( t\right) $ as those formal WKB series (2.46)
whose functional coefficients $\psi _{j}\left( t\right) $ do not contain
terms proportional to $\psi _{0}\left( t\right) .$ The so-called \textbf{%
principal WKB solution}s $U_{\mathrm{princ}}^{\left( l\right) }\left(
t\right) $ are those defined by application of the stationary phase formula
to integral (2.2) (see below).

Recall also that a non-uniqueness of this sort was noted in the WKB approach
to the Schr\"{o}dinger equation. In order to remove this ambiguity one
imposes so-called Born--Sommerfeld quantization conditions (because the wave
functions must decay at $\pm \infty $ see \cite{Sch}).

Finally, we note that in \cite{ZZ2} we considered WKB solutions in some
cases when the corresponding Hamilton--Jacobi equation (see Eq. (2.48)) has
a solution of multiplicity greater than one.

\subsubsection{Stokes phenomenon}

Like in the case of a confluent equation the formal solutions (2.46) are
divergent. But they are asymptotic series for some analytic solutions in
suitable domains.

In the beginning of this section we imposed the assumptions (2.44)--(2.45).
Thus the parameter $A$ is \textbf{real} and large. But the time $t$ can take
complex values, $t\in \mathbb{C}$. The analytic solutions are defines in
some domains in the complex plane; we denote these domains also by $S$
(although they are not sectors). The corresponding result was proved by G.
D. Birkhoff \cite{Bir}. \footnote{%
In \cite{ZZ2} the situation whrre both, $t$ (the time) and $A$ (the
parameter), can be complex. Here we consider the simpler case.}

Anyway we divide the analytic solutions into dominant ones and subdominant
ones and we have corresponding Stokes matrices, like in Subsection 2.3.2.
The dominant solutions are those for which the exponent%
\begin{equation*}
\exp AS(t)
\end{equation*}%
in Eq. (2.46) tends to infinity as $A\rightarrow +\infty ,$ i.e., when $%
\mathrm{Re}S(t)>0$ for $t\in S;$ in fact, these solutions are analytically
continued to larger domains. The Stokes operators in some cases of confluent
hypergeometric equations were analyzed in the book \cite{He} of J. Heading;
it is desirable to study them in general situations (an analogue to the
Duval--Mitschi analysis).

In the below Theorem 4 we express the hypergeometric function in some
`basis' of principal WKB solutions, but the coefficients before the these
solutions are correct only in the dominant case.

\subsubsection{WKB expansions}

Let us apply the stationary phase formula to the integral in Eq. (2.2), with
assumption (2.44) for large $A>0,$ i.e., to%
\begin{equation}
F\left( t\right) =C\int \prod \left( 1-\tau _{i}\right) ^{\beta _{i}-2}%
\mathrm{Res}\left\{ \prod \left( 1-\eta a_{j}\right) ^{-\nu _{j}}\right\}
^{A},  \label{2.52}
\end{equation}%
where%
\begin{equation*}
\eta =\left( \tau _{1}\cdots \tau _{q}t\right) ^{\kappa },
\end{equation*}%
$\kappa =\frac{1}{q+1}.$

(a) The $\mathrm{Res}$ is the residuum of the form $\left\{ \cdot \right\}
^{A}\mathrm{d}\ln ^{q+1}a/\mathrm{d}\ln \left( a_{1}\cdots a_{q+1}\right) $;
it is a mountain pass type integral with the \textbf{phase}%
\begin{equation}
\phi =-\sum \nu _{j}\ln \left( 1-\eta a_{j}\right)   \label{2.53}
\end{equation}%
with the restriction%
\begin{equation*}
\varphi :=a_{1}\cdots a_{q+1}=1.
\end{equation*}%
Let $\rho $ be the \textbf{Lagrange multiplier}; then the critical points of
the phase are given by the equations $\partial \left( \phi -\rho \varphi
\right) /\partial a_{j}=0,$ i.e., we have%
\begin{equation}
a_{j}=\frac{\rho }{\eta \left( \rho +\nu _{j}\right) }.  \label{2.54}
\end{equation}%
Moreover, the condition $\varphi =1$ leads to the algebraic equation for $%
\rho :$%
\begin{equation}
\rho ^{q+1}=\eta ^{q+1}p\left( \rho \right) ,.  \label{2.55}
\end{equation}%
where $p\left( \rho \right) =\prod \left( \rho +\nu _{j}\right) $; note the
similarity of this equation with Eq. (2.48) above.

Let
\begin{equation}
\rho =\rho ^{\left( l\right) }\left( \eta \right) ,\text{ }l=1,\ldots ,q+1,
\label{2.56}
\end{equation}%
be its solutions; thus%
\begin{equation*}
a_{j}^{\left( l\right) }=\frac{\rho ^{\left( l\right) }}{\eta \left( \rho
^{\left( l\right) }+\nu _{j}\right) }.
\end{equation*}%
As $\eta \rightarrow 0$ we have $\rho ^{\left( l\right) }\approx \eta D\zeta
^{l},$ $a_{j}^{\left( l\right) }\approx Dv_{j}^{-1}\zeta ^{l},$ where $%
D=\left( \prod \nu _{j}\right) ^{\kappa },$ $\zeta =\mathrm{e}^{2\pi \mathrm{%
i}\kappa },$ $\kappa =\frac{1}{q+1}.$\bigskip

(b) We put $a_{j}=a_{j}^{\left( l\right) }\mathrm{e}^{\mathrm{i}\theta _{j}},
$ with the restriction $\sum \theta _{j}=1$, and expand the phase. We get%
\begin{equation}
\phi =\phi ^{\left( l\right) }-\frac{1}{2}\sum_{j}\rho ^{\left( l\right)
}\left( 1+\frac{\rho ^{\left( l\right) }}{\nu _{j}}\right) \theta
_{j}^{2}+\ldots .  \label{2.57}
\end{equation}%
where%
\begin{equation}
\phi ^{\left( l\right) }=\phi ^{\left( l\right) }\left( \eta \right) =\sum
\nu _{j}\ln \left( 1+\rho ^{\left( l\right) }/\nu _{j}\right)   \label{2.58}
\end{equation}%
are the corresponding critical values. Of course, the linear term is absent.

Consider the quadratic part, i.e., the quadratic form%
\begin{equation*}
-\frac{1}{2}\sum_{1}^{q+1}\lambda _{j}\theta _{j}^{2},\text{ }\sum \theta
_{j}=0,
\end{equation*}%
with the eigenvalues $\lambda _{j}=\rho ^{\left( l\right) }\left( 1+\rho
^{\left( l\right) }/\nu _{j}\right) .$ By Lemma 2.1 above this form is
equivalent to $-\frac{1}{2}\sum_{1}^{q}\tilde{\lambda}_{j}\tilde{\theta}%
_{j}^{2},$ where
\begin{equation*}
\mathrm{Det}=\mathrm{Det}\left( \rho ^{\left( l\right) }\left( \eta \right)
\right) =\prod \tilde{\lambda}_{j}=e_{q}\left( \lambda _{1},\ldots ,\lambda
_{q+1}\right)
\end{equation*}%
(elementary symmetric polynomial in $\lambda ).$ Let us evaluate this
quantity.

Recall that it is the coefficient before $\lambda $ in the polynomial $%
r\left( \lambda \right) =\prod \left( \lambda +\lambda _{j}\right) ,$ i.e.,
in%
\begin{eqnarray*}
\prod \left( \lambda +\rho +\rho ^{2}/\nu _{j}\right)  &=&\prod \frac{\rho
+\lambda }{\nu _{j}}\times \prod \left( \frac{\rho ^{2}}{\rho +\lambda }+\nu
_{j}\right)  \\
&=&\prod \frac{\rho }{\nu _{j}}\times \left( 1+\frac{\lambda }{\rho }\right)
^{q+1}\times p\left( \rho \left( 1+\frac{\lambda }{\rho }\right)
^{-1}\right) ,
\end{eqnarray*}%
where $\rho =\rho ^{\left( l\right) }.$ We get%
\begin{equation*}
r\left( \lambda \right) \approx \frac{\rho ^{q+1}}{p\left( 0\right) }\times
\left( 1+\left( q+1\right) \frac{\lambda }{\rho }\right) \times p\left( \rho
\right) \left\{ 1-\frac{p^{\prime }\left( \rho \right) }{p\left( \rho
\right) }\lambda \right\}
\end{equation*}%
and hence%
\begin{equation}
\mathrm{Det}^{\left( l\right) }\left( \eta \right) =\rho ^{q+1}\frac{p\left(
\rho \right) }{p\left( 0\right) }\left\{ \frac{q+1}{\rho }-\frac{p^{\prime
}\left( \rho \right) }{p\left( \rho \right) }\right\} ,\text{ }\rho =\rho
^{\left( l\right) }\left( \eta \right) .  \label{2.59}
\end{equation}%
It follows that the interior integral $\mathrm{Res}$ near $a=a^{\left(
l\right) }$ approximately equals%
\begin{equation*}
\left( \frac{1}{2\pi }\right) ^{q/2}A^{-q/2}\left( \mathrm{Det}^{\left(
l\right) }\left( \eta \right) \right) ^{-1/2}\mathrm{e}^{A\phi ^{\left(
l\right) }\left( \eta \right) },
\end{equation*}%
where $\phi ^{\left( l\right) }$ is given in Eq. (2.58) and $\mathrm{Det}%
^{\left( l\right) }$ is given in Eq. (2.59).\bigskip

(c) Let us now consider the exterior integral over $\tau _{i}$'s. Recall
that the quantities $\phi ^{\left( l\right) }$ and $\mathrm{Det}^{\left(
l\right) }$ depend on $\tau $ via $\eta =\left( \tau _{1}\cdots \tau
_{q}t\right) ^{\kappa },$ $\kappa =\frac{1}{q+1}.$ Now we assume that the
principal solution, with the index $(l)$ or $l,$ is dominant; thus the
function $\phi ^{\left( l\right) }$ has positive real part in a crucial
domain.

But, by the similarity of Eqs. (2.48) and (2.55) we see that the action $%
S^{\left( l\right) }(t)$ equals $\int_{0}^{t}\rho ^{\left( l\right) }|_{\eta
=\tau ^{\kappa }}\mathrm{d}\tau .$ It means that $S^{\left( l\right) }\left(
t\right) =\phi ^{\left( l\right) }|_{\tau _{1}=\ldots -=\tau _{q}}=1.$ Thus
the term $\left( \mathrm{Det}^{\left( l\right) }\left( \eta \right) \right)
^{-1/2}$ is taken at $\eta =t^{\kappa }.$

Moreover, with $\tilde{\tau}_{i}=1-\tau _{i}$, we have%
\begin{eqnarray}
\eta  &=&t^{\kappa }\left\{ 1-\kappa \sum \tilde{\tau}_{i}+\ldots \right\} ,
\notag \\
\rho ^{\left( l\right) } &=&\rho ^{\left( l\right) }\left( \eta \right)
=\rho ^{\left( l\right) }\left( t^{\kappa }\right) -\kappa t^{\kappa }\frac{%
\mathrm{d}\rho ^{\left( l\right) }}{\mathrm{d}\eta }\left( t^{\kappa
}\right) \cdot \sum \tilde{\tau}_{i}+\ldots ,  \notag \\
\phi ^{\left( l\right) } &=&\phi ^{\left( l\right) }\left( \eta \right)
=\phi ^{\left( l\right) }\left( t^{\kappa }\right) -\frac{t^{\kappa }}{%
1+\rho ^{\left( l\right) }\left( t^{\kappa }\right) }\frac{\mathrm{d}\rho
^{\left( l\right) }}{\mathrm{d}\eta }\left( t^{\kappa }\right) \cdot \sum
\tilde{\tau}_{i}+\ldots .  \notag \\
&=&\phi ^{\left( l\right) }\left( t^{\kappa }\right) -\xi ^{\left( l\right)
}\cdot \sum \tilde{\tau}_{i}+\ldots .  \label{2.60}
\end{eqnarray}%
Therefore we get $C\left( \frac{1}{2\pi }\right) ^{q/2}A^{-q/2}\left(
\mathrm{Det}^{\left( l\right) }\left( t^{\kappa }\right) \right) ^{-1/2}%
\mathrm{e}^{A\phi ^{\left( l\right) }\left( t^{\kappa }\right) }$ times%
\begin{eqnarray*}
\prod\limits_{i}\int_{0}^{1}\tilde{\tau}_{i}^{\beta _{i}-2}\exp \left( -A\xi
^{\left( l\right) }\tilde{\tau}_{i}\right) \mathrm{d}\tilde{\tau}_{i} &\sim
&\prod \Gamma \left( \beta _{i}-1\right) \cdot \left( A\xi ^{\left( l\right)
}\right) ^{1-\beta _{i}} \\
&=&\left( \prod \Gamma \left( \beta _{i}-1\right) \right) \cdot \left( A\xi
^{\left( l\right) }\right) ^{-\beta },
\end{eqnarray*}%
where $\beta =\sum \left( \beta _{i}-1\right) $. Since $C=\prod \left( \beta
_{i}-1\right) ,$ we get that the contribution to $F\left( t\right) $ from a
neighborhood of $a^{\left( l\right) }$ equals%
\begin{equation*}
\left( \prod \Gamma \left( \beta _{i}\right) \right) \left( \frac{1}{2\pi }%
\right) ^{q/2}\cdot A^{-q/2-\beta }\cdot \mathrm{e}^{A\phi ^{\left( l\right)
}\left( t^{\kappa }\right) }\cdot \left( \mathrm{Det}^{\left( l\right)
}\left( t^{\kappa }\right) \right) ^{-1/2}\cdot \left( \xi ^{\left( l\right)
}\right) ^{-\beta }.
\end{equation*}%
\bigskip

(d) We can summarize this section in the following statement, where the
constants $E_{l}$ in the formula (2.60) are correct only in the cases the
corresponding exponents are dominating.

\begin{theorem}
\label{t4}We have the expansion%
\begin{equation}
F\sim \sum_{l=1}^{q+1}E_{l}\cdot A^{-q/2-\beta }\cdot U_{\mathrm{princ}%
}^{\left( l\right) }\sim \sum_{l=1}^{q+1}E_{l}\cdot A^{-q/2-\beta }\cdot
\mathrm{e}^{AS^{\left( l\right) }\left( t\right) }\psi _{0}^{\left( l\right)
}\left( t\right)   \label{2.61}
\end{equation}

in principal WKB solutions, where%
\begin{eqnarray*}
E_{l} &=&\left( \prod \Gamma \left( \beta _{i}\right) \right) \left( \frac{1%
}{2\pi }\right) ^{q/2}, \\
S^{\left( l\right) }\left( t\right)  &=&\phi ^{\left( l\right) }\left(
t^{\kappa }\right) , \\
\psi _{0}^{\left( l\right) }\left( t\right)  &=&\left( \mathrm{Det}^{\left(
l\right) }\left( t^{\kappa }\right) \right) ^{-1/2}\cdot \left( \xi ^{\left(
l\right) }\right) ^{-\beta }
\end{eqnarray*}%
and $\phi ^{\left( l\right) },$ $\mathrm{Det}^{\left( l\right) }$ and $\xi
^{\left( l\right) }$ are defined in Eqs. (2.58)--(2.60) above. In
particular, the functions $\phi ^{\left( l\right) }\left( t^{\kappa }\right)
$ and $\left( \mathrm{Det}^{\left( l\right) }\left( t^{\kappa }\right)
\right) ^{-1/2}\cdot \left( \xi ^{\left( l\right) }\right) ^{-\beta }$ serve
as evaluations of the above integrals (2.50)--(2.51) defining the `actions' $%
S^{\left( l\right) }$ and the `amplitudes' $\psi _{0}^{\left( l\right) }.$
\end{theorem}

\begin{example}
\label{e29}(Example 2.1 revisited) We consider the integral in Eq. (2.12)
for $u_{0}(t;\lambda )=F\left( -\lambda ,\epsilon \lambda ,\bar{\epsilon}%
\lambda ;1,1;t\right) ,$ $\epsilon =\mathrm{e}^{\mathrm{i}\pi /3}$, for
large $\lambda =A$. Since $\beta _{1}=\beta _{2}=1$ in the integral
representation we deal only with the residuum.

The Hamilton--Jacobi equation, i.e., $\left( 1-t\right) \left( \mathcal{D}%
_{t}S\right) ^{3}=-t,$ has three solutions%
\begin{equation*}
S^{\sigma }\left( t\right) =\sigma S_{0}\left( t\right) =\sigma
\int_{0}^{t}\tau ^{-2/3}\left( 1-\tau \right) ^{-1/3}\mathrm{d}\tau ,\text{ }%
\sigma =-1,\epsilon ,\bar{\epsilon}.
\end{equation*}%
The transport equation, i.e., $\left( \mathcal{D}_{t}S\right) ^{3}\cdot
\mathcal{D}_{t}\psi _{0}+\left( \mathcal{D}_{t}S\right) ^{2}\cdot \mathcal{D}%
_{t}^{2}S\cdot \mathcal{D}\psi _{0}=0,$ has the particular solution $\psi
_{0}\left( t\right) =\left( \mathcal{D}_{t}S\right) ^{-1}$, i.e.,
\begin{eqnarray*}
\psi _{0}(t) &=&1/\vartheta , \\
\vartheta &=&\left( 1/\left( 1-t\right) \right) ^{1/3}.
\end{eqnarray*}

The phase in integral (2.12) equals $\phi \left( a\right) =\ln \left(
1-t^{1/3}a_{1}\right) -\ln \left( 1-t^{1/3}a_{2}\right) -\bar{\epsilon}\ln
\left( 1-t^{1/3}a_{3}\right) ,$ subject to the restriction $%
a_{1}a_{2}a_{3}=1,$ and the measure is $\mathrm{d}^{3}\ln a/\mathrm{d}\ln
\left( a_{1}a_{2}a_{3}\right) .$ There are three conditional critical points
of the phase defined by
\begin{equation*}
t^{1/3}a_{1}^{\sigma }=\frac{-\sigma \vartheta }{1-\sigma \vartheta },\text{
\ }t^{1/3}a_{2}^{\sigma }=\frac{\bar{\epsilon}\sigma \vartheta }{1+\bar{%
\epsilon}\sigma \vartheta },\text{ \ }t^{1/3}a_{3}^{\sigma }=\frac{\epsilon
\sigma \vartheta }{1+\epsilon \sigma \vartheta },\text{ \ }\sigma
=-1,\epsilon ,\bar{\epsilon}.
\end{equation*}%
Near these points we can write $a_{j}=a_{j}^{\sigma }\mathrm{e}^{\mathrm{i}%
\theta _{j}},$ where the angles $\theta _{j}$ are subject to the restriction
$\theta _{1}+\theta _{2}+\theta _{3}=0.$ The expansion of the phase gives
the leading term%
\begin{equation*}
\phi ^{\sigma }=-\ln (1-\sigma \vartheta )+\epsilon \ln (1+\bar{\epsilon}%
\sigma \vartheta )+\bar{\epsilon}\ln (1+\epsilon \sigma \vartheta )=\sigma
S_{0}(t).
\end{equation*}%
Of course, the linear part $\varphi _{1}^{\sigma }$\ vanishes and the
quadratic part equals%
\begin{equation*}
-\frac{1}{2}Q^{\sigma }\left( \theta \right) =-\frac{1}{2}\left\{ \sigma
\vartheta \left( \theta _{1}^{2}+\theta _{2}^{2}+\theta _{3}^{2}\right)
-(\sigma \vartheta )^{2}(-\theta _{1}^{2}+\bar{\epsilon}\theta
_{2}^{2}+\epsilon \theta _{3}^{2})\right\} .
\end{equation*}%
By Lemma 2.1 the determinant of the corresponding matrix equals $%
\det^{\sigma }=3\left( \sigma \vartheta \right) ^{2}$.

As in other oscillatory type integrals (or mountain pass integrals) the
leading part of the hypergeometric function (2.12) arising from a
neighborhood of the point $a^{\sigma },$\ for large $\left\vert \lambda
\right\vert ,$ equals
\begin{equation*}
\mathrm{e}^{\sigma \lambda S_{0}(t)}\times \left( \frac{1}{2\pi \mathrm{i}}%
\right) ^{2}\int \int \mathrm{e}^{-\lambda Q^{\sigma }/2}\mathrm{i}^{2}%
\mathrm{d}\theta _{1}\mathrm{d}\theta _{2}=\mathrm{e}^{\sigma \lambda
S_{0}(t)}\times \frac{1}{4\pi ^{2}}\times \frac{2\pi }{\lambda \sqrt{%
\det^{\sigma }}}=\frac{\mathrm{e}^{\sigma \lambda S_{0}(t)}}{2\pi \sqrt{3}%
\sigma \vartheta \lambda }.
\end{equation*}%
Therefore%
\begin{equation*}
u_{0}\left( t;\lambda \right) \sim \frac{1}{2\pi \sqrt{3}\lambda }\left(
\frac{1-t}{t}\right) ^{1/3}\left\{ -\mathrm{e}^{-\lambda \widetilde{S}(t)}+%
\bar{\epsilon}\mathrm{e}^{\epsilon \lambda \widetilde{S}(t)}+\epsilon
\mathrm{e}^{\bar{\epsilon}\lambda \widetilde{S}(t)}\right\} .
\end{equation*}%
According to Remark 2.7 the summands in the above formula correspond to
principal WKB solution $U_{\mathrm{princ}}^{\sigma }\sim \vartheta ^{-1}%
\mathrm{e}^{\sigma \lambda S_{0}(t)}$. More detailed calculations, which use
a Bessel type approximation of the hypergeometric equation for large
parameter $\lambda $ and so-called Wick formula for Gaussian integrals (see
Example 3.2 below), show the relation%
\begin{equation*}
U_{\mathrm{princ}}^{\sigma }\sim \left( 1-\frac{8}{3^{6}}\lambda
^{-3}+\ldots \right) U_{\mathrm{test}}^{\sigma }
\end{equation*}%
between the principal and testing WKB solutions (which are defined in
Subsection 2.3.3 above).

In the case of the integral formula (2.11) for the function $F\left( \lambda
,-\lambda ;1;t\right) $ we have the relations $U_{\mathrm{princ}}^{\pm }\sim
\left( 1+\frac{5}{2^{8}}\lambda ^{-2}+\ldots \right) U_{\mathrm{test}}^{\pm }
$ for the corresponding WKB series $U^{\pm }\left( t\right) \sim \left(
\frac{1-t}{t}\right) ^{1/4}\mathrm{e}^{\pm \mathrm{i}\lambda \widetilde{S}%
(t)}.$
\end{example}

\section{Variations of hypergeometric functions}

Consider the following perturbation of the hypergeometric equation (2.11):%
\begin{equation}
\left( \mathcal{Q}-t\mathcal{P}-\varepsilon t^{2}\mathcal{R}\right) u=0,
\label{3.1}
\end{equation}%
where $\mathcal{R}=R\left( \mathcal{D} _{t}\right) $ with a polynomial $R$
of degree $\leq q+1$ and $\varepsilon $ is a small complex parameter. Here
the length of the recurrence for the coefficients in the power series
solutions equals 3, so we do not have direct formulas for these coefficients.

We look for solutions to Eq. (3.1) of the form%
\begin{equation}
u\left( t\right) =u\left( t;\varepsilon \right) =u_{0}\left( t\right)
+\varepsilon u_{0,1}\left( t\right) +\varepsilon ^{2}u_{0,2}\left( t\right)
+\ldots ,  \label{3.2}
\end{equation}%
where $u_{0}\left( t\right) =u_{0,0}\left( t\right) =F\left( \alpha
_{1},\ldots ,\alpha _{p};\beta _{1},\ldots ,\beta _{q};t\right) $ is our
hypergeometric functions. The functions $u_{0,j}\left( t\right) $ are called
the \textbf{variations of the hypergeometric function}.

The functions $u_{0,j}$ satisfy the following recurrence:%
\begin{equation}
\left( \mathcal{Q}-t\mathcal{P}\right) u_{0,k+1}=t^{2}\mathcal{R}u_{0,k},%
\text{ }k\geq 0,  \label{3.3}
\end{equation}%
with the initial conditions $u_{0,k}\left( 0\right) =0.$

\begin{theorem}
\label{t5}The functions $u_{0,k}$ are defined by the series%
\begin{equation}
u_{0,k}\left( t\right) =\sum_{n_{0}\geq 0,\ldots ,n_{k}\geq 0}\Omega \left(
m_{k}\right) \left( \prod\limits_{j=0}^{k-1}S\left( m_{j}\right) \right)
\frac{t^{m_{k}}}{m_{k}!},  \label{3.4}
\end{equation}%
where%
\begin{equation*}
m_{j}=n_{0}+\ldots +n_{i}+2j
\end{equation*}%
and%
\begin{equation}
\Omega \left( n\right) =\frac{\left( \alpha _{1}\right) _{n}\cdots \left(
\alpha _{p}\right) _{n}}{\left( \beta _{1}\right) _{n}\cdots \left( \beta
_{q}\right) _{n}},\text{ }S\left( n\right) =\frac{R\left( n\right) Q\left(
n+1\right) }{P\left( n\right) P\left( n+1\right) }.  \label{3.5}
\end{equation}

Moreover, these variations admit integral representations, following from
the below proof.
\end{theorem}

\textit{Proof}. The operator $\mathcal{Q}-t\mathcal{P}$ acts on the space $t%
\mathbb{C}\left[ \left[ t\right] \right] $ of power series in $t$ without
constant term. Assuming $Q\left( n\right) \not=0$ for $n>0,$ the operator $%
\mathcal{Q}$ is invertible. Moreover, the operator $t\mathcal{P}$ is `small'
in relation with $\mathcal{Q}$.\footnote{%
One can introduce a suitable norm in the space $t\mathbb{C}\left[ \left[ t%
\right] \right] ,$ such that the series of finite norm are convergent in a
disc of given small radius. Then the operator $t\mathcal{P}$ will have small
norm.} Therefore, $\mathcal{Q}-t\mathcal{P}=\left( I-t\mathcal{PQ}%
^{-1}\right) \mathcal{Q}$ is invertible with the inverse defined by the von
Neumann type series%
\begin{eqnarray*}
\left( \mathcal{Q}-t\mathcal{P}\right) ^{-1} &=&\mathcal{Q}^{-1}\left( I-t%
\mathcal{PQ}^{-1}\right) ^{-1} \\
&=&\mathcal{Q}^{-1}+\mathcal{Q}^{-1}t\mathcal{PQ}^{-1}+\mathcal{Q}%
^{-1}\left( t\mathcal{PQ}^{-1}\right) ^{2}+\ldots .
\end{eqnarray*}%
In particular, we have%
\begin{equation*}
\left( \mathcal{Q}-t\mathcal{P}\right) ^{-1}t^{n}=\frac{1}{Q\left( n\right) }%
t^{n}+\frac{P\left( n\right) }{Q\left( n\right) Q\left( n+1\right) }t^{n+1}+%
\frac{P\left( n\right) P\left( n+1\right) }{Q\left( n\right) Q\left(
n+1\right) Q\left( n+2\right) }t^{n+2}+\ldots .
\end{equation*}%
Also%
\begin{equation*}
t^{2}\mathcal{R}t^{n}=R\left( n\right) t^{n+2}.
\end{equation*}%
Using this Eq. (3.4) is proved by induction. In particular, the term with
given $n_{0},\ldots ,n_{k}$ equals%
\begin{equation*}
\mathcal{Q}^{-1}\left( t\mathcal{PQ}^{-1}\right) ^{n_{k}}\circ t^{2}\mathcal{%
R}\circ \ldots \circ \mathcal{Q}^{-1}\left( t\mathcal{PQ}^{-1}\right)
^{n_{1}}\circ t^{2}\mathcal{R}\frac{\left( \alpha _{1}\right) _{n_{0}}\cdots
\left( \alpha _{p}\right) _{n_{0}}}{\left( \beta _{1}\right) _{n_{0}}\cdots
\left( \beta _{q}\right) _{n_{0}}}\frac{t^{n_{0}}}{n_{0}!}.
\end{equation*}%
\medskip

Let us prove the second statement of Theorem 5. We do not write down an
integral formula for the function in the right-hand side of Eq. (3.4). We
present only the steps leading to it.

Firstly, we denote%
\begin{equation*}
p_{0}=n_{0},\text{ }p_{j}=n_{j}+2\text{ for }j>0.
\end{equation*}%
Then in Eq. (3.4) we have the sum over $p_{0}\geq 0$ and $p_{j}\geq 2.$
Introduce the following function of many variables%
\begin{equation*}
w\left( t_{0},t_{1},\ldots ,t_{k}\right) =\sum_{p_{0}\geq 0,p_{1}\geq
2\ldots ,p_{k}\geq 2}\Omega \left( p_{0}+\ldots +p_{k}\right) \frac{%
t_{0}^{p_{0}}\cdots t_{k}^{p_{k}}}{\left( p_{0}+\ldots +p_{k}\right) !}.
\end{equation*}%
We have%
\begin{equation*}
v_{0,k}\left( t\right) =\left\{ \left( \prod\limits_{j=0}^{k-1}S\left(
\mathcal{D} _{t_{0}}+\ldots +\mathcal{D} _{t_{j}}\right) \right) w\right\}
|_{t_{0}=\ldots =t_{k}=t},
\end{equation*}%
where $\mathcal{D} _{t_{j}}=t_{j}\partial /\partial t_{j}$ are Euler
operators. The operators%
\begin{eqnarray*}
&&S\left( \mathcal{D} _{t_{0}}+\ldots +\mathcal{D} _{t_{j}}\right) =P\left(
\mathcal{D} _{t_{0}}+\ldots +\mathcal{D} _{t_{j}}\right) ^{-1}P\left(
\mathcal{D} _{t_{0}}+\ldots +\mathcal{D} _{t_{j}}+1\right) ^{-1} \\
&&\text{ \ \ \ \ \ \ \ \ \ \ \ \ \ \ \ \ \ \ \ \ \ \ \ \ \ }R\left( \mathcal{%
D} _{t_{0}}+\ldots +\mathcal{D} _{t_{j}}\right) Q\left( \mathcal{D}
_{t_{0}}+\ldots +\mathcal{D} _{t_{j}}+1\right)
\end{eqnarray*}%
need an interpretation. There are no problems with the operators $R\left(
\mathcal{D} _{t_{0}}+\ldots +\mathcal{D} _{t_{j}}\right) $ and $Q\left(
\mathcal{D} _{t_{0}}+\ldots +\mathcal{D} _{t_{j}}+1\right) ;$ they are
standard differential operators.

Let us interpret the operator $\left( \mathcal{D} _{t_{0}}+\ldots +\mathcal{D%
} _{t_{j}}+\alpha \right) ^{-1}.$ It amounts to solving the equation%
\begin{equation*}
\left( \mathcal{D} _{t_{0}}+\ldots +\mathcal{D} _{t_{j}}+\alpha \right) f=g,
\end{equation*}%
where we may assume that $g=g\left( t_{0},\ldots ,t_{j}\right) .$ Introduce
the variables $r=t_{0}+\ldots +t_{j},$ $s_{i}=t_{i}/r;$ thus $t=rs.$ If $%
f\left( rs\right) =\sum f_{m}\left( s\right) r^{m}$ and $g=\sum g_{m}\left(
s\right) r^{m},$ then $\left\{ \left( \mathcal{D} _{t_{0}}+\ldots +\mathcal{D%
} _{t_{j}}+\alpha \right) f\right\} \left( rs\right) =\sum \left( m+\alpha
\right) f_{m}\left( s\right) r^{m}$ and hence $f\left( rs\right) =\sum \frac{%
g_{m}\left( s\right) }{m+s}r^{m}$ $=r^{-\alpha }\int_{0}^{r}\tau ^{\alpha
-1}g\left( \tau s\right) \mathrm{d}\tau ,$ or%
\begin{equation*}
f\left( t\right) =\left( t_{0}+\ldots +t_{j}\right) ^{-\alpha
}\int_{0}^{t_{0}+\ldots +t_{j}}\tau ^{\alpha -1}g\left( \frac{\tau }{%
t_{0}+\ldots +t_{j}}t\right) \mathrm{d}\tau .
\end{equation*}%
Here one assumes $\mathrm{Re}\alpha >0;$ in the opposite case (but $\alpha
\not\in \mathbb{Z})$ one applies it to part of $g$ with sufficiently high
order terms and for the few low order terms one uses $g_{m}\left( s\right)
r^{m}\longmapsto \frac{1}{m+s}g_{m}\left( s\right) r^{m}.$

Next, consider the following series:%
\begin{equation*}
w_{0}\left( t_{0},\ldots ,t_{k}\right) =\sum_{p_{0}\geq 0,p_{1}\geq 0\ldots
,p_{k}\geq 0}\Omega \left( p_{0}+\ldots +p_{k}\right) \frac{%
t_{0}^{p_{0}}\cdots t_{k}^{p_{k}}}{\left( p_{0}+\ldots +p_{k}\right) !},
\end{equation*}%
i.e., a hypergeometric function of many variables; some of them are known as
Appel functions or Lauricella functions (see \cite[Section 5.7]{BE1} and
\cite[Chapter 6.5]{Koh}). The difference between $w_{0}$ and $w$ is a sum of
functions expressed via hypergeometric series of smaller number of
variables; we do not analyze it in detail.

In order to deal with the series $w_{0}$ we need generalizations of formulas
before the proof of Theorem 2.1. One is the generalized binomial formula%
\begin{equation*}
\left( 1-x_{1}-\ldots -x_{l}\right) ^{-\alpha }=\sum_{n_{1}\geq 0,\ldots
,n_{l}\geq 0}\frac{\left( \alpha \right) _{n_{1}+\ldots +n_{l}}}{%
n_{1}!\cdots n_{l}!}x_{1}^{n_{1}}\cdots x_{l}^{n_{l}}.
\end{equation*}%
Another is the generalized Beta-function%
\begin{equation*}
B\left( \gamma _{0},\ldots ,\gamma _{l}\right) =\int \tau _{0}^{\gamma
_{0}-1}\cdots \tau _{l}^{\gamma _{l}}\frac{\mathrm{d}^{l+1}\tau }{\mathrm{d}%
\left( \tau _{0}+\ldots +\tau _{l}\right) }=\frac{\Gamma \left( \gamma
_{0}\right) \cdots \Gamma \left( \gamma _{l}\right) }{\Gamma \left( \gamma
_{0}+\ldots +\gamma _{l}\right) },
\end{equation*}%
where the integration runs over the simplex $\left\{ \tau _{j}\geq 0,\tau
_{0}+\ldots +\tau _{l}\leq 1\right\} .$

Finally we use the residuum formula, but with increased number of complex
integration arguments. In particular, in the case $p=q+1,$ the number of
arguments is $\left( k+1\right) \times \left( 2q-1\right) :$ $\left(
q+1\right) $ of $c_{i}$'s and $q=p-1$ of $b_{i}$'s for each variable $t_{j}$
(compare the proof of Theorem 2.1). The integral--residua with respect to $%
c_{i}$'s are evaluated using the residuum theorem and there remain only
residua with respect to $b_{i}$'s and the integrations over $\tau _{i}$'s.

In the case $p<q+1$ the integral is suitably modified, with exponents
replacing some binomial formula. \qed\medskip

The integral formulas from the thesis of Theorem 5 are quite complicated and
not illuminating; so we do not write them down. In the below examples we
present them in some cases which are relatively simple.

\begin{example}
\label{e31}Consider the perturbation%
\begin{equation*}
\left\{ \mathcal{D}_{t}-t\left( \mathcal{D}_{t}+\alpha \right) -\varepsilon
t^{2}\left( \mathcal{D}_{t}+\gamma \right) \right\} u=0
\end{equation*}%
of the hypergeometric equation with $p=q+1=1;$ it probably a simples
example. Recall that for $\varepsilon =0$ the solution $u_{0}\left( t\right)
=\left( 1-t\right) ^{-\alpha },$ see Eq. (2.7). The above equation is of
first order, $\left( 1-t-\varepsilon t^{2}\right) \dot{u}=\left( \alpha
+\varepsilon \gamma t\right) u,$ and has solution%
\begin{eqnarray*}
u_{0}\left( t;\varepsilon \right) &=&\exp \int_{0}^{t}\frac{\alpha
+\varepsilon \gamma s}{1-s-\varepsilon s^{2}}\mathrm{d}s \\
&=&\left( 1-t\right) ^{-\alpha }\left\{ 1+\varepsilon \left[ \left( 2-\gamma
\right) \ln \left( 1-t\right) +t\frac{2-t}{1-t}-t\right] +O\left(
\varepsilon ^{2}\right) \right\} .
\end{eqnarray*}%
\medskip
\end{example}

\begin{example}
\label{e32}The function $F\left( -\lambda ,\epsilon \lambda ,\bar{\epsilon}%
\lambda ;1,1;t\right) $ from Example 2.1 for large $\lambda $ and small $t$
such that $y=\lambda ^{3}t$ is finite can be rewritten as $F=U_{0}\left(
y\right) +\lambda ^{-3}U_{0,1}\left( y\right) +\ldots ,$ where $U_{0}\left(
y\right) =\sum_{n\geq 0}\frac{\left( -y\right) ^{n}}{\left( n!\right) ^{3}}$
is Bessel like confluent hypergeometric function and $U_{0,1}=\sum_{m\geq
0,n\geq 0}\frac{n^{3}}{\left( \left( m+n\right) !\right) ^{3}}\left(
-y\right) ^{m+n},$ $U_{0,2},$ etc. are variations of $U_{0}.$ Here the
analogue of the equation (3.1) is simpler. But the integral formulas for $%
U_{0,j}\left( y\right) $ can be expressed directly from Eq. (2.12). We have%
\begin{eqnarray*}
U_{0}\left( y\right) &=&\mathrm{Res}_{a=0}\exp \left( -y^{1/3}\phi \right)
\frac{\mathrm{d}^{3}\ln a}{\mathrm{d}\ln \left( a_{1}a_{2}a_{3}\right) }, \\
U_{0,1}\left( y\right) &=&\mathrm{Res}_{a=0}\exp \left( -y^{1/3}\phi \right)
A\frac{\mathrm{d}^{3}\ln a}{\mathrm{d}\ln \left( a_{1}a_{2}a_{3}\right) },
\end{eqnarray*}%
where the phase $\phi =\phi _{0},$%
\begin{equation*}
\phi _{k}=\frac{1}{k+1}\left( a_{1}^{k+1}-\epsilon a_{2}^{k+1}-\bar{\epsilon}%
a_{3}^{k+1}\right) ,
\end{equation*}%
and the amplitude%
\begin{equation*}
A=-\frac{1}{3!}y^{2}\phi _{1}^{3}+y^{5/3}\phi _{1}\phi _{2}-y^{4/3}\phi _{3}.
\end{equation*}%
\medskip

Consider variations of the function $V_{2}\left( z\right) $ from Example 2.6
(see Eq. (2.22)). Recall that we have $v_{2}\left( s;\lambda \right)
=\lambda ^{3}s^{2}+\ldots =V_{2}\left( z\right) +\ldots $ is one of the
basic solutions of the corresponding hypergeometric equation near $s=1-t=0.$
Moreover, we assume small $s$ and large $\lambda $ such that $z=\lambda
^{3}s^{2}$ is finite. In fact, the perturbed equation becomes more
complicated than Eq. (3.1); we have%
\begin{equation}
\left( \mathcal{Q}-z-\lambda ^{-3/2}z^{1/2}\mathcal{R}+\lambda ^{-3}z%
\mathcal{S}\right) V=0,  \label{3.6}
\end{equation}
where $\mathcal{Q}=2\mathcal{D} _{z}\left( 2\mathcal{D} _{z}-1\right) \left(
2\mathcal{D} _{z}-2\right) ,$ $\mathcal{R}=2\mathcal{D} _{z}\left( 2\mathcal{%
D} _{z}-1\right) \left( 4\mathcal{D} _{z}-1\right) $ and $\mathcal{R}=\left(
2\mathcal{D} _{z}\right) ^{3}.$

Therefore $v_{2}=V_{2}\left( z\right) +\lambda ^{-3/2}V_{2,1}\left( z\right)
+\lambda ^{-3}V_{2,2}\left( z\right) +\ldots .$ $\ $Calculations show that%
\begin{eqnarray*}
V_{2,1} &=&\frac{8}{\sqrt{z}}\sum_{m,n\geq 1}\frac{\left( 4n-1\right) \left(
1/2\right) _{n}}{\left( 2m+2n-1\right) !\left( 1/2\right) _{m+n-1}\left(
n-1\right) !}\left( \frac{z}{2}\right) ^{m+n} \\
&=&\frac{8}{z}\left( 4\mathcal{D}_{z_{2}}-1\right) |_{z_{j}=z}\mathrm{Res}%
_{a=0}\Lambda _{1}\left( \theta \right) \Lambda _{2}\left( \vartheta \right)
\mathrm{d}^{2}\ln a,
\end{eqnarray*}%
where%
\begin{equation*}
\theta _{j}=z_{j}^{1/3}a_{j},\text{ }\vartheta _{j}=z_{j}^{1/3}/2a_{j}^{2},%
\text{ }j=1,2,
\end{equation*}%
and%
\begin{eqnarray*}
\Lambda _{1} &=&\sum_{m,n\geq 1}\frac{\theta _{1}^{2m}\theta _{2}^{2n}}{%
\left( 2m+2n-1\right) !}=\sum_{p\geq 2}\frac{h_{p}\left( \theta
_{1}^{2},\theta _{2}^{2}\right) }{\left( 2p-1\right) !} \\
&=&\frac{\theta _{1}^{3}\sinh \theta _{1}-\theta _{2}^{3}\sinh \theta _{2}}{%
\theta _{1}^{2}-\theta _{2}^{2}}-\left( \theta _{1}^{2}+\theta
_{2}^{2}\right) ,
\end{eqnarray*}%
\begin{eqnarray*}
\Lambda _{2} &=&\sum_{m,n\geq 1}\frac{\left( 1/2\right) _{n}}{\left(
1/2\right) _{m+n-1}\left( n-1\right) !}\vartheta _{1}^{m}\vartheta _{2}^{n}
\\
&=&\vartheta _{1}\vartheta _{2}\left( \mathcal{D}_{\vartheta _{1}}+\mathcal{D%
}_{\vartheta _{2}}+\frac{3}{2}\right) \int_{0}^{1}\tau ^{1/2}\mathrm{e}%
^{\left( 1-\tau \right) \vartheta _{1}+\tau \vartheta _{2}}\mathrm{d}\tau
\end{eqnarray*}%
and $h_{p}\left( x,y\right) $ is the complete symmetric polynomial of degree
$p.$ Also the second variation $V_{2,2}$ was found, but we do not present
its rather complicated expression.
\end{example}

\section{Reflections about the Wasow theorem}

W. Wasow in \cite[Theorem 3.3]{Was} studied perturbations of the Airy
equation from the point of view of asymptotic expansions and proved its
`equivalence' with the very Airy equation. The aim of this section is to
take somewhat different point of view to this subject. We think that our
approach is simpler and more correct.\medskip

Recall that the Airy equation (2.17), i.e., $\ddot{u}=tu$ (or $\left(
\mathcal{D} _{t}^{2}-\mathcal{D} _{t}-t^{3}\right) u=0)$ has as basic
solutions either $u_{1}\left( t\right) =F\left( \emptyset
;2/3;t^{3}/9\right) $ and $u_{2}\left( t\right) =tF\left( \emptyset
;4/3;t^{3}/9\right) ,$ or $\mathrm{Ai}\left( t\right) $ and $\mathrm{Ai}%
\left( \mathrm{e}^{2\pi \mathrm{i}/3}t\right) $ (see Eq. (2.18)).

Wasow applies the change
\begin{equation}
t=\varepsilon ^{-2/3}x  \label{4.1}
\end{equation}%
and rewrites the Airy equation as%
\begin{equation}
\varepsilon ^{2}u^{\prime \prime }=xu,  \label{4.2}
\end{equation}%
where $^{\prime }=\frac{\mathrm{d}}{\mathrm{d}x}.$ Then, with $u=y_{1}$ and $%
\varepsilon y_{1}^{\prime }=y_{2}$, he gets the system%
\begin{equation}
\varepsilon Y^{\prime }=\boldsymbol{A}_{0}\left( x\right) Y,  \label{4.3}
\end{equation}%
where $\boldsymbol{A}_{0}=\left(
\begin{array}{ll}
0 & 1 \\
x & 0%
\end{array}%
\right) .$

Next, Eq. (4.1) is generalized to%
\begin{equation}
\varepsilon ^{2}u^{\prime \prime }=\left( x\varphi \left( x\right)
+\varepsilon \psi \left( x;\varepsilon \right) \right) u  \label{4.4}
\end{equation}%
(where the functions $\phi $ and $\psi $ are analytic), or to the system%
\begin{equation}
\varepsilon Y^{\prime }=\boldsymbol{A}\left( x,\varepsilon \right) Y,\text{
\ }\boldsymbol{A}\left( x,0\right) =\boldsymbol{A}_{0}\left( x\right)
\label{4.5}
\end{equation}%
(where $\boldsymbol{A}$ is defined by Eq. (4.4) in the same way as $%
\boldsymbol{A}_{0}$ is defined by Eq. (4.2)). Wasow develops tools to prove
that systems (4.3) and (4.5) are equivalent. Firstly he simplifies system
(4.4) for $\varepsilon =0$ (see \cite[Theorem 29.1]{Was}).

\begin{lemma}
\label{l41}Assume the germs $\varphi \left( x\right) =1+\ldots $ and $\psi
\left( x;\varepsilon \right) $ analytic in $\left( \mathbb{C},0\right) $\
and $\left( \mathbb{C}^{2},0\right) $ respectively. Then there exists an
analytic change%
\begin{equation*}
z=a\left( x\right) =x+\ldots ,\text{ }u=b\left( x\right) w=\left( 1+\ldots
\right) w
\end{equation*}%
reducing Eq. (4.4) to%
\begin{equation}
\varepsilon ^{2}\frac{\mathrm{d}^{2}w}{\mathrm{d}z^{2}}=\left( z+\varepsilon
\chi \left( z;\varepsilon \right) \right) w  \label{4.6}
\end{equation}%
with analytic germ $\chi .$
\end{lemma}

\textit{Proof}. Let $c\left( x\right) =a^{\prime }\left( x\right) =\frac{%
\mathrm{d}z}{\mathrm{d}x};$ thus $\frac{\mathrm{d}}{\mathrm{d}x}=c\left(
x\right) \frac{\mathrm{d}}{\mathrm{d}z}$. In particular, $w^{\prime }=c\frac{%
\mathrm{d}w}{\mathrm{d}z}$ and $w^{\prime \prime }=\left( c\frac{\mathrm{d}w%
}{\mathrm{d}z}\right) ^{\prime }=c^{\prime }\frac{\mathrm{d}w}{\mathrm{d}z}+c%
\frac{^{2}\mathrm{d}^{2}w}{\mathrm{d}z^{2}}$. We have%
\begin{equation*}
u^{\prime \prime }=b^{\prime \prime }w+2b^{\prime }w^{\prime }+bw^{\prime
\prime }=b^{\prime \prime }w+\left( 2b^{\prime }c+bc^{\prime }\right) \frac{%
\mathrm{d}w}{\mathrm{d}z}+c^{2}b\frac{\mathrm{d}^{2}w}{\mathrm{d}z^{2}}
\end{equation*}%
and hence we get the equation%
\begin{equation*}
\varepsilon ^{2}c^{2}b\frac{\mathrm{d}^{2}w}{\mathrm{d}z^{2}}+\varepsilon
^{2}\left( 2b^{\prime }c+bc^{\prime }\right) \frac{\mathrm{d}w}{\mathrm{d}z}%
=\left( \varphi b-\varepsilon ^{2}b^{\prime \prime }+\varepsilon \psi
b\right) w.
\end{equation*}%
We do not want the terms with $\frac{\mathrm{d}w}{\mathrm{d}z};$ so, we
assume%
\begin{equation*}
c=b^{-2}
\end{equation*}%
and we arrive at the equation%
\begin{equation*}
\varepsilon ^{2}\frac{\mathrm{d}^{2}w}{\mathrm{d}z^{2}}=\left\{ \varphi
b^{4}-\varepsilon ^{2}b^{3}b^{\prime \prime }+\varepsilon \psi b^{4}\right\}
w.
\end{equation*}%
We want $\varphi b^{4}=z,$ which amounts to the differential equation $%
\left( \frac{\mathrm{d}x}{\mathrm{d}z}\right) ^{2}\varphi \left( x\right)
=z. $ Therefore,%
\begin{equation*}
z=a\left( x\right) =\left( \frac{2}{3}\int_{0}^{x}\sqrt{\varphi \left(
s\right) }\mathrm{d}s\right) ^{2/3}.
\end{equation*}%
\qed\medskip

Next in \cite[Theorem 29.2]{Was} it is proved that there exists a formal
series $\widehat{\boldsymbol{P}}\left( x;\varepsilon \right) =\sum_{n\geq 0}%
\boldsymbol{P}_{n}\left( x\right) \varepsilon ^{n},$ with matrix-valued
functions $\boldsymbol{P}_{n}\left( x\right) $ holomorphic near $x=0,$ such
that the change $Y=\widehat{\boldsymbol{P}}Z$ transforms Eq. (4.5) to Eq.
(4.3). In \cite[Theorem 30.1 and Theorem 30.2]{Was} it is proved that there
exist analytic representations of the series $\widehat{\boldsymbol{P}}$ in
some sectorial domains; this suggests a nontrivial Stokes phenomenon.
Finally in \cite[Theorem 30.3]{Was} one finds the following result.

\begin{theorem}
\label{t6}There exists a matrix-valued functions $\boldsymbol{P}\left(
x;\varepsilon \right) ,$ with the expansion $\widehat{\boldsymbol{P}}$ as
above, but holomorphic in $\left( \mathbb{C},0\right) \times \left( \mathbb{C%
},0\right) ,$ which transforms Eq. (4.5) with $\varphi \left( x\right) =x$
to Eq. (4.3). More precisely, if we put%
\begin{equation*}
Y=P(x;\varepsilon )Z
\end{equation*}%
and assume that $Z=(z_{1}(x),z_{2}(x))^{\top }$ satisfies Eq. (4.5) (with $Y$
replaced with $Z)$, then $Y=\left( y_{1}(x),y_{2}(x)\right) ^{\top }$
satisfies Eq. (4,3).
\end{theorem}

\textit{Proof}. Let us apply the change inverse to the change (4.1) to Eq.
(4.4). We get%
\begin{equation}
\ddot{u}=\left( t+\mu \tilde{\psi}\left( t;\mu \right) \right) u,
\label{4.7}
\end{equation}%
where $\mu =\varepsilon ^{1/3}$ and $\tilde{\psi}\left( t;\mu \right) =\psi
\left( \mu ^{2}t,\mu ^{3}\right) .$ Note that both Eqs. (4.2) and (4.7) have
analytic right-hand sides.

Let $\left\{ v_{1}\left( t\right) ,v_{2}\left( t\right) \right\} $ be the
fundamental system of solutions to Eq. (4.7) defined by the initial
conditions: $v_{1}\left( 0\right) =1,$ $\dot{v}_{1}\left( 0\right) =0$ and $%
v_{2}\left( 0\right) =0,$ $\dot{v}_{2}\left( 0\right) =1;$ recall that the
solutions $u_{1}\left( t\right) $ and $u_{2}\left( t\right) $ to the Airy
equation satisfy the same initial conditions. Define $\boldsymbol{F}\left(
t\right) =\boldsymbol{F}\left( t;\mu \right) =\left(
\begin{array}{ll}
v_{1} & v_{2} \\
\dot{v}_{1} & \dot{v}_{2}%
\end{array}%
\right) $ and $\boldsymbol{F}_{0}\left( t\right) $ as the fundamental
matrices corresponding to these two fundamental systems respectively; they
are analytic and invertible.

Then the matrix%
\begin{equation*}
\boldsymbol{Q}\left( t;\mu \right) =\boldsymbol{FF}_{0}^{-1}
\end{equation*}%
is analytic and transforms the system $\dot{u}=v,$ $\dot{v}=\left( t+\mu
\tilde{\psi}\right) u$ (associated with Eq. (4.7) to the system $\dot{p}=q,$
$\dot{q}=tp$ (associated with the Airy equation $\ddot{p}=tp$).

To obtain the matrix $P$ from the thesis of the theorem one should apply the
change (4.1). With the fundamental matrices%
\begin{eqnarray*}
\widetilde{\boldsymbol{F}}\left( x;\varepsilon \right) &=&\left(
\begin{array}{ll}
\tilde{v}_{1}\left( x\right) & \tilde{v}_{2}\left( x\right) \\
\varepsilon \tilde{v}_{1}^{\prime }\left( x\right) & \varepsilon \tilde{v}%
_{2}^{\prime }\left( x\right)%
\end{array}%
\right) =\left(
\begin{array}{ll}
v_{1}\left( \varepsilon ^{-2/3}x\right) & \varepsilon ^{2/3}v_{2}\left(
\varepsilon ^{-2/3}x\right) \\
\varepsilon ^{1/3}\dot{v}_{1}\left( \varepsilon ^{-2/3}x\right) &
\varepsilon ^{1/3}\dot{v}_{2}\left( \varepsilon ^{-2/3}x\right)%
\end{array}%
\right) , \\
&=&\boldsymbol{CF}\left( x;\varepsilon \right) \boldsymbol{D},
\end{eqnarray*}%
where $\tilde{v}_{1}\left( x\right) =1+O\left( x^{2}\right) ,$ $\tilde{v}%
_{2}\left( x\right) =x+O\left( x^{2}\right) =\varepsilon ^{2/3}t+O\left(
t^{2}\right) ,$ $\boldsymbol{C}=\left(
\begin{array}{ll}
1 & 0 \\
0 & \varepsilon ^{1/3}%
\end{array}%
\right) $ and $\boldsymbol{D}=\left(
\begin{array}{ll}
1 & 0 \\
0 & \varepsilon ^{2/3}%
\end{array}%
\right) $ and analogous fundamental matrix $\widetilde{\boldsymbol{F}}%
_{0}\left( x;\varepsilon \right) ,$ we find%
\begin{equation*}
\boldsymbol{P}=\widetilde{\boldsymbol{F}}\widetilde{\boldsymbol{F}}_{0}^{-1}=%
\boldsymbol{CQC}^{-1}.
\end{equation*}%
This matrix is holomorphic in $x$ and expands into powers of $\varepsilon $
(not of $\varepsilon ^{1/3}$) as the corresponding fundamental systems
expand into powers of $\varepsilon $.

Finally we note that Wasow used the fundamental system $\left\{ \mathrm{Ai}%
\left( t\right) ,\mathrm{Ai}\left( e^{2\pi i/3}t\right) \right\} $, instead
of $\left\{ u_{1}\left( t\right) ,u_{2}\left( t\right) \right\} $. \qed%
\medskip

In system (4.5) the small parameter $\varepsilon $ stands in the left-hand
side; but this system can be rewritten as $\dot{Y}=\varepsilon ^{-1}%
\boldsymbol{A}y,$ i.e., with large parameter in the right-hand side. In \cite%
[Section 29.1]{Was} the author writes that `the point $x=0$ is the turning
point' from the WKB analysis point of view.

So, he claims that the Airy equation and its perturbation have WKB solutions
on one side and are analytically equivalent on other side. Therefore, the
Airy equation and its perturbation should experience the same Stokes
phenomena.

In \cite{ZZ3} we have tried to establish analogous equivalence between
equation $\left( \mathcal{Q}_{0}-z\right) V=0$ in Eq. (2.23) and the
hypergeometric equation ($\mathcal{Q}_{0}-z-\lambda ^{-3/2}z^{1/2}\mathcal{R}
$ $+\lambda ^{-3}z\mathcal{S)}V=0$ in Eq. (3.6); it is associated with the
third order hypergeometric equation (2.19) near $s=1-t=0$ with large
variable $z=\lambda ^{3}s^{2}$ and small parameter $\varepsilon =\lambda
^{-3/2}$ (see Example 2.6). Unfortunately, there is no such equivalence as
the below discussion (in Section 5.2) demonstrates.

Our proof of the Wasow theorem demonstrates that there is no Stokes
phenomena for the Airy equation and its perturbation under assumptions of
this theorem. In fact, the large parameter $\varepsilon ^{-1}$ is introduced
purely artificially and it is wrong to say about some WKB analysis.\bigskip

We will finish this section by calculation of the first variation of the
confluent hypergeometric solution $u_{1}\left( t\right) $ for a perturbation
of the Airy equation.

\begin{example}
\label{e41}Consider the equation%
\begin{equation}
\ddot{u}=\left( t+\varepsilon t^{2}\right) u,  \label{4.8}
\end{equation}%
where $\varepsilon $ is a small parameter. The solution of the initial value
problem $u\left( 0\right) =1,$ $\dot{u}\left( 0\right) =0$ is expanded into
power of $\varepsilon :$%
\begin{equation*}
u\left( t;\varepsilon \right) =u_{1}\left( t\right) +\varepsilon
u_{1,1}\left( t\right) +O\left( \varepsilon ^{2}\right) ,
\end{equation*}%
where $u_{1}\left( t\right) =F\left( \emptyset ;2/3;t^{3}/9\right)
=1+t^{3}/\left( 2\cdot 3\right) +\ldots $ and the first variation $u_{1,1}$
satisfies the equation%
\begin{equation*}
\left( \partial ^{2}-t\right) u_{1,1}=t^{2}u_{1},
\end{equation*}%
$\partial =\partial _{t}=\frac{\partial }{\partial t}.$ So, we have to
define the inverse to the operator $\partial ^{2}-t=\left( 1-t\partial
^{-2}\right) \partial ^{2},$ i.e., in the space of power series. Here the
multiplication operator $t$ is small with respect to $\partial ^{2}$ (we
have $t\left( t^{k}\right) =t^{k+1}$ and $\partial ^{-2}\left( t^{k}\right) =%
\frac{1}{\left( k+1\right) \left( k+2)\right) }t^{k+2});$ hence%
\begin{equation*}
\left( \partial ^{2}-t\right) ^{-1}=\partial ^{-2}+\partial ^{-2}t\partial
^{-2}+\ldots .
\end{equation*}%
We have%
\begin{equation*}
\left( \partial ^{2}-t\right) ^{-1}t^{n}=\frac{t^{n+2}}{\left( n+1\right)
\left( n+2\right) }+\frac{t^{n+5}}{\left( n+1\right) \left( n+2\right)
\left( n+4\right) \left( n+5\right) }+\ldots .
\end{equation*}%
This leads to the formula%
\begin{equation}
u_{1,1}=-\frac{3}{t}\sum_{l\geq 1,m\geq 0}\frac{\left( 1/3\right) _{l}z^{l+m}%
}{\left( -1/3\right) _{l}\left( 1/3\right) _{l+m}\left( l+m\right) !},
\label{4.9}
\end{equation}%
$z=t^{3}/9.$ Above $\sum_{l\geq 1,m\geq 0}=\sum_{l,m\geq 0}-\sum_{l=0,m\geq
0},$ where the second sum equals $F\left( \emptyset ;1/3;z\right) $ (with
known integral representation). So, we consider the series%
\begin{equation}
G\left( z_{1},z_{2}\right) =\sum_{l,m\geq 0}\frac{\left( 1/3\right)
_{l}z_{1}^{l}z_{2}^{m}}{\left( -1/3\right) _{l}\left( 1/3\right)
_{l+m}\left( l+m\right) !};  \label{4.10}
\end{equation}%
we have%
\begin{equation}
u_{1,1}\left( t\right) =-\frac{3}{t}G\left( t^{3}/9,t^{3}/9\right) -\frac{3}{%
t}F\left( \emptyset ;1/3;t^{3}/9\right) .  \label{4.11}
\end{equation}

To get an integral for $G$ we represent its summand as $\Gamma \left(
-1/3\right) ^{-1}$ times%
\begin{equation*}
\frac{\Gamma \left( -1/3\right) l!}{\Gamma \left( -1/3+l\right) }\cdot \frac{%
\Gamma \left( 1/3\right) l!m!}{\Gamma \left( 1/3+l+m\right) }\cdot \frac{%
\Gamma \left( 1/3+l\right) }{\Gamma \left( 1/3\right) l!}\cdot \frac{1}{%
\left( l+m\right) !}\cdot \frac{z_{1}^{l}}{l!}\cdot \frac{z_{2}^{m}}{m!}.
\end{equation*}%
In order to represent is as a residuum, we change it to%
\begin{eqnarray*}
&&\frac{\Gamma \left( -1/3\right) l!\left( c_{1}^{4}\right) ^{l_{1}}}{\Gamma
\left( -1/3+l_{1}\right) }\cdot \frac{\Gamma \left( 1/3\right)
l_{2}!m_{1}!\left( c_{2}^{3}/c_{1}\right) ^{l_{2}}\left( e^{2}\right)
^{m_{1}}}{\Gamma \left( 1/3+l_{2}+m_{1}\right) }\cdot \frac{\Gamma \left(
1/3+l_{3}\right) \left( b_{1}^{2}/c_{1}c_{2}\right) ^{l_{3}}}{\Gamma \left(
1/3\right) l_{3}!} \\
&&\cdot \frac{\left( b_{2}/c_{1}c_{2}b_{1}\right) ^{l_{4}}\left(
d_{1}/e\right) ^{m_{2}}}{\left( l_{4}+m_{2}\right) !}\cdot \frac{\left(
z_{1}b_{3}/b_{1}b_{2}c_{1}c_{2}\right) ^{l_{5}}}{l_{5}!}\cdot \frac{\left(
z_{2}d_{2}/ed_{1}\right) ^{m_{3}}}{m_{3}!}
\end{eqnarray*}%
and sum it up. We get%
\begin{eqnarray*}
&&\int_{0}^{1}\left( 1-\tau \right) ^{-\frac{7}{3}}\frac{\mathrm{d}\tau }{%
1-c_{1}^{4}\tau }\cdot \int \int \rho _{1}^{-\frac{5}{3}}\frac{\mathrm{d}%
^{3}\rho /\mathrm{d}\left( \rho _{1}+\rho _{2}+\rho _{3}\right) }{\left(
1-c_{2}^{3}\rho _{2}/c_{1}\right) \left( 1-e^{2}\rho _{3}\right) }\cdot
\left( 1-b_{1}^{2}/c_{1}c_{2}\right) ^{-\frac{1}{3}} \\
&&\cdot \int_{0}^{1}\mathrm{e}^{x\left( \sigma \right) }\left( x\left(
\sigma \right) +1\right) \mathrm{d}\sigma \cdot \mathrm{e}%
^{z_{1}b_{3}/c_{1}c_{2}b_{1}b_{2}}\cdot \mathrm{e}^{z_{2}d_{2}/ed_{1}},
\end{eqnarray*}%
where the double integral $\int \int $ is on the simplex $\left\{ \rho
_{1}+\rho _{2}+\rho _{3}=1\right\} $ and $x\left( \sigma \right) =\sigma
b_{2}/c_{1}c_{2}b_{1}+\left( 1-\sigma \right) d_{1}/e;$ in fact, the first
two integrals need regularizations, but we will do it later. The integral
over sigma is a result of the formula%
\begin{equation*}
\sum \frac{x^{k}y^{l}}{\left( k+l\right) !}=\sum_{p}\frac{h_{p}\left(
x,y\right) }{p!}=\frac{1}{x-y}\left( xe^{x}-ye^{y}\right)
=\int_{0}^{1}\left( xe^{x}\right) _{x=x\left( \sigma \right) }^{\prime }%
\mathrm{d}\sigma .
\end{equation*}%
Of course, we take the residuum of this expression with respect to
corresponding logarithmic form.

The residuum with respect to $c_{2}$ gives $c_{1}=\tau ^{1/4}$, with respect
to $c_{2}$ gives $c_{2}=\rho _{2}\tau ^{3/4}$ and with respect to $e$ gives $%
e=\rho _{3}^{-1/2}$ (compare the proof of Theorem 1). Thus we are left with
the residuum of%
\begin{eqnarray}
G\left( z_{1},z_{2}\right) &=&\frac{1}{\Gamma \left( -1/3\right) }%
\int_{0}^{1}\left( 1-\tau \right) ^{-7/3}\cdot \int \int \rho
_{1}^{-5/3}\cdot \mathrm{Res}\left( 1-\left( \rho _{2}\tau \right)
^{1/3}b_{1}^{2}\right) ^{-1/3}  \notag \\
&&\cdot \int_{0}^{1}\mathrm{e}^{x\left( \sigma \right) }\left( x\left(
\sigma \right) +1\right) \mathrm{d}\sigma \cdot \mathrm{e}^{\left( \rho
_{2}\tau \right) ^{1/3}z_{1}b_{3}/b_{1}b_{2}}\cdot \mathrm{e}^{\rho
_{3}^{1/2}z_{2}d_{2}/d_{1}}  \label{4.12}
\end{eqnarray}%
with respect to the form $\mathrm{d}^{3}\ln b\mathrm{d}^{2}\ln d;$ here%
\begin{equation*}
x\left( \sigma \right) =\sigma \left( \rho _{2}\tau \right)
^{1/3}b_{2}/b_{1}+\left( 1-\sigma \right) \rho _{3}^{1/2}d_{1}.
\end{equation*}%
The integral $\int_{0}^{1}\left( 1-\tau \right) ^{-7/3}\times \left( \cdot
\right) $ is regularized by taking $\left( 1-\mathrm{e}^{2\pi \mathrm{i}%
/3}\right) ^{-1}\int_{\gamma }\left( 1-\tau \right) ^{-7/3}\times \left(
\cdot \right) ,$ where $\gamma $ is a loop with vertex at $\tau =0$ and
surrounding $\tau =1$ (compare Remark 2.1). In the integral $\int \int \rho
_{1}^{-5/3}\times \left( \cdot \right) $ we have the measure $\int_{0}^{1}%
\mathrm{d}\rho _{2}\int_{0}^{\rho _{2}}\left( 1-\rho _{2}-\rho _{3}\right)
^{-5/3}\mathrm{d}\rho _{3};$ so, the integral of $\mathrm{d}\rho _{3}$ is
unchanged but the integral $\int_{0}^{1}\mathrm{d}\rho _{2}$ is replaced
with $\left( 1-\mathrm{e}^{-2\pi \mathrm{i}/3}\right) ^{-1}$ \linebreak $%
\times $(integral over a loop).

There is another way to get a convergent integral. We can reorganize the sum
in Eq. (4.10) using $\Gamma \left( \beta \right) =\Gamma \left( \beta
+1\right) /\beta ;$ we do not present details.\bigskip
\end{example}

\section{MZVs and new differential equations}

\subsection{Hypergeometric equations related with $\protect\zeta \left(
2\right) $}

This is the equation%
\begin{equation}
\left\{ \left( 1-t\right) \mathcal{D} ^{2}+\lambda ^{2}t\right\} u=0,
\label{5.1}
\end{equation}%
$\mathcal{D} =\mathcal{D} _{t}=t\frac{\mathrm{d}}{\mathrm{d}t}.$ Following
\cite{ZZ1} we look for its basic solutions near $t=0$ in form of series in
powers of $\lambda ^{2}:$%
\begin{equation*}
u_{j}\left( t;\lambda \right) =u_{j,0}\left( t\right) -u_{j,1}\left(
t\right) \lambda ^{2}+u_{j,2}\left( t\right) \lambda ^{4}-\ldots ,\text{ }%
j=0,1.
\end{equation*}%
The coefficient functions satisfy the recurrent equations%
\begin{equation*}
\mathcal{D} ^{2}u_{j,0}=0,\text{ }\left( 1-t\right) \mathcal{D}
^{2}u_{j,k+1}=tu_{j,k}.
\end{equation*}%
Thus we choose%
\begin{equation*}
u_{0,0}\left( t\right) =1,\text{ }u_{1,0}\left( t\right) =\ln \left( \lambda
^{2}t\right) =u_{0,0}\cdot \ln \lambda ^{2}+\ln t
\end{equation*}%
and the other equations are solved by%
\begin{equation*}
u_{j,k+1}\left( t\right) =\int_{0}^{t}\frac{\mathrm{d}t_{1}}{t_{1}}%
\int_{0}^{t_{1}}\frac{u_{j,k}\left( t_{2}\right) }{1-t_{2}}\mathrm{d}t_{2}.
\end{equation*}%
We get $u_{0,1}\left( t\right) =\sum_{n}\int_{0}^{t}\frac{\mathrm{d}t_{1}}{%
t_{1}}\int_{0}^{t_{1}}t_{2}^{n}\mathrm{d}t_{2}=\sum_{n}n^{-2}t^{n}=\mathrm{Li%
}_{2}\left( t\right) $ and, generally, $u_{0,k}\left( t\right) =\mathrm{Li}%
_{2,\ldots ,2}\left( t\right) $ (with $k$ 2's). (Here $\mathrm{Li}_{2,\ldots
,2}\left( t\right) $ is a special case of the polylogarithm is defined in
Note 1 above, in particular, $\mathrm{Li}_{d_{1},\ldots ,d_{k}}(1)=\zeta
\left( d_{1},\ldots ,d_{k}\right) ).$ Of course,%
\begin{equation}
u_{0}\left( t;\lambda \right) =F\left( \lambda ,-\lambda ;1;t\right)
\label{5.2}
\end{equation}%
is the hypergeometric function from Eq. (2.10) in Example 2.1 above.

Next, $u_{1,k}\left( t\right) =u_{0,k}\left( t\right) \ln \left( \lambda
^{2}t\right) -2(\mathrm{Li}_{3,2,\ldots ,2}(t)+\mathrm{Li}_{2,3,2,\ldots
,2}(t)+\ldots +\mathrm{Li}_{2,\ldots ,2,3}(t)).$ Thus the second solution is
of the form%
\begin{equation}
u_{1}\left( t;\lambda \right) =u_{0}\left( t;\lambda \right) \ln \left(
\lambda ^{2}t\right) +\tilde{u}_{1}\left( t;\lambda \right) ,  \label{5.3}
\end{equation}%
where $\tilde{u}_{1}$ is analytic near $t=0.$

From this it follows the following

\begin{lemma}
\label{l51}We have%
\begin{eqnarray}
u_{0}\left( 1;\lambda \right) &=&\Delta _{2}\left( \lambda \right) =1-\zeta
\left( 2\right) \lambda ^{2}+\zeta \left( 2,2\right) \lambda ^{4}-\ldots
=\prod_{n=1}^{\infty }\left( 1-\frac{\lambda ^{2}}{n^{2}}\right) ,
\label{5.4} \\
u_{1}\left( 1;\lambda \right) &=&2\Delta _{2}\left( \lambda \right) \left\{
\ln \lambda +\zeta \left( 3\right) \lambda ^{2}+\zeta \left( 5\right)
\lambda ^{4}+\zeta \left( 7\right) \lambda ^{6}+\ldots \right\} .
\label{5.5}
\end{eqnarray}
\end{lemma}

\textit{Proof}. The first statement is obvious. In \cite[Remark 2.1]{ZZ1} it
is claimed that the second formula follows from `simple resummation' of the
series in Eq. (5.4); in Section 5.3 below we give a direct proof. \qed%
\medskip

Functions (5.2) and (5.3) have other interesting properties. Firstly, they
undergo the following monodromy as $t$ turns around $0$ along a small circle:%
\begin{equation}
u_{0}\left( t;\lambda \right) \longmapsto u_{0}\left( \mathrm{e}^{2\pi
\mathrm{i}}t;\lambda \right) =u_{0}\left( t;\lambda \right) ,\text{ }%
u_{1}\left( t;\lambda \right) \longmapsto u_{1}\left( t;\lambda \right)
+2\pi \mathrm{i\cdot }u_{0}(t;\lambda ).  \label{5.6}
\end{equation}%
As functions of $\lambda ,$ they depend on $\lambda ^{2}$; moreover, $u_{0}$
and $\tilde{u}_{1}$ are entire. They also undergo a simple monodromy as $%
\lambda ^{2}$ turns around $0:$%
\begin{equation}
u_{0}\left( t;\lambda \right) \longmapsto u_{0}\left( t;\mathrm{e}^{\pi
\mathrm{i}}\lambda \right) =u_{0}\left( t;\lambda \right) ,\text{ }%
u_{1}\left( t;\lambda \right) \longmapsto u_{1}\left( t;\lambda \right)
+2\pi \mathrm{i}\cdot u_{0}(t;\lambda ).  \label{5.7}
\end{equation}

In \cite{ZZ1} the hypergeometric equation (5.1) was approximated with the
equation%
\begin{equation}
\mathcal{D}_{t}^{2}u=\lambda ^{2}tu,  \label{5.8}
\end{equation}%
or; with the variable $y=\lambda ^{2}t$, the Bessel type equation%
\begin{equation}
\mathcal{D}_{y}^{2}U=yU.  \label{5.9}
\end{equation}%
Eq. (5.9) has basic solutions%
\begin{equation}
U_{0}\left( y\right) =\sum \frac{\left( -y\right) ^{n}}{\left( n!\right) ^{2}%
}=F\left( \emptyset ;1;-y\right) =J_{0}\left( 2\sqrt{y}\right) ,\text{ }%
U_{1}\left( y\right) =U_{0}\left( y\right) \ln y+\widetilde{U}_{1}\left(
y\right) ,  \label{5.10}
\end{equation}%
where $\widetilde{U}_{1}$ is an entire function. The corresponding functions
$U_{0}\left( \lambda ^{2}t\right) $ and $U_{1}\left( \lambda ^{2}t\right) $
have the same monodromy around $t=0$ and around $\lambda ^{2}=0$ as the
functions $u_{0}$ and $u_{1}.$ Repeating the proof of the Wasow theorem
(Theorem 6 above), we get first part of the following

\begin{lemma}
\label{l52}The differential equations (5.1) and (5.8) are analytically
equivalent in a neighborhood of $t=\lambda =0.$

Moreover, the functions $u_{0}\left( t;\lambda \right) $ and $%
u_{1}(t;\lambda )$ satisfy a second order linear differential equation with
respect to $\lambda ;$ in particular, the quantities $u_{0}\left( 1;\lambda
\right) =\Delta _{2}\left( \lambda \right) $ and $u_{1}\left( 1;\lambda
\right) $ satisfy a second order linear differential equation with finite
singular points including $\lambda =0$ (which are regular) and $\lambda
=\infty $ (apparently irregular).
\end{lemma}

\textit{Proof}. The corresponding differential equation for $u=u(t;\cdot
)=u(t;\lambda )$ is following:%
\begin{equation}
\det \left(
\begin{array}{lll}
u & u_{0} & u_{1} \\
u^{\prime } & u_{0}^{\prime } & u_{1}^{\prime } \\
u^{\prime \prime } & u_{0}^{\prime \prime } & u_{1}^{\prime \prime }%
\end{array}%
\right) =0,  \label{5.11}
\end{equation}%
where $^{\prime }=\frac{\partial }{\partial \lambda }.$ Due to the monodromy
properties (5.7), the determinants, like $\left\vert
\begin{array}{ll}
u_{0}^{\prime } & u_{1}^{\prime } \\
u_{0}^{\prime \prime } & u_{1}^{\prime \prime }%
\end{array}%
\right\vert ,$ are single valued functions. \qed\medskip

One would like to ask about the solutions to Eq. (5.1) near the other
singular point $t=1.$ But it turns out that the corresponding differential
equation near $s=1-t=0$ has solutions $v_{0}\left( s\right) =\mathcal{D}
_{s}u_{0}\left( s\right) $ and $v_{1}\left( s\right) =\mathcal{D}
_{s}u_{1}\left( s\right) ,$ $\mathcal{D} _{s}=s\frac{\mathrm{d}}{\mathrm{d}s}
$ (see \cite{ZZ1}).

In \cite{ZZ1} the authors have applied the stationary phase formula to
integral (2.9) (in Example 2.1) to get another proof of the known formula%
\begin{equation}
\Delta _{2}\left( \lambda \right) =\frac{1}{\Gamma \left( 1-\lambda \right)
\Gamma \left( 1+\lambda \right) }=\frac{\sin \left( \pi \lambda \right) }{%
\pi \lambda }.  \label{5.12}
\end{equation}

We return to Eq. (5.11) in Section 5.3.

\subsection{Hypergeometric equation related with $\protect\zeta \left(
3\right) $}

Recall that it is%
\begin{equation}
\left\{ \left( 1-t\right) \mathcal{D} ^{3}+\lambda ^{3}t\right\} u=0,
\label{5.13}
\end{equation}%
$\mathcal{D} =\mathcal{D} _{t},$ i.e., Eq. (2.19) from Example 2.6 above.

We seek solutions $u_{j}\left( t;\lambda \right) ,$ $j=1,2,3,$ in the form
of series in powers of $\lambda ^{3}:$ $u_{j}\left( t;\lambda \right)
=u_{j,0}\left( t\right) -u_{j,1}(t)\lambda ^{3}+\ldots $. We get equations%
\begin{equation*}
\mathcal{D} ^{3}u_{j,0}=0,\text{ }\left( 1-t\right) \mathcal{D}
^{3}u_{j,k+1}=tu_{j,k}.
\end{equation*}

The following statement is a generalization of Lemma 5.1 (compare also the
beginning of Example 2.6).

\begin{lemma}
\label{l53}The basic solutions $u_{j}$ are such that%
\begin{eqnarray*}
u_{0,k}\left( t\right) &=&\mathrm{Li}_{3,\ldots ,3}\left( t\right) , \\
u_{1,k}\left( t\right) &=&u_{0,k}\left( t\right) \ln \left( \lambda
^{3}t\right) +\tilde{u}_{1,k}(t), \\
u_{2,k}\left( t\right) &=&\frac{1}{2}u_{0,k}\left( t\right) \ln ^{2}\left(
\lambda ^{3}t\right) +\tilde{u}_{1,k}(t)\ln (\lambda ^{3}t)+\tilde{u}%
_{2,k}(t),
\end{eqnarray*}%
where $\tilde{u}_{1,k}(t)=-3(\mathrm{Li}_{4,3,\ldots ,3}(t)+\mathrm{Li}%
_{3,4,3,\ldots ,3}(t)+\ldots +\mathrm{Li}_{3,\ldots ,3,4}(t))$ and $\tilde{u}%
_{2,k}(t)=6(\mathrm{Li}_{5,3,\ldots ,3}(t)+\mathrm{Li}_{3,5,3,\ldots
,3}(t)+\ldots +\mathrm{Li}_{3,\ldots ,3,5}(t)).$

In particular, we have%
\begin{eqnarray*}
u_{0}\left( 1;\lambda \right) &=&\Delta _{3}\left( \lambda \right) =1-\zeta
(3)\lambda ^{3}+\zeta (3,3)\lambda ^{6}-\ldots , \\
u_{1}(1;\lambda ) &=&3\Delta _{3}(\lambda )\left\{ \ln \lambda +\zeta
(4)\lambda ^{3}+\zeta (7)\lambda ^{6}+\zeta (10)\lambda ^{9}+\ldots \right\}
, \\
u_{2}(1;\lambda ) &=&\Delta _{3}(\lambda )\{\frac{1}{2}\ln ^{2}\lambda
^{3}+3\lambda ^{3}\ln \lambda ^{3}\left[ \zeta (4)+3\zeta (7)\lambda
^{3}+\ldots \right] \\
&&-6\lambda ^{3}\left[ \zeta (5)+\zeta (8)\lambda ^{3}+\zeta (11)\lambda
^{6}+\ldots \right] \}.
\end{eqnarray*}
\end{lemma}

Therefore we have the following basic solutions neat $t=0:$%
\begin{eqnarray}
u_{0} &=&F(-\lambda ,\epsilon \lambda ,\bar{\epsilon}\lambda ;1,1;t),
\label{5.14} \\
u_{1} &=&u_{0}\ln \left( \lambda ^{3}t\right) +\tilde{u}_{1}(t;\lambda ),%
\text{ }u_{2}=\frac{1}{2}u_{0}\ln ^{2}(\lambda ^{3}t)+\tilde{u}\ln (\lambda
^{3}t)+\tilde{u}_{2}(t;\lambda ).  \notag
\end{eqnarray}%
As functions of $\lambda ,$ they depend on $\lambda ^{3}$; moreover, $u_{0},$
$\tilde{u}_{1}$ and $\tilde{u}_{2}$ are entire. They also undergo a simple
monodromy as $\lambda ^{2}$ turns around $0$, a suitable generalization of
Eqs. (5.7) (which we do not write down).

In \cite{ZZ3} the third order hypergeometric equation (5.13) was
approximated with the equation%
\begin{equation}
\mathcal{D}_{t}^{3}u=\lambda ^{3}tu,  \label{5.15}
\end{equation}%
or; with the variable%
\begin{equation*}
y=\lambda ^{3}t,
\end{equation*}%
the Bessel type equation%
\begin{equation}
\mathcal{D}_{y}^{3}U=yU.  \label{5.16}
\end{equation}%
Eq. (5.16) has basic solutions%
\begin{eqnarray}
U_{0} &=&\sum \frac{\left( -y\right) ^{n}}{\left( n!\right) ^{3}}=F\left(
\emptyset ;1,1;-y\right) ,  \label{5.17} \\
U_{1} &=&U_{0}\ln y+\widetilde{U}_{1}\left( y\right) ,\text{ }  \notag \\
U_{2} &=&\frac{1}{2}U_{0}\ln ^{2}y+\widetilde{U}_{1}\ln y+\widetilde{U}%
_{2}(y),  \notag
\end{eqnarray}%
where $\widetilde{U}_{1,2}$ are entire functions. The corresponding
functions $U_{j}\left( \lambda ^{3}t\right) $, $j=0,1,2,$ have the same
monodromy around $t=0$ and around $\lambda ^{3}=0$ as the functions $u_{0}$,
$u_{1}$ and $u_{2}.$ Repeating the proof of the Wasow theorem (see Theorem 6
above and \cite[Proposition 3.1]{ZZ1}), we get the first part of the
following

\begin{lemma}
\label{l54}The differential equations (5.13) and (5.16) are analytically
equivalent in a neighborhood of $t=\lambda =0.$

Moreover, the functions $u_{j}\left( t;\lambda \right) $, $j=0,1,2$, satisfy
a third order linear differential equation with respect to $\lambda ;$ in
particular, the quantities $u_{0}\left( 1;\lambda \right) =\Delta _{3}\left(
\lambda \right) ,$ $u_{1}\left( 1;\lambda \right) $ and $u_{2}(1;\lambda )$
satisfy a third order linear differential equation with finite regular
singular points (including $\lambda =0)$ and $\lambda =\infty $ (apparently
irregular).
\end{lemma}

Recall also the following analogue of Eq. (5.12):%
\begin{equation}
\Delta _{3}(\lambda )=\prod_{n=1}^{\infty }\left( 1-\frac{\lambda ^{3}}{n^{3}%
}\right) =\frac{1}{\Gamma (1+\lambda )\Gamma (1-\epsilon \lambda )\Gamma (1-%
\bar{\epsilon}\lambda )}.  \label{5.18}
\end{equation}

Now we recall the analysis of basic solutions near $s=1-t=0,$ presented in
Example 2.6. The hypergeometric equation (5.13) takes the form (2.20), i.e.,%
\begin{equation}
\left\{ \mathcal{D} _{s}\left[ \left( 1-s\right) \mathcal{D} _{s}\right]
^{2}-\lambda ^{3}\right\} v=0,  \label{5.19}
\end{equation}%
$\mathcal{D} _{s}=s\frac{\partial }{\partial s},$ and the basic solutions $%
v_{j}(s;\lambda )$ are of the form%
\begin{equation*}
v_{1}=\lambda ^{3/2}s+\ldots ,\text{ }v_{2}=\lambda ^{3}s^{2}+\ldots ,\text{
}v_{3}=\frac{1}{4}v_{2}\ln (\lambda ^{3}s^{2})+1+w(s;\lambda ).
\end{equation*}%
where $v_{1},$ $v_{2}$ and $w$ are analytic (compare Eqs. 2.21)).

With the variable%
\begin{equation*}
z=\lambda ^{3}s^{2}
\end{equation*}%
(where $\lambda $ is large but $z$ is still small), one gets approximations%
\begin{equation*}
v_{j}(s;\lambda )=V_{j}(z)+O(\lambda ^{-3/2}),
\end{equation*}%
where $V_{j}$ are basic solutions to the completely confluent hypergeometric
equation%
\begin{equation}
\left\{ 8\mathcal{D}_{z}(\mathcal{D}_{z}-1/2)(\mathcal{D}_{z}-1)-z\right\}
V=0  \label{5.20}
\end{equation}%
(compare Eqs. (2.22)--(2.23)). We have the following analogue of Lemma 5.4.

\begin{lemma}
\label{l55}The differential equations (5.19) and (5.20) are analytically
equivalent in a neighborhood of $s=\lambda =0.$
\end{lemma}

In our paper \cite[Theorem 5.7]{ZZ3} it was stated that the generating
function $\Delta _{3}(\lambda )$ (denoted there $f_{3}(x))$ satisfies a
sixth order linear differential equation near $\lambda =\infty $ and that
admits a WKB type expansion there.

By Lemma 5.4 above we know that this function satisfies a third order linear
differential equation for with regular singularity at $\lambda =0$ (and
plausibly infinitely many other regular finite singularities); so, in this
aspect, the situation is simpler. But the statement about the WKB type
expansion is definitely wrong and we are going to prove it.

In \cite{ZZ3} we have wrongly assumed that the equivalences from Lemmas 5.4
and 5.5 work not only for small $\lambda ,$ but also for large $y=\lambda
^{3}t$ (respectively, for large $z=\lambda ^{3}s^{2}=\lambda ^{3}(1-t)^{2}).$
As we have mentioned in the Introduction, we were mislead by the Wasow's
approach to the perturbations of the Airy equation (see Section 4 above). We
have assumed that the Stokes phenomena for Eq. (5.16) (respectively, (5.19))
are the same as the Stokes phenomena for the corresponding hypergeometric
equation (5.13) (respectively, (5.20)).

Note that the generating function $\Delta _{3}(\lambda )$ equals the
connection coefficient $C(\lambda )$ in the representation%
\begin{equation*}
u_{0}(t;\lambda )=A(\lambda )v_{1}(1-t;\lambda )+B(\lambda
)v_{2}(1-t;\lambda )+C(\lambda )v_{3}(1-t;\lambda )
\end{equation*}%
(because $v_{1}(0;\lambda )=v_{2}(0;\lambda )=0$ and $v_{3}(0;\lambda )=1).$
So, one expands the solution $u_{0}$ for large $y=\lambda ^{3}t$ in WKB type
solutions near $t=0,$ using the mountain pass integral form Eq. (2.11).
Next, one expands the basic solutions $v_{j}$ for large $z=\lambda ^{3}s^{2}$
in some WKB type solutions, using the expansions of the functions $V_{j}(z)$
from Example 2.6 and of their variations. The formal WKB series are
relatively simple; so, one can hope to get the connecting coefficient $%
C(\lambda )$ by evaluating some integrals.

In \cite{ZZ3} it was claimed that $C(\lambda )=\Delta _{3}(\lambda )$ is
expanded in the following WKB type functions%
\begin{equation}
\Delta ^{\left( \sigma \right) }(\lambda )=\lambda ^{-3/2}\mathrm{e}^{2\pi
\sigma \lambda /\sqrt{3}}\omega ^{\left( \sigma \right) }(\lambda ^{-1/2}),%
\text{ }\sigma =-1,\epsilon ,\bar{\epsilon},  \label{5.21}
\end{equation}%
as $\lambda \rightarrow \infty ,$ where $\omega ^{\left( \sigma \right)
}=1+\ldots $ is a formal power series; moreover, these expressions would be
subject to a concrete Stokes phenomenon. To explain this formula we recall
the formula $u_{0}\sim (2\pi \sqrt{3}\lambda
)^{-1}((1-t)/t)^{1/3}\sum_{\sigma =-1,\epsilon ,\bar{\epsilon}}\bar{\sigma}%
\mathrm{e}^{\sigma \tilde{S}(t)},$ $\tilde{S}=\int_{0}^{t}\tau
^{-2/3}(1-\tau )^{-1/3}\mathrm{d}\tau ,$ from Example 2.9. Near $t=0$ we
have $\lambda \tilde{S}(t)\approx 3\lambda t^{1/3}=3y^{1/3}$ and $u_{0}\sim
(2\pi \sqrt{3})^{-1}y^{-1/3}\sum \bar{\sigma}\mathrm{e}^{3\sigma y^{1/3}}.$
Near $s=1-t=0$ we have $\lambda \tilde{S}(t)\approx 2\pi /\sqrt{3}-\frac{3}{2%
}\lambda s^{2/3}=2\pi /\sqrt{3}-\frac{3}{2}z^{1/3}$ and the WKB solutions $%
\sim \frac{\mathrm{i}}{2}\sqrt{2\pi /3}z^{1/6}\mathrm{e}^{-\frac{3}{2}\sigma
z^{1/3}}$ to the corresponding confluent equation (i.e., for $V_{j}(z))$;
compare Example 2.8.

Unfortunately, the Wasow type approach does not work in the cases of a
genuine Stokes phenomena and our hopes to get a significant progress in
understanding the nature of the function $\Delta _{3}(\lambda )$ have turned
out futile. Instead, we have the following result which confirms this.

\begin{proposition}
\label{p51}It is impossible to expand the function $\Delta _{3}(\lambda )$
in functions like in Eq. (5.21) in any section near $\lambda =\infty .$
\end{proposition}

\textit{Proof}. We use the identity%
\begin{equation}
\Delta _{3}(\lambda )\Delta _{3}(-\lambda )=\prod \left( 1-\frac{\lambda ^{6}%
}{n^{6}}\right) =\Delta _{2}(\lambda )\Delta _{2}(-\epsilon \lambda )\Delta
_{2}(-\bar{\epsilon}\lambda ),  \label{5.22}
\end{equation}%
where $\Delta _{2}(\lambda )=\frac{1}{2\pi \mathrm{i}\lambda }\left( \mathrm{%
e}^{\pi \mathrm{i}\lambda }-\mathrm{e}^{-\pi \mathrm{i}\lambda }\right) .$

The right-hand side equals%
\begin{eqnarray}
&\sim &\frac{\mathrm{i}}{(2\pi \lambda )^{3}}\{\left( \mathrm{e}^{2\pi
\mathrm{i}\lambda }-\mathrm{e}^{-2\pi \mathrm{i}\lambda }\right) -\left(
\mathrm{e}^{\left( -\sqrt{3}+\mathrm{i}\right) \pi \lambda }-\mathrm{e}%
^{-\left( -\sqrt{3}+\mathrm{i}\right) \pi \lambda }\right)  \label{5.23} \\
&&-\left( \mathrm{e}^{\left( \sqrt{3}+\mathrm{i}\right) \pi \lambda }-%
\mathrm{e}^{-\left( \sqrt{3}+\mathrm{i}\right) \pi \lambda }\right) \}
\notag
\end{eqnarray}%
and is single valued (no Stokes phenomena).

Suppose that, in some sector near $\lambda =\infty ,$ we have the
representations%
\begin{eqnarray*}
\Delta _{3}(\lambda ) &=&a\Delta ^{\left( -1\right) }(\lambda )+b\Delta
^{\left( \epsilon \right) }(\lambda )+c\Delta ^{\left( \bar{\epsilon}\right)
}(\lambda ) \\
&\sim &\lambda ^{-3/2}\left\{ a\mathrm{e}^{-2\pi \lambda /\sqrt{3}}+b\mathrm{%
e}^{\pi (1+\mathrm{i}\sqrt{3})\lambda /\sqrt{3}}+c\mathrm{e}^{\pi (1-\mathrm{%
i}\sqrt{3})\lambda /\sqrt{3}}\right\} , \\
\Delta _{3}(-\lambda ) &=&\alpha \Delta ^{\left( -1\right) }(-\lambda
)+\beta \Delta ^{\left( \epsilon \right) }(-\lambda )+\gamma \Delta ^{\left(
\bar{\epsilon}\right) }(-\lambda ) \\
&\sim &\mathrm{i}\lambda ^{-3/2}\left\{ \alpha \mathrm{e}^{2\pi \lambda /%
\sqrt{3}}+\beta \mathrm{e}^{-\pi (1+\mathrm{i}\sqrt{3})\lambda /\sqrt{3}%
}+\gamma \mathrm{e}^{-\pi (1-\mathrm{i}\sqrt{3})\lambda /\sqrt{3}}\right\} ,
\end{eqnarray*}%
where the coefficients depend on the sector. We get%
\begin{eqnarray}
-\mathrm{i}\lambda ^{3}\Delta _{3}(\lambda )\Delta _{3}(-\lambda ) &\sim
&\left( a\alpha +b\beta +c\gamma \right) +\left( c\beta e^{2\pi \mathrm{i}%
\lambda }+b\gamma e^{-2\pi \mathrm{i}\lambda }\right)   \label{5.24} \\
&&+\left( a\gamma e^{\left( -\sqrt{3}+\mathrm{i}\right) \pi \lambda
}+c\alpha e^{-\left( -\sqrt{3}+\mathrm{i}\right) \pi \lambda }+\right)
\notag \\
&&+\left( b\alpha e^{\left( \sqrt{3}+\mathrm{i}\right) \pi \lambda }+a\beta
e^{-\left( \sqrt{3}+\mathrm{i}\right) \pi \lambda }\right) .  \notag
\end{eqnarray}%
Comparing Eqs. (5.23)--(5.24) we get the conditions%
\begin{equation*}
a\alpha +b\beta +c\gamma =0,\text{ }c\beta +b\gamma =0,\text{ }b\alpha
+a\beta =0,\text{ }c\alpha +a\gamma =0
\end{equation*}%
and all the coefficients are nonzero. But we get $\beta =-\frac{b}{a}\alpha ,
$ $\gamma =-\frac{c}{a}\alpha $ and%
\begin{equation*}
\gamma =-\frac{c}{b}\beta =-\frac{c}{b}\left( -\frac{b}{a}\alpha \right) =%
\frac{c}{a}\alpha ,
\end{equation*}%
which leads to a contradiction. \qed

\subsection{MZVs and Gamma functions}

Here we consider the second order linear differential equation (5.11)
satisfied by the generating function $u_{0}(1;\lambda )=\Delta _{2}(\lambda
) $ \ and by $u_{1}(1;\lambda )$ (see Lemma 5.2). It takes the form%
\begin{equation}
A(\lambda )\Phi ^{\prime \prime }+B(\lambda )\Phi ^{\prime }+C(\lambda )\Phi
=0,  \label{5.25}
\end{equation}%
where $^{\prime }=\frac{\mathrm{d}}{\mathrm{d}\lambda }$ and $A(\lambda ),$ $%
B(\lambda )$ and $C(\lambda )$ are entire functions. Here we have replaced $%
u(1;\lambda )$ in Eq. (5.11) with $\Phi (\lambda )$ and the coefficients $%
A,B,C$ are defined by corresponding minors of the matrix in the same
equation.

Of course, $\lambda =0$ is a regular singular point; but there can be other
finite regular singular points.

The aim of this section is to look more closely at the basic solutions to
Eq. (5.25) and, in particular, at the Stokes-like phenomenon at $\lambda
=\infty .$ To this aim we will use the representation of these solutions via
the Euler Gamma functions.

Recall that one solution equals%
\begin{equation}
\Phi _{1}(\lambda )=\Delta _{2}(\lambda )=\frac{1}{\Gamma (1+\lambda )\Gamma
(1-\lambda )}=\frac{\sin (\pi \lambda )}{\pi \lambda }.  \label{5.26}
\end{equation}

\begin{lemma}
\label{l56}The second solution equals%
\begin{equation}
\Phi _{2}(\lambda )=u_{1}(1;\lambda )=\frac{\sin (\pi \lambda )}{\pi \lambda
}\left\{ \ln \lambda ^{2}-2\gamma -\Psi (1+\lambda )-\Psi (1-\lambda
)\right\} ,  \label{5.27}
\end{equation}%
where $\Psi (z)=\Gamma ^{\prime }(z)/\Gamma (z)$ is the logarithmic
derivative of the Gamma function and $\gamma =-\Psi (1)$ is the
Euler--Mascheroni constant.
\end{lemma}

\textit{Proof}. Recall that we evaluate at $t=1$ the second solution to the
hypergeometric equation (5.1); the first solution is the hypergeometric
series.

The second solution is of the form $u_{1}=u_{0}\ln (\lambda ^{2}t)+\tilde{u}%
_{1}(t;\lambda ).$ It is calculated by, firstly, perturbing the
hypergeometric equation (5.1) to the hypergeometric equation%
\begin{equation*}
\left\{ \left( 1-t\right) \mathcal{D} _{t}^{2}+\left( t\lambda ^{2}-\mu
^{2}\right) \right\} u=\left\{ \left( \mathcal{D} _{t}^{2}-\mu ^{2}\right)
-t\left( \mathcal{D} _{t}^{2}-\lambda ^{2}\right) \right\} u=0
\end{equation*}%
and then passing to a corresponding limit as $\mu \rightarrow 0.$ The latter
equation has two solutions:%
\begin{equation*}
\eta _{\mu }=t^{\mu }F\left( \mu +\lambda ,\mu -\lambda ;1+2\mu ;t\right)
\end{equation*}%
and $\eta _{-\mu }$. One finds%
\begin{equation*}
\eta _{\mu }=u_{0}+\mu u_{0}\ln t+\mu \tilde{u}_{1}\left( t;\lambda \right)
+O\left( \mu ^{2}\right) .
\end{equation*}

Below we put $t=1$ and use the known formula for the value of the
hypergeometric function at $t=1$ (see \cite{BE1} and Eq. (2.9) above): $%
F(\alpha _{1},\alpha _{2};\beta ;1)=\Gamma (\beta )\Gamma (\beta -\alpha
_{1}-\alpha _{2})/\Gamma (\beta -\alpha _{1})\Gamma (\beta -\alpha _{2}).$
Thus%
\begin{equation*}
F\left( \mu +\lambda ,\mu -\lambda ;1+2\mu ;1\right) =\frac{\Gamma \left(
1+2\mu \right) }{\Gamma \left( 1+\lambda +\mu \right) \Gamma \left(
1-\lambda +\mu \right) }.
\end{equation*}

Next, we have%
\begin{equation*}
\Gamma (x+\mu )=\Gamma (x)\left( 1+\mu \frac{\Gamma ^{\prime }(x)}{\Gamma (x)%
}+\ldots \right) =\Gamma (x)\left( 1+\mu \Psi (x)+\ldots \right) ,
\end{equation*}%
where $\Psi (x)$ is the logarithmic derivative of the Gamma function.

Hence%
\begin{equation*}
\tilde{u}_{1}\left( 1;\lambda \right) =u_{0}(1;\lambda )\left\{ 2\Psi
(1)-\Psi (1+\lambda )-\Psi (1-\lambda )\right\} ,
\end{equation*}%
which give the formula (5.27). \qed\bigskip

Using the formula (see \cite[Eq. 1.17(5)]{BE1})%
\begin{equation}
\Psi (1+z)=-\gamma +\zeta (2)z-\zeta (3)z^{2}+\zeta (4)z^{3}-\ldots \text{,
as }\lambda \rightarrow 0,  \label{5.28}
\end{equation}%
we get Eq. (5.5) from Lemma 5.1 above.\footnote{%
From the formulas in Lemma 5.3 above evaluating at $t=1$ the basic solution
of the third order hypergeometric equation one can guess that also in that
case one can express the basic solutions via $\Delta _{3}(\lambda )$ and the
Euler Psi-function. Probably, the values at $t=1$ of the hypergeometric
functions $F(\alpha _{1},\alpha _{2},\alpha _{3};\beta _{1},\beta _{2};t)$
are expressed via a product of Gamma functions. We did not succeed in
proving a corresponding formula.}\medskip

Now we consider the case with $\left\vert \lambda \right\vert \rightarrow
\infty .$ Here we use the following formulas:%
\begin{equation}
\Psi (1+z)-\Psi (1-z)=\frac{1}{z}-\pi \tan ^{-1}(\pi z)  \label{5.29}
\end{equation}%
(which follows from the formula for $\Gamma (1+z)\Gamma (1-z),$ see \cite[%
Section 1.7.1]{BE1}) and%
\begin{equation}
\Psi (z)=\ln z-\frac{1}{2z}-\sum_{n\geq 1}\frac{B_{2n}}{2n}\frac{1}{z^{2n}},%
\text{ as }\left\vert z\right\vert \rightarrow \infty \text{ and }\left\vert
\arg z\right\vert <\pi ,  \label{5.30}
\end{equation}%
where $B_{2n}$ are the Bernoulli numbers (which follows from the stationary
phase formula for the mountain pass integral defining the Gamma function or
the Stirling formula, see \cite[Eq. 1.18(7)]{BE1}). Thus%
\begin{eqnarray*}
\Psi (1+z) &\approx &\ln z+\frac{1}{2z}+\left\{ \ln \left( 1+\frac{1}{z}%
\right) -\frac{1}{z}+\frac{1}{2z(z+1)}-\sum \frac{B_{2n}}{2n}\frac{1}{%
(z+1)^{2n}}\right\} \\
&=&\ln z+\frac{1}{2z}+\Omega \left( \frac{1}{z}\right) ,
\end{eqnarray*}%
where $\Omega =O\left( \frac{1}{z^{2}}\right) ,$ and%
\begin{equation*}
\Psi (1-z)=\Psi \left( 1+\left( -z\right) \right) =-\mathrm{i}\pi +\ln (z)-%
\frac{1}{2z}+\Omega \left( -\frac{1}{z}\right) ,
\end{equation*}%
as $\left\vert z\right\vert \rightarrow \infty $ and $\pi <\arg z<2\pi .$

For $\left\vert z\right\vert \rightarrow \infty $ and $\left\vert \arg
z\right\vert <\pi $ we get%
\begin{eqnarray}
\Phi _{2}(\lambda ) &=&\frac{\sin (\pi \lambda )}{\pi \lambda }\left\{ 2\ln
\lambda -2\gamma -2\Psi (1+\lambda )+\frac{1}{\lambda }-\pi \tan ^{-1}(\pi
\lambda )\right\}  \notag \\
&=&2\frac{\sin (\pi \lambda )}{\pi \lambda }\left\{ -\gamma -\Omega \left(
\frac{1}{z}\right) \right\} -\frac{\cos (\pi \lambda )}{\lambda }.
\label{5.31}
\end{eqnarray}

For $\left\vert z\right\vert \rightarrow \infty $ and $\pi <\arg z<2\pi $ we
get%
\begin{eqnarray}
\Phi _{2}(\lambda ) &=&\frac{\sin (\pi \lambda )}{\pi \lambda }\left\{ 2\ln
\lambda -2\gamma -2\Psi (1-\lambda )-\frac{1}{\lambda }+\pi \tan ^{-1}(\pi
\lambda )\right\}  \notag \\
&=&2\frac{\sin (\pi \lambda )}{\pi \lambda }\left\{ -\gamma -\Omega \left( -%
\frac{1}{z}\right) -\mathrm{i}\pi \right\} +\frac{\cos (\pi \lambda )}{%
\lambda }.  \label{5.32}
\end{eqnarray}

We observe a kind of Stokes type phenomenon. It seems to be strange.
However, Eq. (5.25) has infinite number of finite regular singular points.
Indeed, by Eq. (5.11) the coefficient $A(\lambda )$ in Eq. (5.25) equals%
\begin{equation*}
A=\det \left(
\begin{array}{ll}
\Phi _{1} & \Phi _{2} \\
\Phi _{1}^{\prime } & \Phi _{2}^{\prime }%
\end{array}%
\right) =\Phi _{1}^{2}\Theta ,
\end{equation*}%
where the function $\Theta (\lambda )=\frac{\mathrm{d}}{\mathrm{d}\lambda }%
\left\{ \ln \lambda ^{2}-2\gamma -\Psi (1+\lambda )-\Psi (1-\lambda
)\right\} $ has infinitely many poles (compensating zeroes of $\Phi _{1}$ at
nonzero integers) and apparently it has also infinitely many poles.

Thus $\lambda =\infty $ is an accumulation point of finite regular singular
points of Eq. (2.25). This a somewhat new phenomenon.\bigskip

Finally we note that the generating function $\Delta _{3}(\lambda )=\frac{1}{%
\Gamma (1+\lambda )\Gamma (1-\epsilon \lambda )\Gamma (1-\bar{\epsilon}%
\lambda )}$ experiences very complicated Stokes type behavior as $\lambda
\rightarrow \infty $; i.e., when we expand the factors $\Gamma (1+\sigma
\lambda )$ (when $\mathrm{Re}(\sigma \lambda )\geq 0)$ or the factors $%
\Gamma (1+\sigma \lambda )=\frac{\pi \sigma \lambda }{\sin (\pi \sigma
\lambda )}\Gamma ^{-1}(1-\sigma \lambda )$ (where $\mathrm{Re}(\sigma
\lambda )\leq 0)$ using the Stirling formula. Plausibly, it is caused by by
the property of accumulation of finite regular singular points of the
corresponding third order equation.

We think that these equations deserve further investigations.

\section{Appendices}

\subsection{Proof of Lemma 2.1}

Substituting $\theta _{1}=-\theta _{2}-\ldots -\theta _{q+1}$ we get the form%
\begin{equation*}
\sum_{j\geq 2}e_{1}\left( \lambda _{1},\lambda _{j}\right) \theta _{j}^{2}+2%
\frac{\lambda _{1}}{e_{0}}\sum_{j<k}\theta _{j}\theta _{k},
\end{equation*}%
where $e_{0}=1.$ By induction we will prove that the form equals%
\begin{eqnarray*}
&&e_{1}\left( \lambda _{1},\lambda _{2}\right) \tilde{\theta}_{2}^{2}+\frac{%
e_{2}\left( \lambda _{1},\ldots ,\lambda _{3}\right) }{e_{1}\left( \lambda
_{1},\lambda _{2}\right) }\tilde{\theta}_{3}^{2}+\ldots +\frac{e_{k-1}\left(
\lambda _{1},\ldots ,\lambda _{k}\right) }{e_{k-2}\left( \lambda _{1},\ldots
,\lambda _{k-1}\right) }\tilde{\theta}_{k}^{2} \\
&&+\sum_{j>k}\frac{e_{k}\left( \lambda _{1},\ldots ,\lambda _{k},\lambda
_{j}\right) }{e_{k-1}\left( \lambda _{1},\ldots ,\lambda _{k}\right) }\theta
_{j}^{2}+2\frac{\lambda _{1}\cdots \lambda _{k}}{e_{k-1}\left( \lambda
_{1},\ldots ,\lambda _{k}\right) }\sum_{k<i<j}\theta _{i}\theta _{j},
\end{eqnarray*}%
where $\tilde{\theta}_{j}=\theta _{j}+\ldots $.

We put%
\begin{equation*}
\tilde{\theta}_{k+1}=\theta _{k+1}+\frac{\lambda _{1}\cdots \lambda _{k}}{%
e_{k}\left( \lambda _{1},\ldots ,\lambda _{k+1}\right) }\sum_{j>k+1}\theta
_{j}
\end{equation*}%
and get the form with initial terms with $\tilde{\theta}_{i}^{2},$ $i\leq
k+1,$ like above and with%
\begin{equation*}
\frac{1}{e_{k-1}\left( \lambda _{1},\ldots ,\lambda _{k}\right) }\left(
e_{k}\left( \lambda _{1},\ldots ,\lambda _{k},\lambda _{j}\right) -\frac{%
\left( \lambda _{1}\cdots \lambda _{k}\right) ^{2}}{e_{k}\left( \lambda
_{1},\ldots ,\lambda _{k+1}\right) }\right) \theta _{j}^{2},\text{ }j>k+1,
\end{equation*}%
\begin{equation*}
\frac{2\cdot \lambda _{1}\cdots \lambda _{k}}{e_{k-1}\left( \lambda
_{1},\ldots ,\lambda _{k}\right) }\left( 1-\frac{\left( \lambda _{1}\cdots
\lambda _{k}\right) ^{2}}{e_{k}\left( \lambda _{1},\ldots ,\lambda
_{k+1}\right) }\right) \theta _{i}\theta _{j},\text{ }k+1<i<j.
\end{equation*}%
The induction step follows from the following relations between elementary
symmetric polynomials:%
\begin{eqnarray*}
&&e_{k}\left( \lambda _{1},\ldots ,\lambda _{k},\lambda _{j}\right)
e_{k}\left( \lambda _{1},\ldots ,\lambda _{k+1}\right) -\left( \lambda
_{1}\cdots \lambda _{k}\right) ^{2} \\
&=&e_{k-1}\left( \lambda _{1},\ldots ,\lambda _{k}\right) e_{k+1}\left(
\lambda _{1},\ldots ,\lambda _{k+1},\lambda _{j}\right) ,
\end{eqnarray*}%
\begin{equation*}
e_{k}\left( \lambda _{1},\ldots ,\lambda _{k+1}\right) -\lambda _{1}\cdots
\lambda _{k}=\lambda _{k+1}e_{k-1}\left( \lambda _{1},\ldots ,\lambda
_{k}\right) .
\end{equation*}

To justify the first relation we note that in left-hand side we have
monomials: either of the form (i) $\lambda ^{I}\lambda ^{J}\lambda
_{j}\times \lambda ^{I}\lambda ^{K}\lambda _{k+1},$where the sets $%
I,J,K\subset A=\left\{ 1,\ldots ,k\right\} $ are disjoint; or of the form
(ii) $\lambda ^{A}\times \lambda ^{K}\lambda _{k+1}$; or of the form (iii) $%
\lambda ^{J}\lambda _{j}\times \lambda ^{A}$. On the right-hand side we
have: either the monomials (i') $\lambda ^{I}\lambda ^{L}\times \lambda
^{I}\lambda ^{M}\lambda _{k+1}\lambda _{j},$ where $I,L,M\subset A$ are
disjoint; or (ii') $\lambda ^{A}\times \lambda ^{L}\lambda _{k+1};$ or
(iii') $\lambda ^{A}\times \lambda ^{M}\lambda _{j}.$

In the cases (i) and (i') we deal with divisions $J\cup K=L\cup M$ of the
same subset of $A$ into two disjoint subsets of the same cardinality. It is
also clear that the terms (ii) and (ii') (respectively, (iii) and (iii'))
correspond one to another.

The second relation is quite obvious. \qed

\subsection{Remark 2.5 continued}

In this case the phase equals%
\begin{equation*}
\phi =-\sum_{1}^{p}\alpha _{j}\ln \left( 1-a_{j}\eta \right) +\eta t^{\kappa
}\sum_{k=p+1}^{q+1}a_{k},
\end{equation*}%
$\kappa =\frac{1}{q+1-p},$ where $a_{i}$'s are subject to the restriction $%
\varphi =a_{1}\cdots a_{q+1}=1.$ The critical points of the phase are given
by the equations%
\begin{eqnarray*}
\left( \phi -\rho \varphi \right) _{a_{j}}^{\prime } &=&\frac{\eta \alpha
_{j}}{1-a_{j}\eta }-\frac{\rho }{a_{j}}=0,\text{ }j=1,\ldots ,p, \\
\left( \phi -\rho \varphi \right) _{a_{k}}^{\prime } &=&\eta t^{\kappa }-%
\frac{\rho }{a_{k}}=0,\text{ }k=p+1,\ldots ,q+1,
\end{eqnarray*}%
where $\rho $ is the Lagrange multiplier. We find%
\begin{eqnarray*}
a_{j} &=&\frac{\rho }{\eta \left( \rho +\alpha _{j}\right) },\text{ }%
j=1,\ldots ,p, \\
a_{k} &=&\rho /\eta t^{\kappa },\text{ }k=p+1,\ldots ,q+1,
\end{eqnarray*}%
and hence the multiplier $\rho $ satisfies the equation%
\begin{equation*}
\left( \frac{\rho }{\eta }\right) ^{q+1}=t\prod\limits_{1}^{p}\left( \rho
+\alpha _{j}\right) .
\end{equation*}%
As $\left\vert t\right\vert \rightarrow \infty $ we have approximate
solutions for $\rho $ of two sorts:%
\begin{eqnarray*}
\rho ^{\left( l\right) } &\approx &-\alpha _{l}+b_{l}t^{-1},\text{ }%
l=1,\ldots ,p, \\
\rho ^{\left( m\right) } &\approx &\zeta ^{m-p}\eta ^{\left( q+1\right)
\kappa }t^{\kappa },\text{ }m=p+1,\ldots ,q+1,
\end{eqnarray*}%
where $b_{l}=\left( -\alpha _{l}/\eta \right)
^{q+1}/\prod\limits_{n\not=l}\left( \alpha _{n}-\alpha _{l}\right) ,$ $\zeta
=\mathrm{e}^{2\pi \mathrm{i}\kappa }.$

Let us expand the phase near a critical point corresponding to $\rho =\rho
^{\left( l\right) }\approx -\alpha _{l},$ $l=1,\ldots ,p.$ Putting $%
a_{l}=a_{l}^{\left( l\right) }\mathrm{e}^{\mathrm{i}\theta _{l}}\approx
\left( -\alpha _{l}t/\eta b_{l}\right) \mathrm{e}^{\mathrm{i}\theta _{l}}$, $%
a_{j}=a_{j}^{\left( l\right) }\mathrm{e}^{\mathrm{i}\theta _{j}}\approx
\left( -\alpha _{l}/\eta \left( \alpha _{j}-\alpha _{l}\right) \right)
\mathrm{e}^{\mathrm{i}\theta _{j}}$ for $j\not=l,$ and $a_{k}=a_{k}^{\left(
l\right) }\mathrm{e}^{\mathrm{i}\theta _{k}}\approx \left( -\alpha
_{l}t^{\kappa }/\eta \right) \mathrm{e}^{\mathrm{i}\theta _{k}}$ we get the
phase%
\begin{equation*}
\phi =-\alpha _{l}\ln t+\mathrm{const}-\mathrm{i}\alpha _{l}\theta
_{l}-\sum_{j\not=l}\alpha _{j}\ln \left( 1+\frac{\alpha _{l}}{\alpha
_{j}-\alpha _{l}}\mathrm{e}^{\mathrm{i}\theta _{j}}\right) -\alpha
_{l}\sum_{k}\mathrm{e}^{\mathrm{i}\theta _{k}}+o\left( 1\right)
\end{equation*}%
as $\left\vert t\right\vert \rightarrow \infty .$ We see that the leading
part of the phase is $-\alpha _{l}\ln t$, which implies the asymptotic $%
\approx \mathrm{const}\cdot t^{-\alpha _{l}}$ (as expected). But the
`oscillating' part, as function of the angles $\theta _{j}$ and $\theta
_{k}, $ is quite irregular (not large). So, the corresponding integral is
not easy and we do not have a formula for the constant.

Consider the case of a critical point corresponding to $\rho =\rho ^{\left(
m\right) }.$ We put $a_{j}=a_{j}^{\left( m\right) }\mathrm{e}^{\mathrm{i}%
\theta _{j}}\approx \frac{1}{\eta }\mathrm{e}^{\mathrm{i}\theta _{j}}$ and $%
a_{k}=a_{k}^{\left( m\right) }\mathrm{e}^{\mathrm{i}\theta _{k}}\approx
\zeta ^{m-p}\eta ^{p\kappa }\mathrm{e}^{\mathrm{i}\theta _{k}}$ and get%
\begin{equation*}
\phi =-\sum_{j\leq p}\alpha _{j}\ln \left( 1-\mathrm{e}^{\mathrm{i}\theta
_{j}}\right) +\zeta ^{m-p}\eta ^{\left( q+1\right) \kappa }t^{\kappa
}\sum_{k>p}\mathrm{e}^{\mathrm{i}\theta _{k}}.
\end{equation*}%
Here we have the dominant part with $t^{\kappa }$ (before the second sum),
which can be treated via the stationary phase formula. But there remains the
finite part leading to an integral of the form $\int \prod \left( 1-\mathrm{e%
}^{\mathrm{i}\theta _{j}}\right) ^{-\alpha _{j}}\mathrm{d}\theta _{j}.$%

\end{document}